\documentclass[a4paper,12pt]{amsart}
\usepackage{latexsym,amsmath,amssymb,amsthm}
\newtheorem{thm}{Theorem}[section]
\newtheorem{lemma}[thm]{Lemma}
\newtheorem{cor}[thm]{Corollary}
\newtheorem{rem}[thm]{Remark} 

\newtheorem{define}[thm]{Definition}
\newtheorem{prop}[thm]{Proposition}

\numberwithin{equation}{section}

\def\o{\omega}

\def\R{\mathbb R}

\input{amssym.def}

\title[A general regularity theory for weak MCF]{A general regularity theory \\ for weak mean curvature flow}
\author[K. Kasai]{Kota Kasai}

\address{Department of Mathematics, Hokkaido University, Sapporo 060-0810 Japan.}
\email{kasai@math.sci.hokudai.ac.jp}

\author[Y. Tonegawa]{Yoshihiro Tonegawa}

\address{Department of Mathematics, Hokkaido University, Sapporo 060-0810 Japan.}
\email{tonegawa@math.sci.hokudai.ac.jp}

\date{}


\thanks{Y.Tonegawa is partially supported by JSPS Grant-in-aid for scientific research (B) $\#$21340033 and (S) $\#$21224001.
He also thanks Neshan Wickramasekera for various useful comments on the preliminary version of the paper.}

\begin{document}
\setlength\parskip{5pt}
\begin{abstract}
We give a new proof of Brakke's partial regularity theorem up to
$C^{1,\varsigma}$ for
weak varifold solutions of mean curvature flow by utilizing
parabolic monotonicity formula,
parabolic Lipschitz approximation and blow-up technique. 
The new proof extends to a general flow whose
velocity is the sum of the mean curvature and any given background
flow field in a dimensionally sharp integrability class. 
It is a natural parabolic generalization of 
Allard's regularity theorem in the sense that the 
special time-independent case reduces to Allard's theorem. 
\end{abstract}

\maketitle
\makeatletter
\@addtoreset{equation}{section}                  
\renewcommand{\theequation}{\thesection.\@arabic\c@equation}
\makeatother

\section{Introduction}
A family $\{M_t\}_{t\geq 0}$of $k$-dimensional surfaces in ${\mathbb R}^n$ 
is called the
mean curvature flow (hereafter abbreviated MCF) 
if the velocity of $M_t$ is equal to its
mean curvature at each point and time. 
As one of the most fundamental geometric evolution problems, 
the MCF has been the subject of intensive research since
1980's.
The earliest study of MCF goes back to the seminal work of Brakke 
\cite{Brakke} who used the notion of varifold \cite{Allard} in
geometric measure theory to show
the existence of weak solutions for general initial surfaces. 
More precisely, given any $k$-dimensional integral varifold $V_0$ with
some mild finiteness assumptions, he showed the existence of a family of
varifolds $\{V_t\}_{t\geq 0}$ 
each of which satisfies the MCF equation in a distributional 
sense for all $t\geq 0$. For a.e$.$ time, Brakke additionally proved that
$V_t$ is integral. Assuming that the density function is equal to 1 
almost everywhere, Brakke also claimed that the weak MCF 
varifold solutions
are supported by smooth $k$-dimensional
surfaces almost everywhere. The hypothesis is called the `unit density
hypothesis' and it is a natural assumption since even the
time-independent case of Allard's regularity theory has essentially the same hypothesis.
Brakke's regularity 
theory introduced remarkably ingenious tools such as 
the `clearing-out lemma' and the `popping soap film lemma', and it has a 
significant influence on the analysis of the related singular perturbation
problems such as the Allen-Cahn equation \cite{Ilmanen} and the parabolic
Ginzburg-Landau equation \cite{Ambrosio-Soner,BOS,Jerrard-Soner,Lin}. 
However, the detail of Brakke's regularity theory is 
technically involved and a clear accessible proof is desired due to its 
importance.  

The aim of the present paper is twofold. The first aim is to give a new and
self-contained proof of Brakke's
partial regularity theorem up to $C^{1,\varsigma}$. Brakke's proof relies on a long 
chain of graphical approximations which is complicated and hard to
follow (see \cite[6.9, `Flattening out']{Brakke}). Utilizing the `$L^2$-$L^{\infty}$
lemma' (which we explain below) we replace this part with 
Allard-like Lipschitz approximations. 
The second aim is to generalize the result so that the
velocity of motion may have an additional transport term which belongs to a certain 
integrability class. The existence of such flow has been studied by 
Liu-Sato-Tonegawa \cite{LST} and it motivated the authors to investigate the 
generalization of Brakke's theorem. Also, such extra term may naturally arise
when one considers the MCF in Riemannian manifolds via Nash's imbedding
theorem. From a view
point of regularity theory,
the result of the present paper 
reduces essentially to Allard's regularity theorem
\cite{Allard} in the special case of time-independent case. 
We note that simple modifications of Brakke's original proof
do not seem to yield the theorem of the present paper 
when one puts the general transport term. Building upon the
present work and assuming that the transport term $u$ (see the
assumption) is $\alpha$-H\"{o}lder continuous, the second author established $C^{2,\alpha}$ regularity theorem
in the forthcoming paper \cite{regTonegawa}. In particular
\cite{regTonegawa} proves that the unit density Brakke's MCF without 
transport term (i.e. $u=0$) is for a.e$.$ time and a.e$.$ everywhere
$C^{\infty}$, which was originally claimed in \cite{Brakke} (see Section \ref{conrem1}
for a further comment). 

The major difference of our proof compared to \cite{Brakke} is the use of
the `$L^2$-$L^{\infty}$ lemma' of Section \ref{moform}. 
It shows that the smallness
of measurement of `height' of MCF in the $L^2$ sense (hereafter called $L^2$-height)
guarantees that the whole support of MCF lies in a narrow
region close to a $k$-dimensional plane. The proof of this fact utilizes
an analogy of the well-known parabolic monotonicity
formula due to Huisken \cite{Huisken} (and see \cite{Chen-Struwe,Struwe} for
similar formula in case of harmonic heat flow).
Note that Brakke's work preceded the discovery 
of parabolic monotonicity formula and only the `elliptic' monotonicity
formula was available at the time.
The $L^2$-$L^{\infty}$ lemma may be 
considered as a robust version of Brakke's clearing-out lemma in 
the sense that the former may accommodate general transport terms 
while the latter apparently has some limitations doing so. With this new
input the proof may be outlined as follows. We make
a full use of the so-called popping soap film lemma with some modifications from \cite[6.6]{Brakke}, 
which gives a control of the 
space-time
$L^2$ norm of the mean curvature in terms of the smallness of $L^2$-height (corresponding 
to Section \ref{enees}). 
Then the
rest of the proof proceeds more or less like Allard's regularity proof with
parabolic modifications. Namely we approximate the support of moving
varifolds by a Lipschitz
graph (with respect to the parabolic metric) utilizing parabolic
monotonicity formula. Then through a contradiction 
argument, we employ the blow-up technique. 
While the blow-up limit is a harmonic
function in Allard's case, it is a solution 
of the heat equation in the present case. 
This gives a decay estimate of $L^2$-height with respect to a slightly 
tilted $k$-dimensional plane in a smaller scale. 
Iteration argument gives H\"{o}lder estimate
of the spacial gradient of the graph just like Allard's case. Though the 
proof of the present paper is lengthy and
technically involved, we expect that the researchers familiar with regularity
theory would find the details more tractable and natural than
Brakke's original one. 

There have been a large amount of research on the MCF
and we may mention only a small fraction of them related to 
Brakke's varifold solutions. For existence of generalized solutions
we mention \cite{AWT,CGG,ES1,ES2,LS}. In particular 
Evans-Spruck \cite{ES2} proved that
the almost every level set of viscosity solution of MCF in the context of
level set method is MCF in Brakke's sense. These level sets naturally
satisfy the aforementioned unit density hypothesis, hence Brakke's 
theorem applies. White discovered a simple proof of
local regularity theorem \cite{White1} up to and including the time 
at which singularities first occur in any classical MCF. He showed the
local regularity when the Gaussian density is close to 1. This result
is sufficient for many interesting cases (see the introduction of 
\cite{White1}), but nevertheless, does not
replace Brakke's theorem in full. There have been a significant advance
of knowledge in the case of mean convex hypersurface, where the salient features
of singularities of MCF have been established (\cite{Metzger,White2,White3}). 
The good reference on the subject is \cite{Ecker}, where one can find
many useful estimates for Brakke's flow written in smooth setting.

The organization of the paper is as follows. Section 2 lists notations and recalls
some well-known results from geometric measure theory. Section 3 contains the 
assumptions and main result. The content of Section 4, the Lipschitz approximation
of varifolds for fixed time, is used in the subsequent Section 5, which deals with 
the local energy estimate, or tilt-excess estimate in terms of $L^2$-height if one looks
for an analogy with Allard's regularity theory. Section 6 is independent of the previous
two sections, and establishes monotonicity-type formulae and $L^2$-$L^{\infty}$ estimates. 
Section 7 is independent and gives construction of parabolic Lipschitz graph which
approximates moving varifold with good error bounds.
Based on the estimates and constructions obtained
through Section 5-7, Section 8 shows the decay estimates of $L^2$-height via blow-up
techniques. Theorem \ref{regth} is the main local regularity theorem. Section 9 shows
that Theorem \ref{regth} is applicable for a.e$.$ time and concludes the proof of main
partial regularity theorem. Section 10 lists a few concluding remarks and 
Section 11 collects basic results from \cite{Allard} and \cite{Brakke} for the convenience of
the reader. 

\section{Preliminaries}
\subsection{Basic notations}
Throughout this paper, $k$ and $n$ will be positive integers with $0<k<n$.
By regarding $n$ and $k$ as absolute constants, dependence of constants on $k$ and $n$ may not be stated explicitly. 
Let ${\mathbb N}$ be the natural number and ${\mathbb R}^+:=\{x\geq 0\}$. 
For $0<r<\infty$ and $a\in \R^n$ (or $\R^k$) let
\begin{equation*}
B_r(a):=\{x\in {\mathbb R}^n\, :\, |x-a|<r\},\, \, B_r^k(a):=\{x\in 
{\mathbb R}^k\, :\, |x-a|<r\}
\end{equation*}
and when $a=0$ let $B_r:=B_r(0)$ and $B^k_r:=B^k_r(0)$. We also define
$\tau(r)(x):=rx$. For $s\in {\mathbb R}$ 
define a parabolic cylinder
\begin{equation*}
P_r(a,s):=\{(x,t)\in {\mathbb R}^n\times{\mathbb R}\, :\, |x-a|<r,\,|t-s|<r^2\}.
\end{equation*}
Note that it is not an often-used parabolic cylinder with $t$ located at the top end.
Let $\overline{A}$
be the closure of $A\subset{\mathbb R}^n$. For $a\in \overline{A}$ let
\begin{equation*}
{\rm Tan}\, (A,a)
\end{equation*}
be the closed cone whose intersection with $\{x\in {\mathbb R}^n\, :\, |x|=1\}$ is given by
\begin{equation}
\cap_{0<r<\infty}\overline{\{(x-a)/|x-a|\,:\, x\in A\cap B_r(a)\setminus\{a\}\}}.
\label{deftan}
\end{equation}
We denote
by ${\mathcal L}^n$ the Lebesgue measure on ${\mathbb R}^n$ and by
${\mathcal H}^k$ the $k$-dimensional Hausdorff measure on ${\mathbb R}^n$.
The restriction of ${\mathcal H}^k$ to a set $A$ is denoted by ${\mathcal H}^k
\lfloor_{A}$. 
We let $\omega_k:={\mathcal L}^k(B_1^k)$.
For an open subset $U\subset {\mathbb R}^n$ let $C_c(U)$ be the set of all compactly 
supported continuous functions on $U$ and let $C_c(U;{\mathbb R}^n)$ be 
the set of all compactly supported, continuous vector fields. The upper subscript
of $C_c^l(U)$ and $C_c^l(U;{\mathbb R}^n)$ indicates 
continuous $l$-th order differentiability. For $g\in C^1(U;{\mathbb R}^n)$,
we regard $\nabla g(x)$ as an element of ${\rm Hom}({\mathbb R}^n,{\mathbb R}^n)$.
Similarly for $g\in C^1(U)$, we regard the Hessian matrix $\nabla^2 g(x)$ as an 
element of ${\rm Hom}({\mathbb R}^n,{\mathbb R}^n)$. $\nabla$ always indicates
differentiation with respect to space variables, and not with respect to time variable. 

For any Radon measure $\mu$ on 
${\mathbb R}^n$ and $\phi\in C_c({\mathbb R}^n)$ we often write 
$\mu(\phi)$ for $\int_{{\mathbb R}^n}\phi\, d\mu$. Let ${\rm spt}\, \mu$ be
the support of $\mu$, i.e., $x\in {\rm spt}\, \mu$ if $\forall r>0,\,
\mu(B_r(x))>0$. 
Let $\Theta^k(\mu,x)$ be the $k$-dimensional density of $\mu$ at $x$, i.e., 
$\lim_{r\rightarrow 0}\mu(B_r(x))/(\o_k r^k)$, when the limit exists. 
For $\mu$ a.e$.$ defined 
function $u$, and $1\leq p\leq \infty$, $u\in L^p(\mu)$ means $\left(\int |u|^p\, d\mu\right)^{1/p}<\infty$. 

For $-\infty<t<s<\infty$ and $x,\, y\in {\mathbb R}^n$, define
\begin{equation}
\rho_{(y,s)}(x,t):=\frac{1}{(4\pi(s-t))^{k/2}}\exp\left(-
\frac{|x-y|^2}{4(s-t)}\right).
\label{hkneldef}
\end{equation}
$\rho_{(y,s)}$ is the $k$-dimensional backward heat kernel which is often
used throughout this paper.

\subsection{The Grassmann manifold and varifolds}
Let ${\bf G}(n,k)$ be the space of $k$-dimensional subspaces of ${\mathbb R}^n$
and let ${\bf A}(n,k)$ be the space of $k$-dimensional affine planes of ${\mathbb R}^n$.
For $S\in {\bf G}(n,k)$, we identify $S$ with the corresponding orthogonal 
projection of ${\mathbb R}^n$ onto $S$. Let $S^{\perp}\in {\bf G}(n,n-k)$ be the orthogonal
complement of $S$. For two  elements $A$
and $B$ of ${\rm Hom}\, ({\mathbb R}^n,{\mathbb R}^n)$, define a scalar product $A\cdot B:={\rm trace}\, (A^*\circ B)$ where $A^*$ is the transpose of $A$ and $\circ$ indicates the usual composition. The identity of ${\rm Hom}\, ({\mathbb R}^n,{\mathbb R}^n)$ is denoted by $I$. Let $a\otimes
b\in {\rm Hom}\, ({\mathbb R}^n,{\mathbb R}^n)$ be the tensor product of $a,\, b\in {\mathbb R}^n$. For $A\in {\rm Hom}\, ({\mathbb R}^n,{\mathbb R}^n)$ define
\begin{equation*}
|A|:=\sqrt{A\cdot A},\hspace{1cm}
\|A\|:=\sup\{|A(x)|\, :\, x\in {\mathbb R}^n,\, |x|=1\}.
\end{equation*}
For $T\in {\bf G}(n,k)$, $a\in {\mathbb R}^n$ and $0<r<\infty$ we define the cylinder
\begin{equation*}
C(T,a,r):=\{x\in {\mathbb R}^n\, :\, |T(x-a)|<r\},\,\, C(T,r):=C(T,0,r).
\end{equation*}
We recall some notions related to varifold and refer to \cite{Allard,Simon} for more
details. For any open set $U\subset {\mathbb R}^n$, define $G_k(U):=U\times {\bf G}(n,k)$. 
A general $k$-varifold in $U$ is a Radon measure on 
$G_k(U)$. Set of all general $k$-varifolds in $U$ is denoted by ${\bf V}_k(U)$. 
For $V\in {\bf V}_k(U)$, let $\|V\|$ be the mass measure of $V$, 
namely, 
\begin{equation*}
\|V\|(\phi):=\int_{G_k(U)}\phi(x)\, dV(x,S),\,\,\,\forall \phi\in C_c(U).
\end{equation*}
For
proper map $f\in C^1({\mathbb R}^n;{\mathbb R}^n)$ define
$f_{\#}V$ as the push-forward of varifold (see \cite{Allard} for the definition).
Given any ${\mathcal H}^k$ measurable 
countably $k$-rectifiable set $M\subset U$ with
locally finite ${\mathcal H}^k$ measure, there
is a natural $k$-varifold $|M|\in {\bf V}_k(U)$ defined by
\begin{equation*}
|M|(\phi):=\int_{M}\phi(x,{\rm Tan}_x M)\, d{\mathcal H}^k(x),\,\,\,\forall \phi\in C_c(G_k(U)),
\end{equation*}
where ${\rm Tan}_x M\in {\bf G}(n,k)$ is the approximate tangent space which exists
${\mathcal H}^k$ a.e$.$ on $M$. In this case, $\||M|\|={\mathcal H}^k\lfloor_M$.
We say $V\in {\bf V}_k(U)$ is integral if 
\begin{equation*}
V(\phi)=\int_{M}\phi(x,{\rm Tan}_x M)\theta(x)\, d{\mathcal H}^k(x),\,\,\,\forall \phi\in C_c(G_k(U)),
\end{equation*}
with some ${\mathcal H}^k$ measurable 
countably $k$-rectifiable set $M\subset U$ and ${\mathcal H}^k$
a.e$.$ integer-valued integrable 
function $\theta$ defined on $M$. Note that for such varifold, $\Theta^k(\|V\|,x)
=\theta(x)\in {\mathbb N}$, ${\mathcal H}^k$ a.e$.$ on $M$. Set of all 
integral $k$-varifolds in $U$ is denoted by ${\bf IV}_k(U)$.
We say $V$ is a unit density 
$k$-varifold if $V$ is integral and $\theta=1$ a.e$.$ on $M$, that is, $V=|M|$. 
When $V$ is integral, we often
write
$\int_{U}(g(x))^{\perp}\, d\|V\|(x)$ for $\int_{G_k(U)}
S^{\perp}(g(x))\, dV(x,S)$, for example, since there should be no ambiguity. 

\subsection{First variation and generalized mean curvature}
For $V\in {\bf V}_k(U)$ let $\delta V$ be the first variation of $V$, namely,
\begin{equation*}
\delta V(g):=\int_{G_k(U)}\nabla g(x)\cdot S\, dV(x,S)
\end{equation*}
for $g\in C_c^1(U;{\mathbb R}^n)$. 
Let $\|\delta V\|$ be the total variation when it exists, and if 
$\|\delta V\|$ is absolutely continuous with respect to $\|V\|$, we have
for some $\|V\|$ measurable vector field $h(V,\cdot)$
\begin{equation}
\delta V(g)=-\int_{U}g(x)\cdot h(V,x)\, d\|V\|(x).
\label{fvf}
\end{equation}
The vector field $h(V,\cdot)$ is called the generalized mean curvature of $V$. 
We say $V$ is stationary if $h(V,\cdot)=0$, $\|V\|$ a.e$.$ in $U$, or equivalently, 
$\delta V(g)=0$ for all $g\in C^1_c(U;{\mathbb R}^n)$. 
For any $V\in {\bf IV}_k(U)$ with integrable $h(V,\cdot)$, Brakke's 
perpendicularity theorem of generalized mean curvature
\cite[Chapter 5]{Brakke} says that we have
\begin{equation}
\int_{U} (g(x))^{\perp}\cdot h(V,x)\, d\|V\|(x)=\int_{U} g(x)\cdot h(V,x)\, d\|V\|(x)
\label{fvf2}
\end{equation}
for all $g\in C_c(U;{\mathbb R}^n)$.

\subsection{The right-hand side of MCF equation}
For any $V\in {\bf V}_k(U)$, $u\in L^2(\|V\|)$
and $\phi\in C^1_c(U;{\mathbb R}^+)$, define
\begin{equation}
{\mathcal B}(V,u,\phi):=\int_{U} (-\phi(x)h(V,x)+\nabla\phi(x))\cdot(h(V,x)+
(u(x))^{\perp})\, d\|V\|(x)
\label{Bdef}
\end{equation}
when $V\in {\bf IV}_k(U)$, $\|\delta V\|$ is locally finite and absolutely continuous
with respect to $\|V\|$, and $h(V,\cdot)\in L^2(\|V\|)$. Otherwise we define
${\mathcal B}(V,u,\phi)=-\infty$. 
Formally, if a family of smooth $k$-dimensional surfaces $\{M_t\}$ moves
by the velocity equal to the mean curvature plus smooth $u$, then, one
can check that $V_t=|M_t|$ satisfies
\begin{equation}
\frac{d}{dt}\|V_t\|(\phi)\leq {\mathcal B}(V_t,u,\phi),\,\,\,\forall\phi\in C_c^1(U;{\mathbb R}^+).
\label{formal}
\end{equation}
In fact, \eqref{formal}
holds with equality. Conversely, if \eqref{formal} is satisfied, then one can 
prove that the velocity is equal to the mean curvature plus $u$. The inequality
in \eqref{formal} allows the sudden loss of surface and it is the source of 
general non-uniqueness of Brakke's formulation. 

\section{Assumptions and main result}
\subsection{Assumptions}
For an open set $U\subset {\mathbb R}^n$ and $0<\Lambda\leq \infty$ 
suppose that we have a 
family of $k$-varifolds $\{V_t\}_{0\leq 
t< \Lambda}$ and a family of $n$-vector valued functions
$\{u(\cdot,t)\}_{0\leq t< \Lambda}$ both on $U$ satisfying the followings.
\newline
{\bf (A1)} For a.e$.$ $t\in [0,\Lambda)$, $V_t$ is a unit density $k$-varifold.
\newline
{\bf (A2)} There exists $1\leq E_1<\infty$ such that 
\begin{equation}
\|V_t\|(B_r(x))
\leq \omega_k r^k E_1,\,\,\,\forall B_r(x)\subset U,\,\, \forall t\in [0,\Lambda).
\label{dbd}
\end{equation}
\newline
{\bf (A3)} Suppose $2\leq p<\infty$ and $2<q<\infty$ satisfy
\begin{equation}
\varsigma:=1-\frac{k}{p}-\frac{2}{q}>0
\label{expcond}
\end{equation}
and assume that $u$ satisfies
\begin{equation}
\|u\|_{L^{p,q}(U\times(0,\Lambda))}:=
\left(\int_0^{\Lambda}\left(\int_{U}|u(x,t)|^p\, d\|V_t\|(x)\right)^{\frac{q}{p}}\, dt\right)^{\frac{1}{q}}<\infty.
\label{ubd}
\end{equation}
\newline
{\bf (A4)}
For all $\phi\in C^1(U\times[0,\Lambda);{\mathbb R}^+)$ with $\phi(\cdot,t)\in C^1_c(U)$ and $0\leq t_1<t_2<\Lambda$,
\begin{equation}
\|V_{t_2}\|(\phi(\cdot,t_2))-\|V_{t_1}\|(\phi(\cdot,t_1)) 
\leq \int_{t_1}^{t_2}
{\mathcal B}(V_t,u(\cdot,t),\phi(\cdot,t))\, dt+
\int_{t_1}^{t_2}\int_{U}\frac{\partial\phi}{\partial t}(\cdot,t)\, d\|V_t\| dt\label{maineq}
\end{equation}
holds. 
\label{assumptions}

As an immediate consequence,  for $\phi\in C^2_c(U;{\mathbb R}^+)$, we have from \eqref{maineq} and Cauchy-Schwarz
inequality
\begin{equation}
\begin{array}{ll}
\|V_{t_2}\|(\phi)-\|V_{t_1}\|(\phi)& \leq \int_{t_1}^{t_2}\int_U -\phi|h|^2+|\nabla\phi| |h|
+\phi|h||u|+|\nabla\phi||u|\, d\|V_t\|dt \\
& \leq \int_{t_1}^{t_2}\int_U -\frac{|h|^2\phi}{2}+\frac{|\nabla\phi|^2}{\phi}+|u|^2\phi+|\nabla\phi||u|\, d\|V_t\|dt \\
&\leq \int_{t_1}^{t_2}\int_U -\frac{|h|^2\phi}{2}\, d\|V_t\|dt+c(\|\phi\|_{C^2},E_1,\|u\|_{L^{p,q}},{\rm spt}\,\phi).
\end{array}
\label{hhhp}
\end{equation}
 Thus we have from \eqref{hhhp} that
 \begin{equation}
 h(V_t,\cdot)\in L_{loc}^2(\|V_t\|)
 \label{hhh}
 \end{equation}
 for a.e$.$ $t\in (0,\Lambda)$.
\subsection{Remarks on the assumptions}
Before stating our main theorem, we motivate our assumptions. 

The setting that Brakke proved the local regularity theorem 
corresponds to $u= 0$ case. We note that the weak 
solution he constructed in \cite{Brakke} satisfies (A2)-(A4)
locally in space and time. In this special case, we may relax the assumption (A2)
by, for example, 
\begin{equation*}
\sup_{0\leq t<\Lambda}\|V_t\|(U)\leq C,
\end{equation*}
or 
\begin{equation*}
U=B_R,\,\, \|V_0\|(U)<\infty \,\,\mbox{and}\,\,\, {\rm spt}\,\|V_0\|\subset \subset U.
\end{equation*}
The assumption (A2) (restricted for $U'\times(\delta,\Lambda)$ with $\forall U'\subset\subset U$ and $\forall \delta>0$) follows 
from the parabolic monotonicity formula of Section \ref{monosubsec}, 
and \eqref{maineq} is the integral form of Brakke's formulation,
except that we need time varying test functions.
See the discussion in \cite[3.5]{Brakke} for the derivation.   
The varifolds are also integral
for a.e$.$ time if $V_0$ is assumed to be integral. 
With parabolic monotonicity formula it is not difficult to prove that there is an initial
short time interval during which the varifolds remain unit density
if the initial varifold $V_0$ has (A2) with $E_1$ close to 1. 

As we noted in Section 1, in \cite{ES2} Evans-Spruck proved
that almost all level set of MCF viscosity solution is a unit density
MCF in the sense of Brakke. Thus almost all level set satisfies (A1)-(A4) with $u=0$.

Ilmanen \cite{Ilmanen} proved that the singular perturbation limit of 
the Allen-Cahn equation is a rectifiable Brakke's MCF under mild conditions on
initial data. In addition the second author \cite{Tonegawa} proved that the limit varifolds are 
integral modulo division by a surface energy constant for a.e$.$ time. 
Thus the limit varifolds arising from the Allen-Cahn equation, unless
there are portions where higher ($\geq 2$) multiplicities occur, satisfy 
(A1)-(A4). We mention that a simple proof of \cite{Ilmanen,Tonegawa}
was obtained by Sato \cite{Sato} utilizing the result of R\"{o}ger-Sch\"{a}tzle \cite{RS}
(see also \cite{Nagase}) for $k=n-1$
and $n=2,\,3$. 

As for $u\neq 0$ case, we note first
that \eqref{expcond} is dimensionally a sharp condition. A simple
dimension analysis of natural re-scaling $\tilde{x}=\lambda^{-1} x$, $\tilde{t}
=\lambda^{-2} t$ and $\tilde{u}=\lambda u$ for $\lambda>0$ shows
\begin{equation*}
\|\tilde{u}\|_{L^{p,q}}=\lambda^{1-\frac{k}{p}
-\frac{2}{q}}\|u\|_{L^{p,q}}.
\end{equation*}
It is a rule of thumb that \eqref{expcond} is essential to
successfully obtain a local regularity theorem. We also note that
if $V_t=V_0$ for all $t>0$, then $q$ may be regarded as $+\infty$
and \eqref{expcond} requires $p>k$. Time-independence means the velocity is 0. Then 
we have $0=h(V,\cdot)+u^{\perp}$ (see Section 10)
and \eqref{ubd}
gives $h(V,\cdot)\in L^p(\|V\|)$ with $p>k$. It is the 
condition for the application of Allard's regularity theorem. 

The existence
of flow, where the velocity has additional transport term $u$, has been investigated in \cite{LST} for $k=n-1$, $n=2,3$.
There, given an arbitrary $u$ in $L^p_{\rm loc}([0,\infty);W^{1,p}({\mathbb T}^n))$ with
$p>\frac{n+2}{2}$, where ${\mathbb T}^n$ is the $n$-dimensional torus and
$W^{1,p}({\mathbb T}^n)=\{u\in L^p({\mathbb T}^n)\, :\, \nabla u\in L^p({\mathbb T}^n)\}$, and 
$C^1$ initial ${\rm spt}\,\|V_0\|$, it was proved
that there exists a family $\{V_t\}_{t\geq 0}$ of $(n-1)$-varifolds satisfying (A2)-(A4) and a.e$.$ time integrality. 
$u(\cdot,t)$ in this case is defined as a trace on the support of $\|V_t\|$.
The fulfillment of condition \eqref{ubd} is not stated explicitly there, 
but can be verified
using \cite[Theorem 2.1]{LST}. The unit density hypothesis (A1) is 
satisfied for at least short time \cite[Theorem 2.2(c)]{LST}. Thus we can apply the
present regularity theory to this initial short time. Unlike $u=0$ case, we are not
aware of any short time existence result for classical solutions under low regularity 
assumptions on $u$.  
The extension of \cite{LST} for $n>3$ and $k=n-1$ which covers the whole 
range of $p$ and $q$ satisfying \eqref{expcond} is currently under 
investigation with similar approach. We are unaware of any general
method in case $k<n-1$, most relevant being the approximation scheme via
the parabolic Ginzburg-Landau equation \cite{BOS} where $u=0$, $k=n-2$ and $n\geq 3$. 
Note that \cite{BOS} proved that
the limit varifold $V_t$ corresponding to vorticity concentration 
is rectifiable $k$-varifold moving by MCF in Brakke's formulation. 
One the other hand, one apparently does not exclude the possibility of non-integral varifold
due to the non-trivial diffuse part of the limit measure. 
\subsection{Main result}
\begin{define}
A point $x\in U\cap{\rm spt}\, \|V_t\|$ is said to be a $C^{1,\varsigma}$ regular point if there
exists some open neighborhood $O$ in ${\mathbb R}^{n+1}$ containing $(x,t)$ 
such that $O\cap \cup_{0<s<\Lambda}({\rm spt}\, \|V_s\|\times\{s\})$
is an embedded $(k+1)$-dimensional manifold represented as the image of 
$f:O'(:=B^k_R\times (t-R^2,t+R^2))\rightarrow O\subset{\mathbb R}^{n+1}$ with 
\begin{equation*}
\sup_{(y_j,s_j)\in O',\,j=1,2}\frac{\|\nabla f(y_1,s_1)-\nabla f(y_2,s_2)\|}{\max\{|y_1-y_2|^{\varsigma},
|s_1-s_2|^{\varsigma/2}\}}<\infty,\,\,\, \sup_{(y,s_j)\in O',\, j=1,2}\frac{|f(y,s_1)-f(y,s_2)|}{|s_1-s_2|^{(1+\varsigma)/2}}<\infty.
\end{equation*}
Here $\nabla$ is
the differentiation with respect to $y$-variables.
\end{define}
Here is the main partial regularity theorem.
\begin{thm} Under the assumptions (A1)-(A4), for a.e$.$ $t\in (0,\Lambda)$, there exists a (possibly empty) closed set $G_t\subset{\rm spt}\, \|V_t\|$ with 
${\mathcal H}^k(G_t)=0$ such that ${\rm spt}\, \|V_t\|\setminus G_t$ is a set of $C^{1,\varsigma}$ regular points. 
\label{gp}
\end{thm}
Theorem \ref{gp} is the consequence of Theorem \ref{regth} and Theorem 
\ref{prprop} which are the main local regularity theorems. For 
more detailed local estimates and assumptions, see the statement of Theorem \ref{regth}. 
\section{Lipschitz approximation for fixed time}
In this section, the main result is Proposition \ref{lipapxmn} which 
gives a graphical Lipschitz approximation of varifold, and which is one of the
essential ingredients for Theorem \ref{poptheorem} of the next section.  The results of 
this section do not involve time variable. The content of this section 
corresponds to \cite[5.3,\, 5.4]{Brakke}. To avoid unnecessary technical 
complications, we avoid general multiple-sheet situations discussed in \cite{Brakke}
which are not needed for the purpose of this paper.  
\begin{lemma}
Corresponding to each $0<\lambda<1$ there exists $\gamma>0$
with the following property.
For $V\in {\bf IV}_k(B_R)$ assume
\begin{equation}
\Theta^k(\|V\|,0)=1,
\label{sinmono1}
\end{equation}
\begin{equation}
r\|\delta V\|(B_r)\leq \gamma \|V\|(B_r),\hspace{.3cm}0<\forall r\leq R.
\label{sinmono3}
\end{equation}
Then we have
\begin{equation}
\| V\|(B_r)\geq \lambda\omega_k r^k, \hspace{.3cm}0<\forall r\leq R.
\label{specsinm1}
\end{equation}
\label{lemsinmo}
\end{lemma}
\begin{rem}
If the condition \eqref{sinmono3} is replaced by $\|\delta V\|(B_r)\leq\gamma\|V\|(B_r)$,
then the well-known monotonicity formula \cite[5.1]{Allard} says $r^{-k}\|V\|(B_r)\exp(\gamma r)$
is monotone increasing. The reason for using the weaker condition \eqref{sinmono3} is to obtain
a better error estimate for Lipschitz approximation in Proposition \ref{lipapxmn}. 
See Remark \ref{imprem} for further comment.
\end{rem}
{\it Proof}. 
Without loss of generality we may assume $R=1$. Assume that the 
conclusion were false.
Then for each $m\in {\mathbb N}$ there exist $V_m\in {\bf IV}_k(B_1)$ such that 
\eqref{sinmono1} and \eqref{sinmono3} hold for $V=V_m$ and $\gamma=
1/m$, but \eqref{specsinm1} does not hold. Let $R_m$ be the supremum of $r>0$ such
that $\|V_m\|(B_s)\geq \lambda\omega_k s^k$ holds for all $0<s\leq r$. By \eqref{sinmono1},
we have $R_m>0$. Since we are assuming the negation of \eqref{specsinm1}, we have
$R_m<1$ and $\|V_m\|(B_{R_m})=\lambda\omega_k R_m^k$. 
Let $\tilde{V}_m=(\tau(R_m^{-1}))_{\#}V_m$. 
Then we have by \eqref{sinmono3} and by the definition of $R_m$ 
\begin{equation}
\| \delta \tilde{V}_m\|(B_1)=\frac{1}{R_m^{k-1}}\|\delta V_m\|(B_{R_m})
\leq \frac{1}{m R_m^k}\|V_m\|(B_{R_m}) = \frac{\omega_k \lambda}{m},
\label{sinmono7}
\end{equation}
\begin{equation}
\|\tilde{V}_m\|(B_1)=\frac{1}{R_m^k}\|V_m\|(B_{R_m})= \omega_k \lambda,
\label{sinmono8}
\end{equation}
\begin{equation}
\|\tilde{V}_m\|(B_s)=\frac{1}{R_m^k}\|V_m\|(B_{sR_m})\geq \lambda \omega_k s^k, 
\hspace{.3cm}0<\forall s\leq 1.
\label{sinmono9}
\end{equation}
By \eqref{sinmono7} and \eqref{sinmono8} and the compactness theorem of 
integral varifolds (\cite[6.4]{Allard}), there exists a subsequence $\{\tilde{V}_{m_j}\}$ 
and $V\in {\bf IV}_k(B_1)$ such that $\tilde{V}_{m_j}\rightarrow V$ in the sense
of varifolds. By \eqref{sinmono7}, $V$ is stationary, and by \eqref{sinmono9},
$0\in {\rm spt}\, \|V\|$. From \eqref{sinmono8} and \eqref{sinmono9}, we also
have
\begin{equation}
\|V\|(B_1)=\lambda\omega_k.
\label{sinmono11}
\end{equation}
By the upper-semicontinuity of density function 
and the integrality of $V$, we have $\Theta^k(\|V\|,0)\geq 1$. 
By the monotonicity formula of stationary varifold, we have $\|V\|(B_1)
\geq \omega_k$, which contradicts with \eqref{sinmono11} since $\lambda<1$.
This proves \eqref{specsinm1}. 
\hfill{$\Box$}
\begin{lemma} Corresponding to each $0<\lambda<1$ and $1\leq E_1<\infty$,
there exists $\gamma>0$ with the following property. For $V\in {\bf IV}_k(B_R)$
assume \eqref{sinmono1}, \eqref{sinmono3},
\begin{equation}
\|V\|(B_R)\leq  \omega_k R^k E_1,
\label{sinmono4}
\end{equation}
\begin{equation}
\int_{G_k(B_R)}\|S-T\|\, dV(x,S)\leq \gamma \|V\|(B_R).
\label{sinmono2}
\end{equation}
Then for any $b\in T$ with $|b|\leq R/2$, 
\begin{equation}
\|V\|(B_r(b))\geq \lambda \omega_k r^k,\hspace{.3cm}R/10\leq \forall r\leq R/2.
\label{sinmono5}
\end{equation}
\label{lemsinmo2}
\end{lemma}
{\it Proof}.
Assume $R=1$ without loss of generality.
Assume that the conclusion were false. Then for each $m\in {\mathbb N}$
we have a sequence $V_m\in {\bf IV}_k(B_1)$, $b_m\in T$ with $|b_m|\leq 1/2$
and $R_m$ with $1/10\leq R_m\leq 1/2$ such that \eqref{sinmono1}, \eqref{sinmono3},
\eqref{sinmono4} and 
\eqref{sinmono2}
hold for $\gamma=1/m$ and $V=V_m$ but \eqref{sinmono5} fails, that is,
\begin{equation}
\|V_m\|(B_{R_m}(b_m))<\lambda \omega_k R_m^k.
\label{specsinm2}
\end{equation}
Due to \eqref{sinmono4} and \eqref{sinmono3} we have a subsequencial limit $V\in {\bf IV}_k(B_1)$ which is stationary, and
we may assume that $b_m\rightarrow b_*\in T$ with $|b_*|\leq 1/2$ 
and $R_m\rightarrow R_*$ with
$1/10\leq R_*\leq 1/2$.  By \eqref{specsinm2} we also have
\begin{equation}
\|V\|(B_{R_*}(b_*))\leq \lambda \omega_k R_*^k.
\label{specsinm3}
\end{equation}
By \eqref{sinmono2}, we also have
\begin{equation}
\int_{G_k(B_1)}\|S-T\|\, dV(x,S)=0.
\label{specsinm4}
\end{equation}
Since $V$ is stationary and integral, one can show with \eqref{specsinm4} that 
$V$ is invariant in $T$-direction in $B_1$. 
Lemma \ref{lemsinmo} shows that $0\in {\rm spt}\, \|V\|$, 
thus $b_*\in {\rm spt}\, \|V\|$. In particular we have $\Theta^k(\|V\|, b_*)
\geq 1$. By the monotonicity formula we have $\|V\|(B_{R_*}(b_*))\geq
\omega_k R_*^k$. This contradicts with \eqref{specsinm3} and
we may conclude the proof. 
\hfill{$\Box$}
\begin{lemma}
Corresponding to
$1/2<\lambda<1$, $1\leq l<\infty$ and $1\leq E_1<\infty$,
there exists $\gamma>0$ with the following property. For $V\in {\bf IV}_k(C(T,3R))$, 
and two distinct points $y_1,\,y_2\in C(T,R)$, assume
\begin{equation}
|y_1-y_2|\leq l |T^{\perp}(y_1-y_2)|,
\label{twomono1}
\end{equation}
\begin{equation}
\Theta^k(\|V\|,y_i)=1 \mbox{ for } i=1,\, 2,
\label{twomono2}
\end{equation}
\begin{equation}
\int_{G_k(B_r(y_i))}\|S-T\|\, dV(x,S)\leq \gamma\|V\|(B_r(y_i)),
\hspace{.3cm}0<\forall r\leq 2R,\, i=1,\,2,
\label{twomono3}
\end{equation}
\begin{equation}
r\|\delta V\|(B_r(y_i))\leq \gamma \|V\|(B_r(y_i)),
\hspace{.3cm}0<\forall r\leq 2R,\,i=1,\,2,
\label{twomono4}
\end{equation}
\begin{equation}
\| V\|(B_r(y_i))\leq \omega_k r^k E_1,\hspace{.3cm}0<\forall r \leq 2R,\, i=1,\,2.
\label{twomono5}
\end{equation}
For $i=1,\, 2$ set $\tilde{y}_i:=T^{\perp}(y_i)\in {\mathbb R}^n$. 
Then we have
\begin{equation}
\|V\|(B_R(\tilde{y}_1)\cup B_R(\tilde{y}_2))\geq 2\lambda \omega_k R^k.
\label{twomono7}
\end{equation}
\label{lemsinmo3}
\end{lemma}
\begin{rem}
The above Lemma resembles \cite[6.2]{Allard} which is one of the essential
ingredients for graphical Lipschitz approximation needed for Allard's 
regularity theory. Note that the difference of above claim is that $\tilde{y}_1$
and $\tilde{y}_2$ are projected points to $T^{\perp}$ and the two balls centered 
at these points are inside of the cylinder $C(T,R)$. Roughly speaking, the Lemma
claims that if there are two good points placed in vertical positions, there must be
almost two parallel sheets inside of $C(T,R)$. 
\end{rem}
{\it Proof}.
Without loss of generality, we may assume $R=1$.
Suppose that the claim were false. Then we would have a sequence of 
$V_m\in {\bf IV}_k(C(T,3))$ and $y_{1,m}\neq y_{2,m}$ in $C(T,1)$
such that \eqref{twomono1}-\eqref{twomono5} are satisfied with $V=V_m$, 
$y_1=y_{1,m}$, $y_2=y_{2,m}$ and
$\gamma=1/m$ while \eqref{twomono7} does not hold, i.e., 
\begin{equation}
\|V_m\|(B_1(\tilde{y}_{1,m})\cup B_1(\tilde{y}_{2,m}))<2\lambda\omega_k,
\label{twomono8}
\end{equation}
where $\tilde{y}_{i,m}:=T^{\perp}(y_{i,m})$ for $i=1,\, 2$. 
For $m\in {\mathbb N}$, $0< s\leq 1$ and $i=1,\, 2$ set
\begin{equation}
A_{i,m}^s:=\{x\, :\, {\rm dist}\, (y_{i,m}-s T(y_{i,m}),x)<s\},\hspace{.3cm}
A_m^s:=A_{1,m}^s\cup A_{2,m}^s.
\label{twomono9}
\end{equation}
Using $y_{i,m}\in C(T,1)$, one can check that $A^s_{m}\subset C(T,1)$.
One can also see that $A_m^1=B_1(\tilde{y}_{1,m})\cup B_1(\tilde{y}_{2,m})$. 
We note that 
\begin{equation}
\liminf_{s\rightarrow 0}\frac{1}{\omega_k s^k}\|V_m\|(A_{i,m}^s)\geq \frac{\lambda+1}{2}\,\,(>\lambda)
\label{twomono10}
\end{equation}
for all sufficiently large $m$ and $i=1,\, 2$. 
To see this, for small $s>0$ we consider 
$V_{m,s}=(\tau(s^{-1}))_{\#}V_m$ and parallel translate $s^{-1}y_{i,m}$ to the origin. 
Such change of coordinates transforms $A_{i,m}^s$ to $B_1(-T(y_{i,m}))$. 
Then by Lemma \ref{lemsinmo2} with $b=-T(y_{i,m})$, $R=2$, $\lambda$ there
replaced by $(\lambda+1)/2$ and \eqref{twomono2}-\eqref{twomono5}, we obtain that
$\|V_{m,s}\|(B_1(-T(y_{i,m})))\geq \frac{\lambda+1}{2}\omega_k$ for all large $m$. 
This shows \eqref{twomono10} after changing back to the original coordinates. 
Next, let $R_m$ be the supremum of $r$ such that 
\begin{equation}
\|V_m\|(A^s_m)\geq 2\lambda \omega_k s^k
\label{twomono11}
\end{equation}
holds for all $0<s\leq r$. By \eqref{twomono10}, for all sufficiently large $m$, we have
\begin{equation}
\liminf_{s\rightarrow 0}\frac{1}{\omega_k s^k}\|V_m\|(A^s_m)\geq \sum_{i=1,\, 2}\liminf_{s\rightarrow 0}
\frac{1}{\omega_k s^k}\|V_m\|(A^s_{i,m})\geq (\lambda +1)> 2\lambda.
\label{twomono12}
\end{equation}
Thus \eqref{twomono12} shows that $R_m>0$. On the other hand, since 
$A_m^1=B_1(\tilde{y}_{1,m})\cup B_1(\tilde{y}_{2,m})$, \eqref{twomono8} shows $R_m<1$. 
 In particular, by the definition of $R_m$, we have
 \begin{equation}
 \|V_m\|(A_m^{R_m})=2\lambda \omega_k R_m^k
 \label{twomono13}
 \end{equation}
 for all large $m$. Set $\hat{V}_m:=(\tau(R_m^{-1}))_{\#} V_m$ and $\hat{y}_{i,m}
 :=R_m^{-1} y_{i,m}$. The change of variables
 $\tau(R_m^{-1})$ transforms $A_{i,m}^{R_m}$ to $B_1(\hat{y}_{i,m}-T(y_{i,m}))$,
 which we denote by $\hat{B}_{i,m}$. Suppose that 
 $\hat{B}_{1,m_j}\cap \hat{B}_{2,m_j}=\emptyset$ for
 some subsequence $\{m_j\}_{j=1}^{\infty}$ (subsequently omitting the sub-index).
 By Lemma \ref{lemsinmo2} with
 the origin replaced by $\hat{y}_{i,m}$, $b$ replaced by $\hat{y}_{i,m}-T(y_{i,m})$ and 
 $R=2$, for all sufficiently large $m$, we may conclude that 
 $\|\hat{V}_{m}\|(\hat{B}_{i,m})\geq (\lambda+1)\omega_k /2$. Here we have 
 used \eqref{twomono2}-\eqref{twomono5} with $\gamma=1/m$. After
 changing back to the original coordinates, we have $\|V_{m}\|(A_{i,m}^{R_m})
 \geq (\lambda+1)\omega_k R_m^k /2$.  Since we are now assuming $A_{1,m}^{R_m}
 \cap A_{2,m}^{R_m}=\emptyset$, this contradicts with \eqref{twomono13}. 
 Suppose that $\hat{B}_{1,m}\cap\hat{B}_{2,m}\neq \emptyset$
 for all large $m$. By shifting $x\rightarrow x-\hat{y}_{1,m}$ 
 and re-defining $\hat{B}_{i,m}$ as $\hat{B}_{i,m}
 -\hat{y}_{1,m}$ for $i=1,\, 2$, we may assume that
 $\hat{B}_{1,m}\cup\hat{B}_{2,m}\subset B_4$. 
 We have $\hat{B}_{1,m}=B_1(-T(y_{1,m}))$ and $\hat{B}_{2,m}
 =B_1(\hat{y}_{2,m}-\hat{y}_{1,m}-T(y_{2,m}))$. We also shift $\hat{V}_m$ by $x\rightarrow
 x-\hat{y}_{1,m}$.
 Since $\hat{y}_{2,m}-\hat{y}_{1,m}\in B_4$ and $|T(y_{i,m})|\leq 1$ for $i=1,\, 2$, there
 exists a subsequence (denoted by the same index) such that $\hat{y}_{2,m}-\hat{y}_{1,m}
 \rightarrow \hat{y}_*$, $T(y_{i,m})\rightarrow b_i \in T$ for $i=1,\, 2$, and
 $\hat{V}_m\rightarrow V$. We also have $\hat{B}_{1,m}\cup\hat{B}_{2,m}
 \rightarrow B_1(-b_1)\cup B_1(\hat{y}_*-b_2)$ in Hausdorff distance. 
 
 Suppose $T^{\perp}(\hat{y}_*)=0$. Then 
 $\lim _{m\rightarrow \infty} |T^{\perp}(\hat{y}_{2,m}-\hat{y}_{1,m})|=0$, and
 we would have by \eqref{twomono1} $\lim_{m\rightarrow\infty}|\hat{y}_{1,m}-
 \hat{y}_{2,m}|=0$. This implies $\lim_{m\rightarrow \infty}|T(y_{1,m}
 -y_{2,m})|=0$, thus we conclude that $b_1=b_2$, $\hat{y}_*=0$ and 
 $B_1(-b_1)\cup B_1(\hat{y}_*-b_2)=B_1(-b_1)$. By 
 \eqref{twomono2}-\eqref{twomono5}, Lemma \ref{lemsinmo} shows
 $0\in {\rm spt}\, \|V\|$. Also $V$ is stationary integral varifold with  
 the same condition as \eqref{specsinm4}, thus $V$ is invariant in 
 $T$-direction on $B_1(-b_1)$ as in the proof of Lemma \ref{lemsinmo2}. 
 By the definition of $R_m$ as in \eqref{twomono11}, we have
 \begin{equation}
 \|V\|(\{x\, :\, {\rm dist}\, (-s b_1,x)<s\})\geq 2\lambda\omega_k s^k
 \label{twomono14}
 \end{equation}
 for all $0<s\leq 1$ and equality holds for $s=1$. 
Letting $s\rightarrow 0$, due to the invariance of $\|V\|$ in $T$-direction
and $b_1\in T$,
this implies that $\Theta^k(\|V\|,0)\geq 2\lambda$. 
Also note that $\Theta^k(\|V\|,x)\in {\mathbb N}$ is constant
in $T$. Since it is an integer and $\lambda>1/2$, 
we have $\Theta^k(\|V\|,0)\geq 2$. Then this implies that $\Theta^k(\|V\|,
-b_1)\geq 2$ because of $-b_1\in T$. By the
monotonicity formula, then, we have $\|V\|(B_1(-b_1))\geq 2\omega_k$, 
which contradicts with the equality of \eqref{twomono14} when $s=1$. 

Next suppose $T^{\perp}(\hat{y}_*)\neq 0$. By the same line of proof 
using Lemma \ref{lemsinmo}, we have $\Theta^k(\|V\|,0)\geq 1$ and 
$\Theta^k(\|V\|,\hat{y}_*)\geq 1$.
Thus, the distinct affine planes containing $-b_1$ and $\hat{y}_*-b_2$ respectively
have positive integer multiplicities for $V$. This shows 
$\|V\|(B_1(-b_1)\cup B_1(\hat{y}_*-b_2))\geq 2\omega_k$. By the definition 
of $R_m$, we on the other hand have the same quantity $= 2\lambda\omega_k$,
a contradiction. This concludes the proof.
\hfill{$\Box$}
\begin{prop}
Corresponding to $1\leq E_1<\infty$ and $0<\nu<1$, 
there exist $0<\alpha_1<1$, $0<\beta_1<1$ and $1<P_1<\infty$ with the following property. For
$V\in {\bf IV}_k(C(T,6R))$ which is finite and which is of unit density, suppose $V=|M|$ and identify
$T$ with ${\mathbb R}^k\times\{0\}$.
Define
\begin{equation}
\alpha^2:=\frac{1}{R^{k-2}}\int_{C(T,6R)}|h(V,x)|^2\, d\|V\|(x),
\label{slip1}
\end{equation}
\begin{equation}
\beta^2:=\frac{1}{R^k}\int_{G_k(C(T,6R))}\|S-T\|^2\, dV(x,S).
\label{slip2}
\end{equation}
Suppose that
\begin{equation}
\|V\|(B_r(x))\leq \omega_k r^k E_1 ,\hspace{.3cm}0<\forall r\leq 4R,\,\,\forall x\in M\cap C(T,2R),
\label{slip3-5}
\end{equation}
\begin{equation}
\frac{1}{\omega_k (2R)^k}\|V\|(C(T,2R))\leq 2-\nu
\label{slip4}
\end{equation}
and
\begin{equation}
\frac{1}{ \omega_k R^k}\|V\|(C(T,R))\geq \nu.
\label{slip5}
\end{equation}
Then either of the following (1) or (2) holds. 
\newline
(1) We have $\alpha\geq \alpha_1$ or $\beta\geq \beta_1$.\newline
(2) There exist Lipschitz maps $f\,:\, T\rightarrow {\mathbb R}^{n-k}$ and 
$F\,:\, T\rightarrow {\mathbb R}^n$ such that 
\begin{equation}
F(x)=(x,f(x)),\hspace{.3cm}\forall x\in T,
\label{slip6}
\end{equation}
\begin{equation}
|f(x_1)-f(x_2)|\leq |x_1-x_2|,\hspace{.3cm}\forall x_1,\, \forall x_2,\in T.
\label{slip7}
\end{equation}
Set $X:=M\cap {\rm image}\, F\cap C(T,R)$ and $Y:=T(X)$. Then
\begin{equation}
\frac{1}{R^k}\left(\|V\|(C(T,R)\setminus X)+{\mathcal L}^k(B_R^k\setminus Y)\right)\leq
\left\{\begin{array}{ll} P_1 (\alpha^{\frac{2k}{k-2}}+\beta^2) & \mbox{ if }k>2, \\
P_1 \beta^2 & \mbox{ if }k\leq 2.
\end{array}\right.
\label{slip8}
\end{equation}
\label{lipapxmn}
\end{prop}
\begin{rem}
If $h(V,\cdot)$ does not exist or $h(V,\cdot)\notin L^2(\|V\|)$, then we 
define $\alpha=\infty$ and we have trivially (1). Later we apply 
the result to $V_t$ for a.e$.$ time, where $h(V_t,\cdot)\in L^2(\|V_t\|)$.
\label{remh2}
\end{rem}
\begin{rem}
For $k>2$, if (1) holds, we may re-define $P_1$ sufficiently large depending on 
$\alpha_1$ and $\beta_1$ so that
\eqref{slip8} holds trivially. Thus, for $k>2$, we may consider that the case (1) 
is a subset of (2), with $F(x)=(x,0)$, for example. For $k\leq 2$, if $\beta\geq 
\beta_1$, then we may also choose $P_1$ sufficiently large similarly. Thus we 
may drop `or $\beta\geq \beta_1$' from (1). We still find it more convenient 
to write out (1) and (2) as above.
\end{rem}
\begin{rem}
The fact that the power $2k/(k-2)$ of $\alpha$ is strictly larger than $2$ is of
essential importance in proving Theorem \ref{poptheorem}, and this is the
reason that we work with $r\|\delta V\|(B_r)\leq \gamma \|V\|(B_r)$ 
instead of $\|\delta V\|(B_r)\leq \gamma\|V\|(B_r)$, which is done in Allard's
theory. 
\label{imprem}
\end{rem}
{\it Proof}. We may assume $R=1$ without loss of generality, and we may
assume that $\alpha<\infty$ by Remark \ref{remh2}. Let $\gamma>0$
be the smallest constant obtained in Lemma \ref{lemsinmo}-\ref{lemsinmo3} corresponding to 
\begin{equation}
\lambda=1-\frac{\nu}{2^{3k}},
\label{slipsm}
\end{equation}
$l=\sqrt{2}$, and $E_1$ as in the assumption.  
Let $M_g\subset M\cap C(T,2)$ be the set of points $y$ such that the
following \eqref{slip9}-\eqref{slip11} all hold, 
\begin{equation}
\Theta^k(\|V\|,y)=1,\,\,(\mbox{or equivalently, }\Theta^k({\mathcal H}^k\lfloor_{M},y)=1)
\label{slip9}
\end{equation}
\begin{equation}
r\|\delta V\|(B_r(y))\leq \gamma \|V\|(B_r(y)),\hspace{.3cm}0<\forall r\leq 4,
\label{slip10}
\end{equation}
\begin{equation}
\int_{G_k(B_r(y))}\|S-T\|^2\, dV(y,S)\leq \gamma^2 \|V\|(B_r(y)),\hspace{.3cm}
0<\forall r \leq 4.
\label{slip11}
\end{equation}
By Schwarz' inequality, \eqref{slip11} implies
\begin{equation}
\int_{G_k(B_r(y))}\|S-T\|\, dV(y,S)\leq \gamma \|V\|(B_r(y)),\hspace{.3cm}
0<\forall r\leq 4.
\label{slip12}
\end{equation}
Take any $x\in T(M_g)\subset B_2^k$. By definition there exists 
at least one $y\in M_g$ such that $T(y)=x$. Suppose that there are two such
$y$, i.e., there are two distinct $y_1,\, y_2\in M_g$ such that $T(y_1)=T(y_2)=x$.
Then $y_1$, $y_2$ satisfy all the assumptions of Lemma \ref{lemsinmo3},
namely, we have $|y_1-y_2|=|T^{\perp}(y_1-y_2)|$, thus \eqref{twomono1}
is satisfied with $l=\sqrt{2}$, and \eqref{twomono2}-\eqref{twomono5} are satisfied with $R=2$
due to \eqref{slip3-5}, \eqref{slip9}, \eqref{slip10} and \eqref{slip12}. Since $B_2(\tilde{y}_1)
\cup B_2(\tilde{y}_2)\subset C(T,2)$, \eqref{twomono7} and \eqref{slipsm}
shows 
\begin{equation}
\|V\|(C(T,2))\geq 2^{k+1}\lambda\omega_k=2^{k+1}\{1-(\nu/2^{3k})\}\omega_k.
\label{slipsm2}
\end{equation}
This contradicts with
\eqref{slip4}. Since there is only one such $y\in M_g$ with $T(y)=x$,
we may define a map $f\, :\, T(M_g)
\rightarrow {\mathbb R}^{n-k}$ so that $f(x)=T^{\perp}(y)$. Next, 
for any $y_1\neq y_2$ in $M_g$, if $|y_1-y_2|\leq \sqrt{2}|T^{\perp}(y_1-y_2)|$,
similar reasoning using Lemma \ref{lemsinmo3} leads to a contradiction with 
\eqref{slip4}. Thus we have $|y_1-y_2|>\sqrt{2}|T^{\perp}(y_1-y_2)|$. This 
shows that $|f(x_1)-f(x_2)|<|x_1-x_2|$ for all $x_1,\, x_2\in T(M_g)$. By Kirszbraun's
theorem \cite[2.10.43]{Federer}, we have an extension
$\tilde{f}\,:\,T\rightarrow{\mathbb R}^{n-k}$ of $f$ with the equal Lipschitz constant. 
Rewrite this $\tilde{f}$ as $f$. 
We then define $F\,:\, T\rightarrow{\mathbb R}^n$ by $F(x):=(x,f(x))$. 

For any $y\in (M\cap C(T,2))\setminus M_g$, at least one of \eqref{slip9}-\eqref{slip11} fails. The set of points for which
\eqref{slip9} fails is zero set with respect to ${\mathcal H}^k$, 
thus it is zero set with respect to $\|V\|$. 
Let $X'_1$ 
be the set of points for which \eqref{slip10} fails, and let $X'_2$ be the set of points for which 
\eqref{slip11} fails. 
We first estimate $\|V\|(X'_1)$. For $y\in X'_1$, either (A) there exists a sequence 
$R_m\rightarrow 0$, $0<R_m\leq 4$, such that \eqref{slip10} fails for $r=R_m$,
or (B) there exists some $0<R<4$ such that \eqref{slip10} holds for all $0<r\leq R$
and $R$ is the infimum of $r$ such that \eqref{slip10} fails. Either (A) or (B), 
there exists some $0<R(y)<4$ such that both
\begin{equation}
\omega_k R(y)^k/2\leq \|V\|(B_{R(y)}(y)),
\label{slip13}
\end{equation}
\begin{equation}
 \gamma \|V\|(B_{R(y)}(y))< R(y)\|\delta V\|(B_{R(y)}(y))
\label{slip14}
\end{equation}
hold. The reason is that, if (A) holds, we use \eqref{slip9} to have \eqref{slip13} and choose
small enough $R_m$ to have \eqref{slip14}. If (B) holds, by Lemma \ref{lemsinmo}
we have $\|V\|(B_R(y))\geq \{1-(\nu/2^{3k})\}\omega_k R^k$. Thus choosing $R(y)$ to be
slightly larger than $R$ if necessary, we have both \eqref{slip13} and \eqref{slip14}
satisfied. By H\"{o}lder's inequality applied to \eqref{slip14}, we have
\begin{equation}
\gamma^2 \|V\|(B_{R(y)}(y))< R(y)^2 \int_{B_{R(y)}(y)}|h(V,\cdot)|^2\, d\|V\|.
\label{slip15}
\end{equation}
Using \eqref{slip13} and \eqref{slip15}, we obtain
\begin{equation}
(\gamma^2/2)\omega_k R(y)^{k-2}< \int_{B_{R(y)}(y)}|h(V,\cdot)|^2\, d\|V\|.
\label{slip16}
\end{equation}
Consider the case $k\leq 2$. Since $y\in C(T,2)$ and $R(y)<4$, 
we have $B_{R(y)}(y)\subset C(T,6)$ and we have from \eqref{slip1} and \eqref{slip16}
\begin{equation}
\alpha^2> (\gamma^2/2)\omega_k R(y)^{k-2}> (\gamma^2/2^{5-2k})\omega_k.
\label{slip17}
\end{equation}
If $\alpha^2<(\gamma^2/2^{5-2k})\omega_k$ and $X'_1\neq \emptyset$, we would obtain a contradiction.
We set $\alpha_1$ so that $\alpha_1^2 < (\gamma^2/2^{5-2k})\omega_k$. 
Then, either we have $\alpha
\geq \alpha_1$, or else $\alpha<\alpha_1$ and $X'_1=\emptyset$. 
Thus for $k\leq 2$ and for the former case, we may end the proof, with case (1).
For the rest of the proof, for $k\leq 2$, assume that $\alpha<\alpha_1$. Next consider 
$k>2$. From \eqref{slip17}, we have
\begin{equation}
R(y)^2 < (2\alpha^2\gamma^{-2}\omega_k^{-1})^{\frac{2}{k-2}}.
\label{slip18}
\end{equation}
Substituting \eqref{slip18} into \eqref{slip15} gives
\begin{equation}
\|V\|(B_{R(y)}(y))< \gamma^{-2}
(2\alpha^2\gamma^{-2}\omega_k^{-1})^{\frac{2}{k-2}}\int_{B_{R(y)}(y)}|h(V,\cdot)|^2\,
d\|V\|.
\label{slip19}
\end{equation}
From \eqref{slip19}, we may choose $R'(y)< R(y)$ such that 
\begin{equation}
\|V\|(\overline{B_{R'(y)}(y)})<\gamma^{-2}(2\alpha^2\gamma^{-2}\omega_k^{-1})^{\frac{2}{k-2}}\int_{\overline{B_{R'(y)}(y)}}|h(V,\cdot)|^2\, d\|V\|.
\label{slip19-5}
\end{equation}
Thus the Besicovitch covering theorem with \eqref{slip19-5} implies that 
\begin{equation}
\|V\|(X'_1)\leq {\bf B}(n) \gamma^{-2}(2\alpha^2\gamma^{-2}\omega_k^{-1})^{\frac{2}{k-2}}
\int_{C(T,6)}|h(V,\cdot)|^2\, d\|V\|={\bf B}(n)\gamma^{\frac{-2k}{k-2}}
\omega_k^{\frac{-2}{k-2}}\alpha^{\frac{2k}{k-2}}.
\label{slip20}
\end{equation}
Here ${\bf B}(n)$ is the Besicovitch constant which depends only on $n$. 
For the estimate of $\|V\|(X'_2)$ (both $k\leq 2$ and $k>2$), 
the Besicovitch covering theorem and \eqref{slip2} imply that
\begin{equation}
\|V\|(X'_2)\leq {\bf B}(n)\gamma^{-2}\beta^2.
\label{slip21}
\end{equation}
Combining \eqref{slip20} and \eqref{slip21}, we obtain 
\begin{equation}
\|V\|(C(T,2)\setminus M_g)\leq \left\{\begin{array}{ll} {\bf B}(n)(\gamma^{\frac{-2k}{k-2}}
\omega_k^{\frac{-2}{k-2}}\alpha^{\frac{2k}{k-2}}+
\gamma^{-2}\beta^2),
& k>2,\\  {\bf B}(n) \gamma^{-2}\beta^2, & k\leq 2. \end{array}\right.
\label{slip22}
\end{equation}

Next we estimate ${\mathcal H}^k (B_1^k\setminus T(M_g))$. The set
$B_1^k\setminus T(M_g)$ may be regarded as a `hole of $\|V\|$' and we need to 
estimate the size in terms of $\alpha$ and $\beta$. First we claim 
that $T(M_g)\cap B_1^k\neq \emptyset$. Otherwise, $C(T,1)\cap M_g
=\emptyset$. From \eqref{slip5}, we have
\begin{equation}
\nu \omega_k \leq \|V\|(C(T,1)) \leq \|V\|(C(T,2)\setminus M_g).
\label{slip23}
\end{equation}
By \eqref{slip22}, if we restrict $\alpha$ and $\beta$ small depending only on
$n$, $k$, $E_1$ and $\nu$, we would have a contradiction to \eqref{slip23}. Thus we 
either have $\alpha\geq \alpha_1$ or $\beta\geq \beta_1$, where $\alpha_1$
and $\beta_1$ depend only on $n$, $k$, $E_1$ and $\nu$, or we have 
$T(M_g)\cap B^k_1\neq \emptyset$. The former case corresponds to the case (1)
of the conclusion. We assume for the rest of the proof that $\alpha<\alpha_1$
and $\beta<\beta_1$. Thus we have some $x\in M_g\cap C(T,1)$. By Lemma 
\ref{lemsinmo2} with the origin there replaced by $x$, $b$ there replaced by 
$T^{\perp}(x)$ and $R=4$, we conclude that
\begin{equation}
\begin{split}
\big(1-\frac{\nu}{2^{3k}}\big)2^k \omega_k &  \leq \|V\|(B_2(T^{\perp}(x)))
\leq \|V\|(C(T,2)) \\ & =\|V\|(C(T,2)\cap M_g)+\|V\|(C(T,2)\setminus M_g).
\end{split}
\label{slip24}
\end{equation}
On the other hand, since $M_g$ is a Lipschitz graph,
\begin{equation}
\begin{split}
\|V\|(C(T,2)\cap M_g) & =\int_{T(M_g)} |\Lambda_k \nabla F|\, d{\mathcal H}^k \\
   & \leq \int_{T(M_g)} (1+c(n,k)\|{\rm image}\, \nabla F-T\|^2) \, d{\mathcal H}^k  \\
   & \leq {\mathcal H}^k (T(M_g))+c(n,k)\int_{G_k(C(T,2))} \|S-T\|^2\, dV(x,S) \\
   & = 2^k \omega_k - {\mathcal H}^k (B_2^k \setminus T(M_g))+c(n,k)\beta^2.
\end{split}
\label{slip25}
\end{equation}   
Combining \eqref{slip24} and \eqref{slip25}, we obtain
\begin{equation}
{\mathcal H}^k (B_2^k\setminus T(M_g))\leq \frac{\nu}{2^{2k}}\omega_k
+c(n,k)\beta^2 +\|V\|(C(T,2)\setminus M_g).
\label{slip26}
\end{equation}
The estimates \eqref{slip22} and \eqref{slip26} show that
\begin{equation}
{\mathcal H}^k (B_2^k\setminus T(M_g))<\frac{1}{4^k}\omega_k
\label{slip27}
\end{equation}
holds by further restricting $\alpha_1$ and $\beta_1$ if necessary
depending on $1-\nu$ and $k$. 
We set $Q_1=B_1^k\setminus T(M_g)$ and $Q_2=B_2^k \setminus T(M_g)$. 
For any Lebesgue point $x\in Q_1$ of $Q_2$, due to \eqref{slip27}, there
exists $0<r(x)<1$ such that 
\begin{equation}
{\mathcal H}^k (Q_2\cap B_{r(x)}^k(x))=\frac{1}{4^k} \omega_k r(x)^k.
\label{slip28}
\end{equation}
By the Besicovitch covering theorem, there exists a set of mutually disjoint balls
$\{B_{r(x_i)}(x_i)\}_{i\in I}$ ($I$ : a finite set or ${\mathbb N}$) such that
$x_i\in Q_1$ and, denoting $U=\cup_{i\in I}B_{r(x_i)}(x_i)$, 
\begin{equation}
{\mathcal H}^k (Q_1)\leq {\bf B}(k){\mathcal H}^k (U).
\label{slip29}
\end{equation}
We also estimate just like \eqref{slip25} that
\begin{equation}
\begin{split}
\|V\|(T^{-1}(U)) & \leq {\mathcal L}^k(T(M_g)\cap U)
+c(n,k)\beta^2 +\|V\|(C(T,2)\setminus M_g) \\
& = {\mathcal H}^k (U)-{\mathcal H}^k (Q_2\cap U)
+c(n,k)\beta^2 + \|V\|(C(T,2)\setminus M_g).
\end{split}
\label{slip30}
\end{equation}
To obtain a lower bound for \eqref{slip30}, note that \eqref{slip28} guarantees
that there exists $y_i\in M_g\cap C(T,2)$ for each $i\in I$ such that $T(y_i)
\in B^k_{r(x_i)/2} (x_i)$. Then Lemma \ref{lemsinmo2} with the origin there replaced
by $y_i$, $b$ there replaced by $x_i+T^{\perp}(y_i)$ (so that $y_i-x_i-T^{\perp}(y_i)
\in T$) and $R=2r(x_i)$, we conclude that
\begin{equation}
\big(1-\frac{\nu}{2^{3k}}\big)\omega_k r(x_i)^k\leq \|V\|(B_{r(x_i)}(x_i+T^{\perp}
(y_i)))\leq \|V\|(T^{-1}(B_{r(x_i)}^k(x_i))).
\label{slip31}
\end{equation}
Then summing over $i\in I$, and \eqref{slip30} and \eqref{slip31} show that
\begin{equation}
{\mathcal H}^k (Q_2\cap U)\leq \frac{\nu}{2^{3k}} {\mathcal H}^k (U)+c(n,k)\beta^2
+\|V\|(C(T,2)\setminus M_g).
\label{slip32}
\end{equation}
By \eqref{slip28}, we have
\begin{equation}
\frac{1}{4^k} {\mathcal H}^k (U)={\mathcal H}^k (Q_2\cap U),
\label{slip33}
\end{equation}
and combining \eqref{slip32} and \eqref{slip33}, we obtain
\begin{equation}
\big(\frac{1}{4^k}-\frac{\nu}{2^{3k}}\big){\mathcal H}^k (U)\leq c(n,k)\beta^2+\|V\|(C(T,2)\setminus
M_g).
\label{slip34}
\end{equation}
Thus \eqref{slip29} and \eqref{slip34} give
\begin{equation}
{\mathcal H}^k (Q_1)\leq {\bf B}(k)\big(\frac{1}{4^k}-\frac{\nu}{2^{3k}}\big)^{-1}(c(n,k)\beta^2
+\|V\|(C(T,2)\setminus M_g)).
\label{slip35}
\end{equation}
Finally, setting $X=M\cap {\rm image}\, F\cap C(T,1)$ and $Y=T(X)$, we note that
$C(T,1)\cap M_g\subset X$, $\|V\|(C(T,1)\setminus X)\leq \|V\|(C(T,2)\setminus M_g)$,
$B_1^k\setminus Y\subset B_1^k\setminus T(M_g)=Q_1$. Thus \eqref{slip22} and
\eqref{slip35} give the desired estimate \eqref{slip8} for a suitable choice of $P_1$. 
\hfill{$\Box$}
\section{Energy estimates}
The aim of this section is to prove Theorem \ref{poptheorem}, which claims
that the deviation of $k$-dimensional area of moving varifold from that of
flat $k$-dimensional plane can be estimated in terms of $L^2$-height.
In Allard's regularity theory, the analogy of this part may be
the estimate of tilt-excess in terms of $L^2$-height. However, compared to
Allard's theory, the proof is more involved. All we have control is
the rate of change of mass in time. The idea is that when the deviation 
is large, the mean curvature is relatively large. This helps pushing down
the mass more quickly, and the mass becomes very close to that of flat
plane in finite time. The set of estimates in this section is a crucial foundation for 
parabolic Lipschitz approximation and blow-up
technique. The content corresponds to 
\cite[6.6 `Popping soap films']{Brakke} even though the details have been largely
changed and various alterations are made. Some of the results used are relegated
to Appendix for convenience of the readers.

\begin{define}
Fix $\phi\in C^{\infty}([0,\infty))$ such that $0\leq \phi\leq 1$,
\begin{equation}
\phi(x)\left\{\begin{array}{ll}
=1 & {\mbox for }\,\,0\leq x\leq (2/3)^{1/k},\\
>0 & {\mbox for }\,\, 0\leq x<(5/6)^{1/k}, \\
=0 & {\mbox for }\,\, x\geq (5/6)^{1/k}.
\end{array}\right.
\label{pops1}
\end{equation}
For $0<R<\infty$, $x\in {\mathbb R}^n$ and $T\in {\bf G}(n,k)$ define 
\begin{equation}
\phi_{T,R}(x):=\phi\left(R^{-1}|T(x)|\right),\,\,\, \phi_T(x):=\phi_{T,1}(x)=\phi(|T(x)|).
\label{Tphi}
\end{equation}
Set
\begin{equation}
{\bf c}:=\int_T \phi_T^2(x)\, d{\mathcal H}^k (x).
\label{pops1.5}
\end{equation}
\end{define}
In the next Proposition, we use
\begin{equation}
\frac23 \omega_k<{\bf c}<\frac56 \omega_k
\label{pops2}
\end{equation}
which follows easily from \eqref{pops1}. 
\begin{prop}
Corresponding to $1\leq E_1<\infty$ and $0<\nu<1$ there exist 
$0<\alpha_2<1$, $0<\mu_1<1$ and $1<P_2<\infty$ with the following property.
For $V\in {\bf IV}_k(C(T,1))$ which is finite and which is of unit density with
$V=|M|$ and for $T\in {\bf G}(n,k)$, define
\begin{equation}
\alpha^2:=\int_{C(T,1)}|h(V,x)|^2 \phi_T^2(x)\, d\|V\|(x),
\label{pops3}
\end{equation}
\begin{equation}
\mu^2:= \int_{C(T,1)} |T^{\perp}(x)|^2\, d\|V\|(x).
\label{pops5}
\end{equation}
Suppose 
\begin{equation}
\|V\|(B_r(x))\leq \omega_k r^k E_1, \hspace{.3cm}\forall B_r(x)\subset C(T,1).
\label{pops5.5}
\end{equation}
(A) Suppose 
\begin{equation}
\left| \|V\|( \phi_T^2)-{\bf c}\right|\leq \frac{\bf c}{8},\,\,\,\alpha\leq \alpha_2\,\,\,\mbox{and}\,\,\,\mu\leq \mu_1.
\label{pops5.7}
\end{equation}
Then we have
\begin{equation}
\left|\|V\|(\phi_T^2)-{\bf c}\right| \leq
\left\{\begin{array}{ll} P_2(\alpha^{\frac{2k}{k-2}}+\alpha^{\frac32}\mu^{\frac12}
+\mu^2) & \mbox{if }\, k>2, \\
P_2 (\alpha^{\frac32}\mu^{\frac12}+\mu^2) & \mbox{if }\, k\leq 2.
\end{array}\right.
\label{pops6}
\end{equation}
(B) Suppose 
\begin{equation}
\frac{\bf c}{8}<\left|\|V\|(\phi_T^2)-{\bf c}\right|\leq (1-\nu){\bf c}\,\,\mbox{ and }\,\,\,
\mu\leq \mu_1.
\label{popsnu}
\end{equation}
Then we have $\alpha\geq \alpha_2$.
\label{popslem1}
\end{prop}
{\it Proof}. First consider case (A).
Define
\begin{equation}
\beta^2:=\int_{G_k(C(T,1))} \|S-T\|^2 \phi_T^2 \, dV(\cdot,S).
\label{pops4}
\end{equation}
By Proposition \ref{cylgrolem} with $R_1=1/12$ and $R_2=1$, we have
\begin{equation}
12^k \|V\|(\phi_{T,1/12}^2)-\|V\|(\phi_T^2)
\geq -12^k(k\beta^2\log 12+\alpha\beta+\beta^2).
\label{pops7}
\end{equation}
Since $\phi_{T,1/12}\leq 1$ and by \eqref{pops5.7}, \eqref{pops7} and \eqref{pops2}, we have
\begin{equation}
12^k \|V\|(C(T,1/12))\geq \frac78 {\bf c}-c(k)(\alpha\beta+\beta^2)>\frac{7}{12}\omega_k
-c(k)(\alpha\beta+\beta^2).
\label{pops8}
\end{equation}
By restricting $\alpha$ and $\beta$ depending only on $k$, \eqref{pops8} guarantees that
\begin{equation}
\|V\|(C(T,1/12))\geq \frac12 \omega_k\frac{1}{12^k}.
\label{pops9}
\end{equation}
Similarly for $R_1=(3/2)^{1/k}/6$ and $R_2=1$, Proposition \ref{cylgrolem} gives
\begin{equation}
R_1^{-k}\|V\|(\phi_{T,R_1}^2) -\|V\|(\phi_T^2) 
\leq R_1^{-k}(k\beta^2 \log R_1^{-1}+\alpha\beta+\beta^2).
\label{pops10}
\end{equation}
By \eqref{pops1}, we have $\{\phi_{T,R_1}=1\}\supset C(T,1/6)$. Thus \eqref{pops10}, 
\eqref{pops5.7} and \eqref{pops2} give
\begin{equation}
R_1^{-k}\|V\|(C(T,1/6))\leq c(k)(\beta^2+\alpha\beta)+\frac98 {\bf c}
\leq c(k)(\beta^2+\alpha\beta)+\frac{15}{16}\omega_k.
\label{pops11}
\end{equation}
With the definition of $R_1$ and \eqref{pops11}, we have
\begin{equation}
\|V\|(C(T,1/6))\leq c(k) (\beta^2+\alpha\beta)+\frac{45}{32}\frac{1}{6^k}\omega_k.
\label{pops12}
\end{equation}
Since $45/32<3/2$, by restricting $\alpha$ and $\beta$ depending only on $k$, we have
from \eqref{pops12}
\begin{equation}
\|V\|(C(T,1/6))\leq \frac32\omega_k \frac{1}{6^k}.
\label{pops13}
\end{equation}
We now use Proposition \ref{lipapxmn} with $R=1/12$ and $\nu=1/2$. Due to \eqref{pops5.5}, 
\eqref{pops9} and \eqref{pops13} , the assumptions
\eqref{slip3-5}, \eqref{slip4} and \eqref{slip5} are satisfied. The smallness
of $\alpha$ and $\beta$ in \eqref{slip1} and \eqref{slip2} can be guaranteed
by assuming $\alpha$ and $\beta$ in \eqref{pops3} and \eqref{pops4} are
small due to $C(T,1/2)\subset\{\phi_T=1\}$. 
Thus there exists a Lipschitz function $F:T\rightarrow {\mathbb R}^n$ and $\tilde{P}_1$
(which is a suitable multiple of $P_1$ by a constant depending only on $k$)
such that
\begin{equation}
\|V\|(C(T,1/12)\setminus X)+{\mathcal H}^k (B_{1/12}^k\setminus Y)\leq 
\left\{\begin{array}{ll} \tilde{P}_1 (\alpha^{\frac{2k}{k-2}}+\beta^2) & \mbox{if }\, k>2,\\
\tilde{P}_1 \beta^2 & \mbox{if }\, k\leq 2.
\end{array}\right.
\label{pops14}
\end{equation}
Here $X=C(T,1/12)\cap M\cap {\rm image}\, F$ and $Y=T(X)$. 
By using this approximation, we estimate as
\begin{equation}
 \int_{C(T,1/12)}\phi_{T,1/12}^2 \, d\|V\|=\int_{C(T,1/12)\setminus X}\phi_{T,1/12}^2\, d\|V\|
+\int_{Y}\phi^2_{T,1/12}|\Lambda_k \nabla F|\, d{\mathcal H}^k.
\label{pops15}
\end{equation}
The first term on the right-hand side of \eqref{pops15}
is bounded by \eqref{pops14}. The error of replacing $Y$
by $B_{1/12}^k$ for the second term is also bounded similarly. Just like estimating \eqref{slip25},
we may obtain 
\begin{equation}
\left| \int_{C(T,1/12)}\phi_{T,1/12}^2\, d\|V\|-\int_T \phi_{T,1/12}^2 \, d{\mathcal H}^k
\right|\leq \left\{\begin{array}{ll} \tilde{P}_1 (\alpha^{\frac{2k}{k-2}}+\beta^2) & \mbox{if }\, k>2,\\
\tilde{P}_1 \beta^2 & \mbox{if }\, k\leq 2.
\end{array}\right.
\label{pops16}
\end{equation}
where $\tilde{P}_1$ may differ from before only by a multiple of constant depending only on $n$ and
$k$. We now apply Proposition \ref{cylgrolem} again. Observing that 
$12^k\int_T\phi_{T,1/12}^2\, d{\mathcal H}^k ={\bf c}$
and by \eqref{pops16}, yet with another $\tilde{P}_1$, we obtain
\begin{equation}
\left|\|V\|(\phi_T^2(x))-{\bf c}\right|\leq 
\left\{\begin{array}{ll} \tilde{P}_1 (\alpha^{\frac{2k}{k-2}}+\beta^2+\alpha\beta) & \mbox{if }\, k>2,\\
\tilde{P}_1 (\beta^2+\alpha\beta) & \mbox{if }\, k\leq 2.
\end{array}\right.
\label{pops17}
\end{equation}
Finally, with Proposition \ref{tiltexlem}, we have 
\begin{equation}
\beta^2\leq 4\alpha\mu+c(\phi)\mu^2\,\,\mbox{and }\,\, \beta\leq 2\alpha^{\frac12}\mu^{\frac12}+
c(\phi)\mu
\label{pops18}
\end{equation}
where $c(\phi)\geq 1$ depends only on $\sup|\phi'|$ which we may consider to be constant. 
Thus with suitable restrictions on $\alpha$ and $\mu$, $\beta$ is considered small.
By Young's inequality, we also have $\alpha\mu\leq c(\alpha^{3/2}\mu^{1/2}+\mu^2)$
and this combined with \eqref{pops17} and \eqref{pops18} proves \eqref{pops6} with
suitable choices of $P_2$, $\alpha_2$ and $\mu_1$.

For case (B), assume the claim were false. Then for each $m\in{\mathbb N}$ there exists
$V_m=|M_m|$ satisfying \eqref{popsnu} with $V$ there replaced by $V_m$, while
\begin{equation}
\mu_m^2:=\int_{C(T,1)}|T^{\perp}(x)|^2\, d\|V_m\|<\frac{1}{m},\,\,\mbox{and}\,\,\,
\alpha_m^2:=\int_{C(T,1)}|h(V_m,\cdot)|^2\phi_T^2\,d\|V_m\|<\frac{1}{m}.
\label{nu3}
\end{equation}
By Proposition \ref{tiltexlem}, we also have the corresponding $\beta_m\rightarrow 0$.
By \eqref{pops5.5}, there exists a convergent subsequence (denoted by the same 
index) $\{V_m\}$ and its limit $V$ on $C(T,1)$. By \eqref{nu3}, on $\{\phi_T>0\}$,
$V$ is stationary and integral. Since the corresponding $\beta$ and $\mu$ for $V$ vanishes, 
$V\lfloor_{\{\phi_T>0\}}=i|T\cap C(T,1)|$ for some integer $i$. On the other hand, $V$
satisfies \eqref{popsnu} (with possible $\leq $ sign on the left-hand side). But there is no
such integer satisfying ${\bf c}/8\leq |i{\bf c}-{\bf c}|\leq (1-\nu){\bf c}$. Hence we obtain a 
contradiction and we may conclude the proof by choosing smaller $\mu_1$ and $\alpha_2$ 
from (A) and (B).
\hfill{$\Box$}
\begin{cor}
Let $\alpha_2$, $\mu_1$ and $P_2$ be the same constants as in 
Proposition \ref{popslem1}. Set $\mu_2=\min\{\mu_1,({\bf c}/32P_2)^{1/2}\}$. 
For $V$ and $T$ as in Proposition \ref{popslem1} define $\alpha$ and
$\mu$ as in 
\eqref{pops3} and \eqref{pops5}. Define
\begin{equation}
\hat{E}:=\|V\|(\phi_T^2)-{\bf c}.
\label{pops19}
\end{equation}
Assume \eqref{pops5.5}, 
\begin{equation}
\mu\leq \mu_2
\label{popsext2}
\end{equation}
and 
\begin{equation}
2P_2 \mu^2\leq |\hat{E}|\leq (1-\nu){\bf c}.
\label{pops20}
\end{equation}
Then we have
\begin{equation}
\alpha^2\geq \left\{\begin{array}{ll} \min\left\{\alpha_2^2,\,
(4P_2)^{-\frac{k-2}{k}}|\hat{E}|^{\frac{k-2}{k}},\, (4P_2)^{-\frac43} \mu^{-\frac23}
|\hat{E}|^{\frac43}\right\} & \mbox{if }\, k>2, \\
\min\left\{\alpha_2^2,\, (2P_2)^{-\frac43}\mu^{-\frac23}|\hat{E}|^{\frac43}\right\}
& \mbox{if }\, k\leq 2.\end{array}\right.
\label{pops21}
\end{equation}
\label{popslem2}
\end{cor}
\begin{rem} 
Note that the mean curvature square term $\alpha^2$ gives a lower bound for
the rate of mass decrease, and Corollary \ref{popslem2} relates that to the
mass itself. This allows one to introduce a differential inequality for the mass,
see Lemma \ref{popslem3}.
\end{rem}
{\it Proof}. We only prove the case $k>2$ since $k\leq 2$ can be handled similarly.
If $\alpha\geq \alpha_2$, then we obtain
\eqref{pops21} trivially. Thus assume $\alpha<\alpha_2$. The case (B) of Proposition
\ref{popslem1} do not occur under this assumption. The assumptions for
Proposition \ref{popslem1}, \eqref{pops5.5} and \eqref{pops5.7},
are all satisfied due to 
\eqref{pops19}, \eqref{popsext2} and \eqref{pops20}. Then we have 
\begin{equation}
|\hat{E}|\leq P_2(\alpha^{\frac{2k}{k-2}}+ \alpha^{\frac32}\mu^{\frac12}+\mu^2).
\label{pops22}
\end{equation}
By \eqref{pops20} and \eqref{pops22} we conclude that
\begin{equation}
|\hat{E}|\leq 2P_2(\alpha^{\frac{2k}{k-2}}+\alpha^{\frac32}\mu^{\frac12}).
\label{pops23}
\end{equation}
The inequality \eqref{pops23} implies that we have either 
$|\hat{E}|\leq 4P_2\alpha^{2k/(k-2)}$ or $|\hat{E}|\leq 4P_2 \alpha^{3/2}\mu^{1/2}$.
Thus we have
\begin{equation}
(4P_2)^{-\frac{k-2}{k}}|\hat{E}|^{\frac{k-2}{k}}\leq \alpha^2\,\,\mbox{or }\,\,
(4P_2)^{-\frac43}\mu^{-\frac23}|\hat{E}|^{\frac43}\leq \alpha^2.
\label{pops24}
\end{equation}
Thus either $\alpha^2\geq \alpha_2^2$ or \eqref{pops24} holds. This shows
\eqref{pops21}.
\hfill{$\Box$}
\begin{lemma}
Corresponding to $0<P<\infty$, $0<K_1<\infty$
there exists $0<\Lambda<\infty$ with the following property. Given $\mu_*$
with $0<\mu_*<1$, define $f_{\mu_*}\,:\, {\mathbb R}\rightarrow {\mathbb R}^+$
by
\begin{equation}
f_{\mu_*}(t)=\left\{\begin{array}{ll}
P\min\{1,\, |t|^{\frac{k-2}{k}},\, \mu_*^{-\frac23} |t|^{\frac43}\} & \mbox{for }\, k>2, \\
P\min\{1,\,\mu_*^{-\frac23}|t|^{\frac43}\} & \mbox{for }\, k\leq 2.
\end{array}
\right.
\label{pops25}
\end{equation}
Suppose $\Phi\,:\, [0,\Lambda]\rightarrow {\mathbb R}$ satisfies
\begin{equation}
\Phi(t_2)-\Phi(t_1)\leq -\int_{t_1}^{t_2}f_{\mu_*}(\Phi(s))\, ds,
\hspace{.5cm}0\leq\forall t_1<\forall t_2\leq \Lambda.
\label{pops26}
\end{equation}
(1) Suppose $\Phi(0)\leq \omega_k$. Then we have
\begin{equation}
\Phi(\Lambda)\leq K_1\mu_*^2.
\label{pops27}
\end{equation}
(2) Suppose $\Phi(0)<-K_1\mu_*^2$. Then we have
\begin{equation}
\Phi(\Lambda)< -\omega_k.
\label{pops27.2}
\end{equation}
\label{popslem3}
\end{lemma}
\begin{rem}
The point of Lemma \ref{popslem3} is that $\Lambda$ does not
depend on $\mu_*$.
\end{rem}
{\it Proof}. Consider the case $k>2$ and (1). 
We set 
\begin{equation}
\Lambda=\frac{\omega_k-1}{P}+\frac{k}{2P}+\frac{3}{P K_1^{\frac13}}.
\label{pops27.5}
\end{equation}
Since $\Phi$
is monotone decreasing, we may assume $\Phi>0$ on $[0,\Lambda]$ for the
conclusion. We have from \eqref{pops25}
\begin{equation}
f_{\mu_*}(\Phi(s))=\left\{\begin{array}{ll} 
P & \mbox{for }\, 1\leq \Phi(s), \\
P\Phi(s)^{\frac{k-2}{k}} & \mbox{for }\, \mu_*^{\frac{2k}{k+6}}\leq \Phi(s)\leq 1, \\
P\mu_* ^{-\frac23} \Phi(s)^{\frac43} & \mbox{for }\, 0<\Phi(s)\leq \mu_*^{\frac{2k}{k+6}}.
\end{array}
\right.
\label{pops28}
\end{equation}
Suppose that $1\leq \Phi\leq \omega_k$ on $[0,\tilde{t}]$. Then \eqref{pops26} and \eqref{pops28}
show that $\tilde{t}\leq (\omega_k-1)/P$. 
Suppose that $1\geq \Phi\geq \mu_*^{2k/(k+6)}$ on $[\tilde{t},t_*]$. 
From \eqref{pops26} and \eqref{pops28}, we have $\Phi'\leq -P\Phi^{(k-2)/k}$ a.e$.$ on $[\tilde{t},t_*]$. Integrating
over the interval gives
\begin{equation}
\Phi(t_*)^{\frac{2}{k}}\leq \Phi(\tilde{t})^{\frac{2}{k}}-\frac{2P}{k}( t_*-\tilde{t}).
\label{pops29}
\end{equation}
Since $\Phi(\tilde{t})\leq 1$ and $\Phi>0$, \eqref{pops29} shows that $t_*-\tilde{t}< k/(2P)$. 
Thus there exists $\hat{t}\in [0,k/(2P)+(\omega_k-1)/P)$ such that $\Phi(\hat{t})<\mu_*^{2k/(k+6)}$.
From \eqref{pops26} and \eqref{pops28}, we have $\Phi'\leq -P\mu_*^{-2/3}
\Phi^{4/3}$ a.e$.$ on $[\hat{t},\Lambda]$. Integrating over the interval gives
\begin{equation}
\Phi(\Lambda)^{-\frac13}\geq \Phi(\hat{t})^{-\frac13}+\frac13 P \mu_*^{-\frac23}(\Lambda
-\hat{t}).
\label{pops30}
\end{equation}
Thus we obtain from \eqref{pops30}  that $\Phi(\Lambda)\leq 27\mu_*^2 P^{-3}
(\Lambda-\hat{t})^{-3}$. By \eqref{pops27.5}, this shows \eqref{pops27} for $k>2$.
For (2), define $\tilde{\Phi}(t)=-\Phi(\Lambda-t)$ for $t\in [0,\Lambda]$. Then one 
can check that $\tilde{\Phi}$ satisfies \eqref{pops26}. By the first part of the proof, 
if $\tilde{\Phi}(0)\leq \omega_k$, then $\tilde{\Phi}(\Lambda)\leq K_1\mu_*^2$. 
This shows \eqref{pops27.2}. 
The proof for $k\leq 2$ is similar and is omitted. 
\hfill{$\Box$}
\begin{thm} Corresponding to $1\leq E_1<\infty$ and $0<\nu<1$ there exist
$0<\varepsilon_1<1$, $0<\Lambda<\infty$ and $0<K<\infty$ with the following
property. 
For $U=C(T,1)$ suppose $\{V_t\}_{0\leq t\leq 2\Lambda+3}$ and 
$\{u(\cdot,t)\}_{0\leq t\leq 2\Lambda+3}$ satisfy (A1)-(A4). 
Assume
\begin{equation}
0<\exists C<\infty \,\,:\,\, {\rm spt}\, \|V_t\|\subset
C(T,1)\cap \{x\,:\, |T^{\perp}(x)|<C\},\hspace{.5cm}\,0\leq \forall t\leq 2\Lambda+3,
\label{pops30.5}
\end{equation} 
\begin{equation}
\mu_*^2:=\sup_{0\leq t\leq 2\Lambda+3}
\int_{C(T,1)}|T^{\perp}(x)|^2 \, d\|V_t\|(x)\leq \varepsilon_1^2,
\label{pops31}
\end{equation}
\begin{equation}
0\leq \exists \hat{t}_1\leq 1\,\,:\,\, \|V_{\hat{t}_1}\|(\phi_{T}^2)
\leq {\bf c}(2-\nu),
\label{pops31.1}
\end{equation}
\begin{equation}
2\Lambda+2\leq \exists \hat{t}_2\leq 2\Lambda+3\,\,:\,\, \|V_{\hat{t}_2}\|(\phi_{T}^2)
\geq {\bf c}\nu,
\label{pops31.2}
\end{equation} 
\begin{equation}
C(u):=\int_0^{2\Lambda+3}\int_{C(T,1)}2|u|^2\phi_{T}^2\, d\|V_t\|dt\leq 
\varepsilon_1^2.
\label{pops31.3}
\end{equation}
Then, we have
\begin{equation}
\sup_{1+\Lambda\leq t\leq 2+\Lambda}\left|\|V_t\|(\phi_{T}^2)-{\bf c}\right|\leq K (\mu_*^2+C(u))
\label{pops31.4}
\end{equation}
and 
\begin{equation}
\int_{1+\Lambda}^{2+\Lambda}\int_{C(T,1)}|h(V_t,\cdot)|^2\phi_{T}^2\, d\|V_t\|dt\leq 
12K(\mu_*^2+C(u)).
\label{popsmean}
\end{equation}
\label{poptheorem}
\end{thm}
\begin{rem}
Imprecise but helpful picture to keep in mind is to regard the moving
varifold $\|V_t\|$ as a hypersurface graph (restricting to $n=k+1$ case)
$y=f(x,t)$ with moving domain $B(t)\subset B^k_1$, with some extra
pieces around. Assume $C(u)=0$, since it is a perturbative term,
and think $\phi_T$ as a characteristic function of $C(T,1)$. 
Assumption 
\eqref{pops31} says that spatial $L^2$-norm $\int_{B(t)}|f(\cdot,t)|^2\,
dx$ is kept small uniformly for the whole time interval. \eqref{pops31.1}
says that near the initial time, $\int_{B(t)}\sqrt{1+|\nabla f|^2}\,dx
-\omega_k$ plus area of some other pieces is strictly less than $\omega_k$, 
preventing the situation where two flat sheets remaining still.
\eqref{pops31.2} says that near the final time, there is some fixed amount
of mass inside the cylinder. This prevents the complete loss of varifold in the
intermediate time. The conclusion \eqref{pops31.4} gives the energy bound
in the following sense. By Taylor expansion, we may postulate that 
\begin{equation*}
\int_{B(t)}\sqrt{1+|\nabla f|^2}\, dx-\omega_k \approx \int_{B(t)}\frac{|\nabla f|^2}{2}\, dx
+ ({\mathcal L}^k(B(t))-\omega_k).
\end{equation*}
Thus, roughly speaking, 
\eqref{pops31.4} gives a uniform in time Dirichlet energy bound of the
graph in terms of $L^2$-norm of $f$. Since the mean curvature is analogous to 
the Laplacian of $f$, \eqref{popsmean} gives a second derivative $L^2$-norm
(space-time) bound in terms of the $L^2$-norm of $f$.
In the next section we prove that $\mu_*^2$ can be estimated in terms of 
the space-time $L^2$-norm of $f$, $\int dt\int_{B(t)}|f|^2\, dx$. 
\end{rem}
{\it Proof}. 
We prove the case $k>2$. The proof for $k\leq 2$ is similar. 
We set 
\begin{equation}
K_2:=80\sup (5|\nabla\phi_T|^4\phi_T^{-2}+|\nabla|\nabla\phi_T||^2),
\label{pops31.5}
\end{equation}
which is a constant depending only on $k$. 
By Proposition \ref{popslem1} and Corollary \ref{popslem2}, 
corresponding to $E_1$ 
of assumption (A2) and $\nu$ there replaced by $\nu/2$, 
we choose and fix $\alpha_2$, $\mu_2$ and $P_2$. 
We then set 
\begin{equation}
P:=\frac{1}{4\cdot 2^{4/3}}\min\{\alpha_2^2,\, (4P_2)^{-\frac{k-2}{k}},\, (4P_2)^{-\frac43}\}.
\label{pops38}
\end{equation}
With this $P$ and 
\begin{equation}
K_1:=2P_2,
\label{pops39}
\end{equation} 
we use Lemma \ref{popslem3} to fix
$\Lambda$. We then choose 
\begin{equation}
\varepsilon_1^2: = \min\left\{K_2^{-1}(3+2\Lambda)^{-1}\frac{\nu{\bf c}}{4},\, 
\frac{\nu{\bf c}}{4},\, \mu_2^2\right\}
\label{pops40}
\end{equation}
and
\begin{equation}
K:=\max\{2,\, 4P_2,\, 2K_2(2\Lambda+3)\}.
\label{pops40.1}
\end{equation}
Overall, note that all those constants are fixed depending
only on $n$, $k$, $E_1$ and $\nu$. We prove the claim with the constants defined above.

Define
\begin{equation}
E(t):=\|V_t\|(\phi_T^2)- {\bf c}-\int_0^t\int_{C(T,1)}2|u|^2\phi_T^2\, d\|V_s\|ds-K_2 \mu_*^2 t.
\label{pops32}
\end{equation}
We first prove that 
\begin{equation}
E(t_2)-E(t_1)\leq -\frac14 \int_{t_1}^{t_2}\int_{C(T,1)}|h(V_t,\cdot)|^2\phi_T^2\, d\|V_t\|dt
,\,\,0\leq \forall t_1<\forall t_2\leq 2\Lambda+3.
\label{pops33}
\end{equation}
\label{popslem4}
The assumption \eqref{pops30.5} and (A2) show
that $\|V_t\|({\rm spt}\, \phi_T)$ is uniformly bounded in $t$. 
By (A1), (A3) and \eqref{hhh}, for a.e$.$ $t$, we have $V_t\in {\bf IV}_k(U)$,
$h(V_t,\cdot)\in L^2(\|V_t\|)$ and $u(\cdot,\,t)\in L^2(\|V_t\|)$. At such time $t$, by \eqref{Bdef}
and the perpendicularity of mean curvature \eqref{fvf2}, we have (omitting $t$ dependence for simplicity)
\begin{equation}
{\mathcal B}(V,u,\phi_T^2)\leq \int_U -|h|^2\phi_T^2+\phi_T^2|h||u|+|u^{\perp}\cdot\nabla\phi_T^2|
+(\nabla\phi_T^2)^{\perp}\cdot h\, d\|V\|.
\label{pops34}
\end{equation}
Since $\phi_T$ depends only on $|T(x)|$, $\nabla\phi_T\in T$. Thus the last term of
\eqref{pops34} is
\begin{equation}
\begin{split}
\int_{G_k(U)}2\phi_T S^{\perp}(\nabla\phi_T)&\cdot h\, dV
=\int_{G_k(U)}2\phi_T (T-S)(\nabla\phi_T)\cdot h\, dV \\
&\leq \frac14 \int_{U} |h|^2 \phi_T^2\, d\|V\|+4\int_{G_k(U)}|\nabla\phi_T|^2\|S-T\|^2\,
dV.
\end{split}
\label{pops35}
\end{equation}
Similarly,
\begin{equation}
\begin{split}
\int_U |u^{\perp}\cdot\nabla\phi_T^2|\, d\|V\|&=\int_{G_k(U)}2\phi_T |(T-S)(u)\cdot\nabla
\phi_T|\, dV \\
&\leq \int_U \phi_T^2 |u|^2\, d\|V\|+\int_{G_k(U)}|\nabla\phi_T|^2\|S-T\|^2\, dV.
\end{split}
\label{pops35.5}
\end{equation}
To estimate the last term of \eqref{pops35} and \eqref{pops35.5} we slightly modify 
the proof of Lemma \ref{tiltexlem} and obtain
\begin{equation}
\begin{split}
\int_{G_k(U)}|\nabla\phi_T|^2\|S-T\|^2\, dV&\leq 
4\left(\int_U \phi_T^2 |h|^2\, d\|V\|\int_{U} |T^{\perp}(x)|^2 |\nabla\phi_T|^4
\phi_T^{-2}\, d\|V\|\right)^{\frac12} \\ &+16 \int_{U} |T^{\perp}(x)|^2|\nabla
|\nabla\phi_T||^2\, d\|V\|.
\end{split}
\label{pops36}
\end{equation}
Combining \eqref{pops34}-\eqref{pops36}, we obtain
\begin{equation}
{\mathcal B}(V,u,\phi_T^2)\leq \int_{U} -\frac{|h|^2\phi_T^2}{4}+2|u|^2\phi_T^2 
+ 80 |T^{\perp}(x)|^2(5
|\nabla\phi_T|^4\phi_T^{-2}+|\nabla|\nabla\phi_T||^2)\,
d\|V\|.
\label{pops37}
\end{equation}
Since $V_t$ satisfies (A4), integrating \eqref{pops37} over $[t_1,t_2]$
and using \eqref{pops31} and \eqref{pops31.5}, we obtain \eqref{pops33}. 
Next, by \eqref{pops31.1} and \eqref{pops33}, we have for any $t\in [\hat{t}_1,
2\Lambda+3]$
\begin{equation}
E(t)\leq E(\hat{t}_1)\leq {\bf c}(1-\nu).
\label{pops37.1}
\end{equation}
Due to \eqref{pops31}, \eqref{pops31.3} and \eqref{pops40}, \eqref{pops37.1}
shows
\begin{equation}
\|V_t\|(\phi_T^2)-{\bf c}\leq {\bf c}(1-\nu)+2\times \frac{\nu{\bf c}}{4}={\bf c}(1-\nu/2),
\hspace{.5cm} \hat{t}_1\leq \forall t\leq 2\Lambda+3.
\label{pops37.2}
\end{equation}
Similarly, due to \eqref{pops31.2} and \eqref{pops33}, we have
\begin{equation}
E(t)\geq E(\hat{t}_2)\geq {\bf c}(\nu-1)-2\times\frac{\nu{\bf c}}{4}={\bf c}(-1+\nu/2)
\label{37.25}
\end{equation}
and
\begin{equation}
\|V_t\|(\phi_T^2)-{\bf c}\geq {\bf c}(-1+\nu/2)
\label{pops37.3}
\end{equation}
for all $t\in [0,\hat{t}_2]$. 
To obtain a contradiction assume that there exists $t_*\in [\Lambda+1,\Lambda+2]$
with
\begin{equation}
\|V_{t_*}\|(\phi_T^2)-{\bf c}>K(\mu_*^2+C(u)).
\label{pops41}
\end{equation}
By \eqref{pops32}, \eqref{pops40.1} and \eqref{pops33}, \eqref{pops41} implies that for all $t\in [0,t_*]$,
\begin{equation}
\begin{split}
\|V_t\|(\phi_T^2)-{\bf c} & \geq E(t)\geq E(t_*)>K(\mu_*^2+C(u))-C(u)-K_2\mu_*^2 (\Lambda+3) \\
&\geq \frac{K}{2}
\mu_*^2\geq 2P_2 \mu_*^2.
\end{split}
\label{pops42}
\end{equation}
By \eqref{pops37.2} and \eqref{pops42}, we may apply Corollary \ref{popslem2}
for a.e$.$ $t\in [\hat{t}_1, t_*]$ to obtain 
\begin{equation}
\int_{C(T,1)}|h(V_t,\cdot)|^2 \phi_T^2\, d\|V_t\|\geq 4P\min \{1,\, E(t)^{\frac{k-2}{k}},\,
\mu_*^{-\frac23}E(t)^{\frac43}\},
\label{pops43}
\end{equation}
where we also used \eqref{pops38}, \eqref{pops31} and $\|V_t\|(\phi_T^2)-{\bf c}\geq E(t)$.
We then use Lemma \ref{popslem3}. With $\Phi(t)=E(t+\hat{t}_1)$, we have $\Phi(0)\leq 
{\bf c}(1-\nu)<\omega_k$ and \eqref{pops26} is satisfied on $[0,\Lambda]$
due to \eqref{pops33} and \eqref{pops43}. We thus conclude that $E(\Lambda+\hat{t}_1)
\leq 2P_2\mu_*^2$. On the other hand, from \eqref{pops42} and $\Lambda+\hat{t}_1\in
[0,t_*]$, we have $E(\Lambda+\hat{t}_1)>2P_2\mu_*^2$, which is a contradiction. 
Similarly, for a contradiction, assume that there exists $t_*\in [\Lambda+1,\Lambda+2]$
with 
\begin{equation}
 \|V_{t_*}\|(\phi_T^2)-{\bf c}\leq -K(\mu_*^2+C(u)).
\label{pops44}
\end{equation}
Similar computations to \eqref{pops42} using \eqref{pops32}, \eqref{pops33}, \eqref{pops40.1}
and \eqref{pops44} show 
\begin{equation}
\begin{split}
\|V_t\|(\phi_T^2)-{\bf c}& \leq E(t)+C(u)+K_2\mu_*^2 (2\Lambda+3)
\leq E(t_*)+C(u)+K_2\mu_*^2 (2\Lambda+3) \\
&\leq \|V_{t_*}\|(\phi_T^2)-{\bf c}+C(u)+K_2\mu_*^2(2\Lambda+3) \leq 
-\frac{K}{2}(\mu_*^2+C(u))
\leq -2P_2 \mu_*^2
\end{split}
\label{pops45}
\end{equation}
for $t\in [t_*,2\Lambda+3]$. Thus, by  \eqref{pops37.3}, \eqref{pops45} and
\eqref{pops38}, Corollary \ref{popslem2}
shows
\begin{equation}
\int_{C(T,1)}|h(V_t,\cdot)|^2\phi_T^2\, d\|V_t\|\geq 4\cdot 2^{4/3}P
\min\{1,\, | \|V_t\|(\phi_T^2)-{\bf c}|^{\frac{k-2}{k}},\, \mu_*^{-\frac23}|\|V_t\|(\phi_T^2)-{\bf c}|^{\frac43}
\}
\label{pops46}
\end{equation}
for a.e$.$ $t\in [t_*,\hat{t}_2]$. Since $|\|V_t\|(\phi_T^2)-{\bf c}|\geq K(\mu_*^2+C(u))/2$ as in \eqref{pops45},
we have by \eqref{pops40.1}
\begin{equation}
|E(t)|\leq |\|V_t\|(\phi_T^2)-{\bf c}|+C(u)+K_2\mu_*^2(2\Lambda+3)\leq 2|\|V_t\|(\phi_T^2)-{\bf c}|.
\label{pops47}
\end{equation}
Combining \eqref{pops46} and \eqref{pops47}, we obtain
\begin{equation}
\int_{C(T,1)}|h(V_t,\cdot)|^2\phi_T^2\, d\|V_t\|\geq 4P\min\{1,\, |E(t)|^{\frac{k-2}{k}},\,
\mu_*^{-\frac23}|E(t)|^{\frac43}\}.
\label{pops48}
\end{equation}
Then applying Lemma \ref{popslem3}, we obtain $E(t_*+\Lambda)<-\omega_k$. Since $t_*+\Lambda
\in [0,2+2\Lambda]$, this is a contradiction to \eqref{37.25}.
Lastly one observes that \eqref{pops33}, \eqref{pops32}
and \eqref{pops31.4} show \eqref{popsmean}.
\hfill{$\Box$}

For parabolic Lipschitz approximation, we also need the following estimate.
\begin{cor} 
Under the same assumptions of Theorem \ref{poptheorem},
there exists a constant $\tilde{K}$ depending only on $E_1$ and $\nu$ such that 
\begin{equation}
\int_{1+\Lambda}^{2+\Lambda}\left|2^k\|V_t\|(\phi_{T,1/2}^2 )-{\bf c}\right|
\, dt\leq \tilde{K}(\mu_*^2+C(u)).
\label{anopop}
\end{equation}
\label{popscor}
\end{cor}
{\it Proof}. 
By Lemma \ref{tiltexlem}, \eqref{pops31} and \eqref{popsmean}, 
we have
\begin{equation}
\int_{1+\Lambda}^{2+\Lambda}\int_{G_k(C(T,1))}\|S-T\|^2\phi_T^2 \, dV_t(\cdot,S)dt
\leq 4(12K(\mu_*^2 +C(u)))^{1/2}\mu_*+c(k)\mu_*^2.
\label{anopop1}
\end{equation}
Then the difference of $\|V\|(\phi_T^2)$ and $2^k\|V\|(\phi_{T,1/2}^2)$
can be estimated using Theorem \ref{cylgrolem}, \eqref{anopop1} and
\eqref{popsmean}. This gives \eqref{anopop}.
\hfill{$\Box$}
\label{enees}
\section{Parabolic monotonicity formula and $L^2$-$L^{\infty}$ estimate}
In this section we prove two important tools which are used
in the subsequent sections. They are both derived by modifying
the monotonicity formula due to Huisken \cite{Huisken}. 
Note that Brakke did not utilize these
powerful tools since it was not known at the time. 
He instead used a certain vanishing theorem 
\cite[6.3. `Clearing out']{Brakke} to obtain a uniform $L^{\infty}$ estimate
of the height of varifold.
\subsection{Local monotonicity formula}
In this subsection, let $\eta\in C^{\infty}_c(B_{14/15})$
be a non-negative radially symmetric function such that $\eta\equiv
1$ on $B_{13/15}$, $|\nabla\eta|\leq 15$ and
$\|\nabla^2\eta\|\leq c(n)$. Then define 
${\hat \rho}_{(y,s)}(x,t)=\eta(x)\rho_{(y,s)}(x,t)$, where $\rho_{(y,s)}(x,t)$ is defined
as in \eqref{hkneldef}. 
\begin{lemma}
There exists $c_1=c_1(n,k)$ with the following property.
For $0<s<t<\infty$, $x\in B_1$, $y\in B_{4/5}$ and $S\in {\bf G}(n,k)$, we have
\begin{equation}
\left|\frac{\partial{\hat{\rho}}_{(y,s)}(x,t)}{\partial t}
+S\cdot \nabla^2 {\hat{\rho}}_{(y,s)}(x,t)
+\frac{|S^{\perp}(\nabla{\hat{\rho}}_{(y,s)}(x,t))|^2}{{\hat{\rho}}_{(y,s)}(x,t)}\right|\leq c_1.
\label{monest}
\end{equation}
\end{lemma}
{\it Proof}. Here we write $\rho$ for $\rho_{(y,s)}(x,t)$ and 
$\hat{\rho}$ for ${\hat{\rho}}_{(y,s)}(x,t)$. First assume that
$x\in B_{13/15}$ so that the derivatives of $\eta$ are 0 and
${\hat{\rho}}=\rho$ in the neighborhood.
We have
\begin{equation*}\frac{\partial{\hat{\rho}}}{\partial t}
=\left(\frac{k}{2(s-t)}-\frac{|x-y|^2}{4(s-t)^2}\right){\hat{\rho}},
\hspace{.5cm} \nabla{\hat{\rho}}=-\frac{x-y}{2(s-t)}{\hat{\rho}}, 
\end{equation*}
\begin{equation*}
\nabla^2{\hat{\rho}}=-\frac{\hat{\rho}}{2(s-t)}I+
\frac{\hat{\rho}}{4(s-t)^2}(x-y)\otimes(x-y).
\end{equation*}
Since $S\cdot I=k$, $S\cdot(x-y)\otimes(x-y)=|S(x-y)|^2$ and
$|S^{\perp}(x-y)|^2+|S(x-y)|^2=|x-y|^2$, 
$$\frac{\partial{\hat{\rho}}}{\partial t}+S\cdot\nabla^2{\hat{\rho}}
+\frac{|S^{\perp}(\nabla{\hat{\rho}})|^2}{\hat{\rho}}=0.$$
When $x\in B_1\setminus B_{13/15}$, $|x-y|>\frac{1}{15}$. 
The similar computation in this case produces extra terms which involves
$\nabla\eta$ and $\nabla^2\eta$. Since they are 
multiplied by $\exp\left(-\frac{|x-y|^2}{4(s-t)}\right)<\exp\left(-
\frac{1}{900(s-t)}\right)$, there exists some constant $c_1
=c_1(n,k)$ satisfying the desired estimate. 
\hfill{$\Box$}
\begin{prop}
For $U=B_R$, $0<\Lambda\leq \infty$ and $1\leq E_1<\infty$
suppose that $\{V_t\}_{0\leq t<\Lambda}$
and $\{u(\cdot,t)\}_{0\leq t<\Lambda}$ satisfy (A1)-(A4). Define
$\hat{\rho}_{(y,s)}(x,t)=\eta(x/R) \rho_{(y,s)}(x,t)$. 
For $0<t_1<t_2<\Lambda$, $t_2<s$ and $y\in B_{4R/5}$, we have
\begin{equation}
\begin{split}
&\int_{B_R}{\hat{\rho}}_{(y,s)}(\cdot,t_2)\, d\|V_{t_2}\|-
\int_{B_R}{\hat{\rho}}_{(y,s)}(\cdot,t_1)\, d\|V_{t_1}\|
\\&\leq c_2 \|u\|^2_{L^{p,q}(B_R\times(t_1,t_2))}E_1^{1-\frac{2}{p}}
(t_2-t_1)^{\varsigma}+c_1 \o_k E_1 R^{-2}(t_2-t_1).
\end{split}
\label{monoest}
\end{equation}
Here $c_2=c_2(k,p,q)$ and $\varsigma=1-k/p-2/q$. 
\label{cam1}
\end{prop}
{\it Proof}. After a change of coordinates we may assume $R=1$ without 
loss of generality. 
We compute for a.e$.$ $t\in (t_1,t_2)$ (when $V_t$ is integral and $h(V_t,\cdot)
\in L^2(\|V_t\|)$ by (A1) and \eqref{hhh})
 and $t<s$ using the definition \eqref{Bdef} as
\begin{equation}
\begin{split}
&{\mathcal B}(V_t,u(\cdot,t),{\hat{\rho}}(\cdot,t))=\int_{B_1}(-h{\hat{\rho}}+\nabla{\hat{\rho}})\cdot(h+u^{\perp})\, d\|V_t\| \\ &
=\int_{B_1}-|h|^2{\hat{\rho}}+2\nabla{\hat{\rho}}\cdot h+u^{\perp}\cdot(-h{\hat{\rho}}+\nabla{\hat{\rho}})-\nabla{\hat{\rho}}\cdot h\, d\|V_t\|.
\end{split}
\label{Bcom1}
\end{equation}
We then complete the square by adding and subtracting 
$|(\nabla{\hat{\rho}})^{\perp}|^2/{\hat{\rho}}$. Due to the
perpendicularity of mean curvature \eqref{fvf2}, we have for $\|V_t\|$ a.e$.$ $h\cdot\nabla{\hat{\rho}}
=h\cdot(\nabla{\hat{\rho}})^{\perp}$. 
We then continue the computation of \eqref{Bcom1} as 
\begin{equation}
=\int_{B_1}-{\hat{\rho}}\left|h-\frac{(\nabla{\hat{\rho}})^{\perp}
}{\hat{\rho}}\right|^2+\frac{|(\nabla{\hat{\rho}})^{\perp}|^2}{
\hat{\rho}}+u\cdot(-h{\hat{\rho}}+(\nabla{\hat{\rho}})^{\perp})+S\cdot
\nabla^2{\hat{\rho}}\, d\|V_t\|.
\label{Bcom2}
\end{equation}
By the Cauchy-Schwarz inequality applied to the third term of \eqref{Bcom2}, 
we obtain
\begin{equation}
\leq \int_{B_1}\frac{|(\nabla{\hat{\rho}})^{\perp}|^2}{\hat{\rho}}+
|u|^2\hat{\rho}+S\cdot\nabla^2\hat{\rho}\, d\|V_t\|.
\label{Bcom3}
\end{equation}
Substituting \eqref{Bcom3} into \eqref{maineq} and using \eqref{dbd} and
\eqref{monest}, we 
obtain
\begin{equation}
\left.\int_{B_1}\hat{\rho}(\cdot,t)\, d\|V_t\|\right|_{t=t_1}^{t_2}
\leq \int_{t_1}^{t_2}\left( \int_{B_1}|u|^2\hat{\rho}\, d\|V_t\|+
c_1 \o_k E_1\right)\, dt.
\label{Bcom4}
\end{equation}
For the first term of \eqref{Bcom4}, by the H\"{o}lder inequality, we have
\begin{equation}
\begin{split}
&\int_{t_1}^{t_2}\int_{B_1}|u|^2\hat{\rho}\, d\|V_t\|dt\leq \int_{t_1}^{t_2}
\left(\int_{B_1}|u|^p\hat{\rho}\, d\|V_t\|\right)^{\frac{2}{p}}
\left(\int_{B_1}\hat{\rho}\, d\|V_t\|\right)^{1-\frac{2}{p}}\, dt\\
&\leq E_1^{1-\frac{2}{p}}\int_{t_1}^{t_2}\frac{1}{(4\pi(s-t))^{\frac{k}{p}}}
\left(\int_{B_1}|u|^p\, d\|V_t\|\right)^{\frac{2}{p}}\, dt .
\end{split}
\label{Bcom5}
\end{equation}
We used \eqref{dbd} in the last line of \eqref{Bcom5}.
Again using the H\"{o}lder inequality, \eqref{Bcom5} may be estimated as
\begin{equation}
\begin{split}
&\leq \left(\int_{t_1}^{t_2}\left(\int_{B_1}|u|^p\, d\|V_t\|
\right)^{\frac{q}{p}}\, dt\right)^{\frac{2}{q}}\left(\int_{t_1}^{t_2}
(4\pi(s-t))^{-\frac{kq}{p(q-2)}}\, dt
\right)^{\frac{q-2}{q}} E_1^{1-\frac{2}{p}}\\
&\leq \|u\|^2_{L^{p,q}(B_1\times(t_1,t_2))}(4\pi)^{-\frac{k}{p}}
E_1^{1-\frac{2}{p}}\left(\left.-\left(1-\frac{kq}{p(q-2)}\right)^{-1}
(s-t)^{1-\frac{kq}{p(q-2)}}\right|_{t=t_1}^{t_2}\right)^{\frac{q-2}{q}}.
\end{split}
\label{Bcom6}
\end{equation}
Here we used \eqref{expcond} for the integrability with respect to $t$. 
Since $(s-t_1)^{1-\frac{kq}{p(q-2)}}-(s-t_2)^{1-\frac{kq}{p(q-2)}}\leq (t_2-t_1)^{1-\frac{kq}{p(q-2)}}$,
by setting $c_2=c_2(k,p,q)$, we obtain \eqref{monoest}.
\hfill$\Box$
\begin{cor}
Corresponding to $1\leq E_1<\infty$, $p$ and $q$, there exist $0<c_{18}<1$ and $0<c_{19}<1$ with the following property.
Suppose $\{V_t\}_{0\leq t\leq \Lambda}$ and $\{u(\cdot,t)\}_{0\leq t\leq \Lambda}$ satisfy (A1)-(A4) and assume
$B_R(a)\times (t,t+c_{18}R^2)\subset\subset U\times(0, \Lambda)$. In addition assume
\begin{equation}
R^{\varsigma} \|u\|_{L^{p,q}(B_R(a)\times(t,t+c_{18}R^2))}\leq 1,\,\,\, \|V_t\|(B_{14R/15}(a))\leq c_{19}R^k.
\label{cing1}
\end{equation}
Then we have $\|V_{t+c_{18}R^2}\|(B_{4R/5}(a))=0$.
\label{clearing}
\end{cor}
{\it Proof}. After a change of variables, we may assume $R=1$, $t=0$ and $a=0$. First
By Proposition \ref{cam1} with $t_2=c_{18}$ to be chosen, $s=t_2+\epsilon$ and $t_1=0$, we have for 
any $y\in B_{4/5}$
\begin{equation}
\left.\int_{B_1}\hat{\rho}_{(y,t_2+\epsilon)}(\cdot, t)\, d\|V_{t}\| \right|_{t=0}^{t_2}\leq c_2\|u\|_{L^{p,q}(B_1\times(0,t_2))}^2 E_1^{1-\frac{2}{p}}t_2^{\varsigma}+c_1\omega_k
E_1 t_2.
\label{cing2}
\end{equation}
Here $\epsilon>0$ is arbitrary. 
Since $\hat{\rho}$ has support in $B_{14/15}$ and $\eta\leq 1$, and by \eqref{cing1} and \eqref{cing2}, we have
\begin{equation}
\int_{B_1}\hat{\rho}_{(y,t_2+\epsilon)}(\cdot,t_2)\, d\|V_{t_2}\|\leq (4\pi t_2)^{-k/2}c_{19}+c_2 E_1^{1-\frac{2}{p}}
t_2^{\varsigma}+c_1\omega_k
E_1 t_2.
\label{cing3}
\end{equation}
Suppose that $V_{t_2}$ is unit density, and assume for a contradiction that 
$\|V_{t_2}\|(B_{4/5})>0$. Then there exists some $y\in B_{4/5}$ such that $\Theta^k (\|V_{t_2}\|,y)=1$ 
and one can show that 
\begin{equation}
\lim_{\epsilon\rightarrow 0}\int_{B_1}\hat{\rho}_{(y,t_2+\epsilon)}(\cdot, t_2)\, d\|V_{t_2}\|=1.
\label{cing4}
\end{equation}
Combining \eqref{cing3} and \eqref{cing4}, if we choose $t_2=c_{18}$ appropriately small and then 
choose $c_{19}$ small depending only on $E_1$, $p,\, q$, we obtain a contradiction. If $V_{t_2}$ 
is not unit density, we may find a sequence $\lim_{m\rightarrow\infty}t_m= t_2$ with $t_m<t_2$ such that $V_{t_m}$
is unit density by (A1). Then the previous argument shows that $\|V_{t_m}\|(B_{4/5})=0$ for all large $m$. By
\eqref{maineq}, we may conclude that $\|V_{t_2}\|(B_{4/5})=0$.
\hfill{$\Box$}
\label{monosubsec}
\subsection{$L^2$-$L^{\infty}$ estimate}
Proposition \ref{el2elinf} gives $L^{\infty}$ bound of the height of varifold away
from $T$ in terms of space-time $L^2$-height plus some perturbative 
term. This has been already observed by Ecker \cite{Ecker} in a smooth setting for
mean curvature flow while we are unaware of other applications in the setting of
geometric measure theory. The result bridges space-time $L^2$-height smallness
and time-uniform $L^2$-height smallness, the latter being a necessary condition
for Theorem \ref{poptheorem}. Also $L^{\infty}$ estimate makes a good use in the
blow-up analysis of Section \ref{holder}. 
\begin{prop} For $U=B_R$, $R^2\leq \Lambda< \infty$ and $1\leq E_1<\infty$
suppose that $\{V_t\}_{0\leq t<\Lambda}$
and $\{u(\cdot,t)\}_{0\leq t<\Lambda}$ satisfy (A1)-(A4). 
Then there exists a constant
$c_3=c_3(n,k)$ with the 
following property. For $T\in {\bf G}(n,k)$ set 
\begin{equation}
\mu^2:=\frac{c_3}{R^{k+4}}\int_0^{\Lambda}\int_{B_R}
|T^{\perp}(x)|^2\, d\|V_t\|dt+ c_2 \|u\|^2_{L^{p,q}(B_R\times (0,\Lambda ))}
E_1^{1-\frac{2}{p}}\Lambda^{\varsigma}(2+\Lambda/R^2).
\label{Kdef}
\end{equation}
Then for all $t\in (R^2/5,\Lambda)$ we have
\begin{equation}
{\rm spt}\, \|V_t\|\cap B_{4R/5}\subset\{|T^{\perp}(x)| \leq  \mu R \}.
\label{l2linf}
\end{equation}
\label{el2elinf}
\end{prop}
{\bf Proof}. Without loss of generality we may assume $R=1$ and redefine
$\Lambda$ as $\Lambda/R^2$. In the proof
let $\eta\in C^{\infty}(B_1\times(0,\Lambda))$ be a non-negative
function with $\eta\equiv 1$ on $B_{13/15}\times [2/15,
\Lambda)$, $\eta\equiv 0$ on $(B_1\times(0,\Lambda))\setminus(B_{14/15}
\times[1/15,\Lambda))$, $0\leq \eta\leq 1$ and 
$|\nabla\eta|,\, \|\nabla^2\eta\|, 
|\frac{\partial\eta}{\partial t}|\leq c(n)$. For $(y,s)\in B_{4/5}
\times (1/5,\infty)$, we use $\phi(x,t)=|T^{\perp}(x)|^2\rho_{(y,s)}(x,t)\eta(x,t)$
in \eqref{maineq}, over the time interval $t_1=0$ and $0<t_2<\Lambda$ with $t_2<s$. 
We then obtain (writing $\rho_{(y,s)}(x,t)$ as $\rho$)
\begin{equation}
\begin{split}
&\left.\int_{B_1}|T^{\perp}(x)|^2\rho\eta\, d\|V_t\|\right|_{t=t_2} \\ 
&\leq
\int_0^{t_2}\int_{B_1}\{-h\rho\eta|T^{\perp}(x)|^2+\nabla(\rho\eta
|T^{\perp}(x)|^2)\}\cdot(h+u^{\perp})+|T^{\perp}(x)|^2
\frac{\partial}{\partial t}(\rho\eta)\, d\|V_t\|dt
\end{split}
\label{l21}
\end{equation}
since $\eta=0$ for $t=0$. The integrand of the right-hand side of \eqref{l21} 
can be computed as follows for a.e$.$ $t$. 
We use the perpendicularity of the mean curvature
vector \eqref{fvf2}.
\begin{equation}
\begin{split}
&-|h|^2\rho\eta|T^{\perp}(x)|^2+(\nabla\rho\cdot h)\eta|T^{\perp}(x)|^2+
\rho\nabla(\eta|T^{\perp}(x)|^2)\cdot h+(-h\rho\eta|T^{\perp}(x)|^2 \\ &+
\eta|T^{\perp}(x)|^2\nabla\rho)\cdot u^{\perp}+\rho\nabla(\eta
|T^{\perp}(x)|^2)\cdot u^{\perp}+|T^{\perp}(x)|^2\frac{\partial}{\partial t}
(\rho\eta) \\
&\leq -\rho\left|h-\frac{(\nabla\rho)^{\perp}}{\rho}\right|^2\eta|T^{\perp}(x)|^2
-(\nabla\rho\cdot h)\eta|T^{\perp}(x)|^2+\frac{|(\nabla\rho)^{\perp}|^2}{\rho}
\eta|T^{\perp}(x)|^2 \\
&+\rho\nabla(\eta|T^{\perp}(x)|^2)\cdot h+\rho\left|h-\frac{(\nabla\rho)^{\perp}}
{\rho}\right|^2\eta|T^{\perp}(x)|^2 
+\rho\eta|T^{\perp}(x)|^2|u|^2 \\
&+\rho\nabla(\eta|T^{\perp}(x)|^2)
\cdot u^{\perp}+|T^{\perp}(x)|^2\frac{\partial}{\partial t}(\rho\eta).
\end{split}
\label{l22}
\end{equation}
Thus we obtain from \eqref{l21} and \eqref{l22} that
\begin{equation}
\begin{split}
& \left.\int_{B_1}|T^{\perp}(x)|^2\rho\eta\, d\|V_t\|\right|_{t=t_2}
\leq \int_0^{t_2}\int_{B_1}-(\nabla\rho\cdot h)\eta|T^{\perp}(x)|^2
+\rho\nabla(\eta|T^{\perp}(x)|^2)\cdot h \\ &+\frac{|(\nabla\rho)^{\perp}|^2}{\rho}
\eta|T^{\perp}(x)|^2 
+ \rho\eta|T^{\perp}(x)|^2|u|^2+\rho\nabla(\eta|T^{\perp}(x)|^2)
\cdot u^{\perp}+|T^{\perp}(x)|^2\frac{\partial}{\partial t}(\rho\eta)\,
d\|V_t\|dt.
\end{split}
\label{l23}
\end{equation}
The first two terms on the right-hand side of \eqref{l23} are 
\begin{equation}
\begin{split}
&\int_0^{t_2}\int_{G_k(B_1)}\nabla(\eta|T^{\perp}(x)|^2\nabla\rho)\cdot S
-\nabla\{\rho\nabla(\eta|T^{\perp}(x)|^2)\}\cdot S\, dV_t(x,S) dt \\
&=\int_0^{t_2}\int_{G_k(B_1)}(\nabla^2\rho\cdot S)\eta|T^{\perp}(x)|^2
-\rho\nabla^2(\eta|T^{\perp}(x)|^2)\cdot S\, dV_t(x,S)dt.
\end{split}
\label{l24}
\end{equation}
Now we use the identity $\frac{\partial\rho}{\partial t}+\frac{|(\nabla
\rho)^{\perp}|^2}{\rho}+S\cdot\nabla^2\rho\equiv 0$ to obtain from 
\eqref{l23} and \eqref{l24} that
\begin{equation}
\begin{split}
& \left.\int_{B_1}|T^{\perp}(x)|^2\rho\eta\, d\|V_t\|\right|_{t=t_2}\leq
\int_0^{t_2}\int_{G_k(B_1)}\Big{\{}-\nabla^2 (\eta|T^{\perp}(x)|^2)\cdot S+
\eta|T^{\perp}(x)|^2 |u|^2 \\
& +|T^{\perp}(x)|^2\nabla\eta \cdot u^{\perp}+\eta\nabla|T^{\perp}(x)|^2
\cdot u^{\perp}+|T^{\perp}(x)|^2 \frac{\partial
\eta}{\partial t}\Big{\}}\rho\, dV_t(x,S)dt.
\end{split}
\label{l25}
\end{equation}
For the first term involving $u$ in \eqref{l25}, we use $|T^{\perp}(x)|\leq 1$ and the similar computation as in the previous subsection then gives
\begin{equation}
\int_0^{t_2}\int_{B_1} \eta|T^{\perp}(x)|^2|u|^2\rho \, d\|V_t\|dt 
\leq c_2(k,p,q)\|u\|^2_{L^{p,q}(B_1\times(0,\Lambda))}E_1^{1-\frac{2}{p}}
\Lambda^{\varsigma}.
\label{l26}
\end{equation}
For the second term involving $u$ in \eqref{l25}, we have (again using $|T^{\perp}
(x)|\leq 1$)
\begin{equation}
\int_0^{t_2}\int_{B_1}|T^{\perp}(x)|^2\nabla\eta\cdot u^{\perp}\rho\,
d\|V_t\|dt\leq \int_0^{t_2}\int_{B_1}|T^{\perp}(x)|^2\frac{|\nabla\eta|^2}{\eta}
\rho+\eta|u|^2\rho\, d\|V_t\|dt.
\label{l265}
\end{equation}
We note that $\eta^{-1}|\nabla\eta|^2\leq c(n)\max\|\nabla^2\eta\|$ and
$\eta^{-1}(x,t)|\nabla\eta(x,t)|^2 \rho_{(y,s)}
(x,t)$ is uniformly bounded for $(y,s)\in B_{4/5}\times(1/5,\infty)$
and $(x,t)\in A:=(B_{14/15}\times[1/15,
\Lambda))\setminus(B_{13/15}\times[2/15,\Lambda))$ since
$|x-y|\geq 1/15$. Outside of $A$, it vanishes. Thus we have some $c(n,k)$ such that
$\eta^{-1}|\nabla\eta|^2\rho\leq c(n,k)$. The second term in the 
integral of \eqref{l265} can be estimated as in \eqref{l26}.
For the third term involving $u$ in \eqref{l25}, 
\begin{equation}
\begin{split}
\int_0^{t_2} &\int_{B_1}\eta(\nabla|T^{\perp}(x)|^2\cdot u^{\perp})\rho\, d\|V_t\|dt
\leq \frac{1}{\Lambda}\int_0^{t_2}\int_{B_1}|T^{\perp}(x)|^2 \eta\rho\, d\|V_t\|dt \\
& +\Lambda \int_0^{t_2}\int_{B_1}\eta|u|^2\rho\, d\|V_t\|dt.
\end{split}
\label{l266}
\end{equation}
The last term of \eqref{l266} can be estimated as in \eqref{l26}. 
For the remaining two terms of \eqref{l25}, we have
\begin{equation}
\begin{split}
-\nabla^2 (\eta|T^{\perp}(x)|^2)&\cdot S+|T^{\perp}(x)|^2\frac{\partial\eta}{
\partial t}\\
&=-(\nabla^2\eta\cdot S)|T^{\perp}(x)|^2-4(\nabla\eta\otimes T^{\perp}(x))\cdot
S-2\eta T^{\perp}\cdot S+|T^{\perp}(x)|^2\frac{\partial\eta}{\partial t}.
\end{split}
\label{l27}
\end{equation}
Note that $T^{\perp}\cdot S\geq 0$ in general (see \eqref{simin2}), and 
\begin{equation}
4|(\nabla\eta
\otimes T^{\perp}(x))\cdot S|\leq 4|T^{\perp}(x)||\nabla\eta|\sqrt{T^{\perp}
\cdot S}\leq 2\eta T^{\perp}\cdot S+2\frac{|\nabla\eta|^2}{\eta}|T^{\perp}(x)|^2.
\label{l275}
\end{equation}
Since $\frac{|\nabla\eta|^2}{\eta}\leq c(n)\|\nabla^2\eta\|$, \eqref{l275} shows 
that the right-hand side of \eqref{l27} can be bounded from above by $c(n)
|T^{\perp}(x)|^2$ on  $A$ and
otherwise \eqref{l27} $\leq 0$. Since $\rho\leq c(k)$ on 
$A$, we have
\begin{equation}
\begin{split}
& \int_0^{t_2}\int_{G_k(B_1)}\Big{\{}-\nabla^2 (\eta|T^{\perp}(x)|^2)\cdot S +|T^{\perp}
(x)|^2\frac{\partial \eta}{\partial t}\Big{\}}\rho\, dV_t(x,S)dt \\ & \leq 
c(n,k)\int_0^{\Lambda} \int_{B_1}|T^{\perp}(x)|^2\, d\|V_t\|dt.
\end{split}
\label{l28}
\end{equation}
We obtain from \eqref{l25}, \eqref{l26}, \eqref{l265}, 
\eqref{l266} and \eqref{l28} 
\begin{equation}
\begin{split}
& \left.\int_{B_1}|T^{\perp}(x)|^2\eta\rho\, d\|V_t\|
\right|_{t=t_2}\leq c_2 \|u\|^2_{L^{p,q}}E_1^{1-\frac{2}{p}}
\Lambda^{\varsigma} (2+\Lambda) 
\\ &+\frac{1}{\Lambda}\int_0^{t_2}
\int_{B_1}|T^{\perp}(x)|^2\eta\rho\, d\|V_t\|dt+c(n,k)
\int_0^{\Lambda}\int_{B_1}|T^{\perp}(x)|^2\, d\|V_t\|dt.
\end{split}
\label{l29}
\end{equation}
holds for all $t_2\in (0,\Lambda)$. It is easy to check that 
the differential inequality $F'(t)\leq C+F(t)/\Lambda$ for $t\in
(0,\Lambda)$ with $F(0)=0$ means $F'(t)\leq C e$
for all $t\in (0,\Lambda)$. 
Hence from \eqref{l29} we obtain
\begin{equation}
\begin{split}
&\left.\int_{B_1}|T^{\perp}(x)|^2\eta\rho\, d\|V_t\|\right|_{
t=t_2} \\
& \leq e\left(c_2\|u\|^2_{L^{p,q}} E_1^{1-\frac{2}{p}}
\Lambda^{\varsigma}(2+\Lambda)+c(n,k)
\int_0^{\Lambda}\int_{B_1}|T^{\perp}(x)|^2\, d\|V_t\|dt\right).
\end{split}
\label{l210}
\end{equation}
We now set the right-hand side of \eqref{l210} as $\mu^2$ in
the statement of the present proposition with an appropriate
choice of $c_3=c_3(n,k)$ and replacing $c_2$ by $e c_2$. 
Finally assume that we have $\|V_{t_0}\|(B_{4/5}
\cap\{|T^{\perp}(x)| > \mu\})>0$ for some
$t_0\in (1/5,\Lambda)$ when $V_{t_0}$ is unit density. Then there exists
some $x_0\in B_{4/5}$ with
$|T^{\perp}(x_0)|> \mu$ such that
$\Theta^k(\|V_{t_0}\|,x_0)=1$ and the approximate tangent 
space exists. One can then show that $\lim_{\varepsilon
\rightarrow 0+}\int_{B_1}|T^{\perp}(x)|^2\eta \rho_{(x_0,t_0+
\varepsilon)}\, d\|V_{t_0}\|=|T^{\perp}(x_0)|^2$. From 
\eqref{l210}, we should have $|T^{\perp}(x_0)|\leq \mu$,
which is a contradiction. We may easily extend the same conclusion 
for all $t\in (1/5,\Lambda)$ instead of a.e$.$ $t$ by \eqref{maineq}.
Thus we conclude the proof.
\hfill{$\Box$}
\label{moform}
\subsection{Time uniform $L^2$ estimate}
The following estimate is used in Section \ref{partialreg}. The proof is similar to
Proposition \ref{el2elinf}.
\begin{prop}
Corresponding to $1\leq E_1<\infty$, $p$ and $1\leq \Lambda<\infty$ there
exists $c_{17}$ with the following property. For $U=B_{LR}$, $2\leq L<\infty$,
suppose that $\{V_t\}_{0\leq t\leq \Lambda R^2}$ and $\{u(\cdot, t)\}_{0\leq t\leq \Lambda R^2}$
satisfy (A1)-(A4). Then for $T\in {\bf G}(n,k)$ and for all $t\in [0,\Lambda R^2]$,
\begin{equation}
\begin{split}
R^{-(k+2)}&\int_{B_R}|T^{\perp}(x)|^2\, d\|V_t\|
\leq \exp(1/(4\Lambda))R^{-(k+2)}\int_{B_{LR}}|T^{\perp}(x)|^2\, d\|V_0\| \\
&+c_{17}\{(R^{2\varsigma}
\|u\|^2+R^{\varsigma}\|u\|)L^2+L^{k+2}\exp(-(L-1)^2/(8\Lambda))\}.
\end{split}
\label{uniel2a}
\end{equation}
\label{uniel2}
Here $\|u\|=\|u\|_{L^{p,q}(B_{LR}\times[0,\Lambda R^2])}$. 
\end{prop}
{\it Proof}. Without loss of generality, we set $R=1$ after a change of variables. 
Let $\eta\in C^{\infty}_c(B_{L})$ be a radially symmetric non-negative function with $\eta\equiv 1$ on 
$B_{L-1}$, $0\leq \eta\leq 1$, and $|\nabla\eta|, \, \|\nabla^2\eta\|\leq c(n)$. 
In \eqref{maineq}, we use $\phi(x,t)=|T^{\perp}(x)|^2 \rho_{(0,2\Lambda)}(x,t)\eta(x)$ over
the time interval $[0,t_1]$, $t_1\leq \Lambda$. The computations are analogous to Proposition
\ref{el2elinf}, the only difference this time is that $\eta$ does not depend on time. Writing
$\rho=\rho_{(0,2\Lambda)}(x,t)$, we have (see \eqref{l25})
\begin{equation}
\begin{split}
&\left.\int |T^{\perp}(x)|^2 \rho\eta\, d\|V_t\|\right|_{t=0}^{t_1}\leq \int_{0}^{t_1}
\int\{-\nabla^2 (\eta|T^{\perp}(x)|^2)\cdot S+\eta|T^{\perp}(x)|^2|u|^2 \\
&+|T^{\perp}(x)|^2\nabla\eta\cdot u^{\perp}+\eta\nabla|T^{\perp}(x)|^2\cdot u^{\perp}\}
\rho\, dV_t(x,S)dt.
\end{split}\label{uniel2b}
\end{equation}
The last three terms involving $u$ can be estimated as in \eqref{Bcom5} and \eqref{Bcom6}, 
using also $|T^{\perp}(x)|\leq L$. We also use \eqref{l27} and \eqref{l275} (without 
$\partial\eta/\partial t$ term) and dropping negative term to obtain from \eqref{uniel2b}
\begin{equation}
\begin{split}
\left.\int |T^{\perp}(x)|^2 \rho\eta\, d\|V_t\|\right|_{t=0}^{t_1}&\leq c(n,k)\int_0^{t_1}\int
(\|\nabla^2 \eta\| L^2+L|\nabla\eta|)\rho\, d\|V_t\|dt \\
&+c(p,k,E_1,\Lambda)(\|u\|+\|u\|^2)L^2.
\end{split}
\label{uniel2c}
\end{equation}
Since $|\nabla\eta|$ and $\|\nabla^2\eta\|$ are zero on $B_{L-1}$, and since
$2\Lambda\geq 2\Lambda-t_1>\Lambda$, we have $\rho_{(0,2\Lambda)}(x,t)\leq (4\pi \Lambda)^{-k/2}
\exp(-(L-1)^2/8\Lambda)$ there. From \eqref{uniel2c}, one obtains
\begin{equation}
\left.\int |T^{\perp}(x)|^2\rho\eta \, d\|V_t\|\right|_{t=0}^{t_1}\leq c(p,n,k,E_1,\Lambda)\{(
\|u\|+\|u\|^2)L^2+L^{k+2}\exp(-(L-1)^2/8\Lambda)\}
\label{uniel2d}
\end{equation}
Since $\rho_{(0,2\Lambda)}(x,0)\leq (8\pi\Lambda)^{-k/2}$ for all $x\in B_L$, and $\rho_{(0,2\Lambda)}(x,t)\geq (8\pi\Lambda)^{-k/2}\exp(-1/(4\Lambda))$ for $|x|\leq 1$ and
$t\in [0,\Lambda]$, we obtain \eqref{uniel2a} from \eqref{uniel2d} with a suitable choice of 
constant. 
\hfill{$\Box$} 

\section{Parabolic Lipschitz approximation}
The main result of this section is Theorem \ref{lipapp} which gives a good 
Lipschitz graph approximation of moving varifold in space-time. 
The similar Lipschitz approximation has to be constructed in Allard's regularity theory
and one can see a clear parallel in that sense. We note that
along with the monotonicity type estimates in the previous section, 
Lipschitz approximation is completely absent
in Brakke's original proof.  Instead Brakke constructed approximate
graphs by substituting test functions weighted by heat kernel in the
flow inequality (see \cite[6.8, 6.9]{Brakke}).

\begin{prop} Given $0<l<\infty$, there exist constants $c_4=c_4(l)$, 
$c_5=c_5(n,k)$ and $c_6=c_6(k,p,q)$ with the following property.
For $T\in {\bf G}(n,k)$, $U=C(T,1)$, $\Lambda=1$ and $1\leq E_1<\infty$
suppose that $\{V_t\}_{0\leq t\leq 1}$
and $\{u(\cdot,t)\}_{0\leq t\leq 1}$ satisfy (A1)-(A4). 
For $\max\{2,l\}<R<\infty$, $0<\gamma<1$ and $(y_1, s_1),\,
(y_2,s_2)\in C(T,1)\times (0,1)$, 
we assume the following.
\begin{equation}
\Theta^k(\|V_{s_i}\|,y_i)=1\,\,\,\mbox{ for $i=1,2$.}
\label{denhyp}
\end{equation}
\begin{equation}
\sigma:=|T^{\perp}(y_1-y_2)|\geq l^{-1}\max\{|T(y_1-y_2)|,\, 
|s_1-s_2|^{\frac12}\}.
\label{notlip}
\end{equation}
Set ${\bar y}:=\frac{y_1+y_2}{2}$, ${\bar s}:=\frac{s_1+s_2}{2}$ and
$R(l,\sigma):=(R^{3/2}+\frac{\sqrt{l^2+1}}{2})\sigma$. Then assume further that
\begin{equation}
P_{R(l,\sigma)}({\bar y},{\bar s})\subset
C(T,1)\times (0,1),
\label{incl}
\end{equation}
and for each $i=1,2$ we have
\begin{equation}
\int_{P_{R^{3/2}\sigma}(y_i,s_i)}\|S-T\|^2 \, dV_t(\cdot, S)dt\leq \gamma
(R^{3/2}\sigma)^{k+2},
\label{tiltde}
\end{equation}
and $R$ and $l$ satisfy
\begin{equation}
R^2+\frac{l^2}{2}<\frac{(R^{3/2}-1)^2}{2k}.
\label{mass}
\end{equation}
Under the above assumptions, we have
\begin{equation}
\begin{split}
&2
\leq \frac{\exp(c_4/\sqrt{R})}{\left(1-\frac{l^2}{2R^2}\right)^{k/2}}
\left.\int_{B_{R(l,\sigma)}({\bar y})}
\rho_{({\bar y},{\bar s}) }(\cdot,t)\, d\|V_t\|
\right|_{t={\bar s}-R^2\sigma^2}+c_5 E_1 (\exp(-\sqrt{R}) +\sigma^{\varsigma} R^2) \\
&+c_5 \gamma R^{3(k+2)/2}+c_6
\|u\|^2_{L^{p,q}(C(T,1)\times(0,1))}E_1^{1-\frac{2}{p}}
R^{2\varsigma}(\sigma^{2\varsigma}+\sigma^{\varsigma}).
\end{split}
\label{finest}
\end{equation}
\label{lip1}
\end{prop}
\begin{rem}
Proposition \ref{lip1} is analogous to \cite[6.1]{Allard}. To describe the
idea, assume $n=k+1$ case. Consider $T$ as a horizontal plane. 
Then $\sigma$ is the vertical distance of two points, which is typically very small. 
Later $R$ will be chosen very large so that the coefficient of the first term of the
right-hand side of \eqref{finest} is very close to 1. In Lemma \ref{liphelp}
we prove that there can be only one horizontal sheet in terms of density ratio 
via monotonicity formula if the first variation is small. Thus we conclude
that there cannot be two very 
close points of density one positioned in a vertical manner. 
Rough idea of the proof is that we `cut' moving varifold by a horizontal
plane to separate 
these two points, and derive
a monotonicity type formula for each piece. 
The cutting naturally produces extra error terms
which are assumed to be small (see \eqref{tiltde}). It is worth pointing out that the
error bound assumption is only made for the tilt-excess
and not for the first variation. In Allard's analogous Lipschitz approximation,
one makes assumptions on both. 
\end{rem}
{\it Proof}.
Let $\xi_1,\,\xi_2\in C^{\infty}({\mathbb R})$ be non-negative functions
such that 
\begin{equation}
\begin{split}
& \xi_1(s)=\left\{\begin{array}{ll}1 & \mbox{for }s\leq \sigma/4,\\
0 & \mbox{for }s\geq \sigma/2,
\end{array}\right.\hspace{1cm}
\xi_2(s)=\left\{\begin{array}{ll}0 & \mbox{for }s\leq \sigma/2,\\
1 & \mbox{for }s\geq 3\sigma/4, 
\end{array}\right. \\
&0\leq \xi_i\leq 1,\,|\xi_i'|\leq 8/\sigma\mbox{ and }
|\xi_i''|\leq 128/\sigma^2\mbox{ for $i=1,2$}.
\end{split}
\label{xicond}
\end{equation}
Let $\eta_1,\, \eta_2\in C_c^{\infty}({\mathbb R}^n)$ be non-negative functions
such that for each $i=1,2$,
\begin{equation}
\begin{split}
& \eta_i(x)=\left\{\begin{array}{ll} 1 &\mbox{on }B_{(R^{3/2}-1)\sigma}
(y_i), \\ 0 &\mbox{on }{\mathbb R}^n\setminus B_{R^{3/2}\sigma}(y_i),
\end{array}\right. \\
& 0\leq \eta_i\leq 1,\,|\nabla\eta_i|\leq 2/\sigma\mbox{ and } 
\|\nabla^2\eta_i\|\leq c(n)/\sigma^2.
\end{split}
\label{chicond}
\end{equation}
Due to \eqref{notlip}, for each $i=1,2$, $|{\bar y}-y_i|=|y_1-y_2|/2
\leq \sigma\sqrt{l^2+1}/2$, and \eqref{incl} shows
\begin{equation}
\{x\in {\mathbb R}^n\, :\, \eta_i(x)\neq 0\}\subset
B_{R^{3/2}\sigma}(y_i)\subset B_{R^{3/2}\sigma+|{\bar y}-y_i|}
({\bar y})\subset B_{R(l,\sigma)}({\bar y})\subset C(T,1).
\label{ball}
\end{equation}
For each $i=1,2$ and arbitrarily small $\varepsilon>0$ 
we compute \eqref{maineq} with 
\begin{equation*}
\phi(x,t)=\rho_{(y_i,s_i+\varepsilon)}
(x,t)\eta_i(x)\xi_i(|T^{\perp}(x-y_1)|)
\end{equation*}
for the time interval 
$t_1={\bar s}-R^2\sigma^2$ and $t_2=s_i$. 
By \eqref{notlip} and $l<R$, we have $t_2-t_1\geq R^2\sigma^2-|s_1-s_2|/2
\geq R^2\sigma^2-l^2\sigma^2/2>0$. In the following, we denote
$\rho_{(y_i,s_i+\varepsilon)}(x,t)$ (and even after setting $\varepsilon=0$) 
by $\rho_i$.
The same computations 
leading towards \eqref{l25} give (we omit writing variables with no fear of
confusions)
\begin{equation}
\left.\int \rho_i\eta_i\xi_i\, d\|V_t\|\right|_{t={\bar s}-R^2\sigma^2}^{
s_i}\leq \int_{{\bar s}-R^2\sigma^2}^{s_i}dt\int
\{-\nabla^2(\eta_i\xi_i)\cdot S+\eta_i\xi_i|u|^2+\nabla(\eta_i
\xi_i)\cdot u^{\perp}\}\rho_i\, dV_t(\cdot,S).
\label{lip2}
\end{equation}
Due to \eqref{xicond} and \eqref{chicond}, $\eta_i(x)\xi_i(|T^{\perp}(x-y_1)|)$
is identically 1 in a small neighborhood of $x=y_i$ for each $i=1,2$. 
Hence using \eqref{denhyp} we have 
\begin{equation*}
\lim_{\varepsilon\rightarrow 0+}
\int\rho_{(y_i,s_i+\varepsilon)}\eta_i\xi_i\, d\|V_{s_i}\|
=\Theta^k(\|V_{s_i}\|,y_i)=1.
\end{equation*}
Thus by letting $\varepsilon\rightarrow 0$, \eqref{lip2} is
\begin{equation}
\begin{split}
1\leq &\left.\int\rho_i\eta_i\xi_i\, d\|V_t\|
\right|_{t={\bar s}-R^2\sigma^2} 
+\int_{{\bar s}-R^2\sigma^2}^{s_i}dt\int\{-\nabla^2
(\eta_i\xi_i)\cdot S \\ &+\eta_i\xi_i|u|^2+\nabla(\eta_i\xi_i)
\cdot u^{\perp}\}\rho_i\, dV_t(\cdot,S) 
=:I_1+I_2+I_3+I_4
\end{split}
\label{lip3}
\end{equation}
where $\rho_i=\rho_{(y_i,s_i)}(x,t)$. We next estimate the right-hand
side of \eqref{lip3}. We start with $I_2$. 

{\bf Estimate for $I_2$}.
\newline
For the integrand of $I_2$ 
we have (writing $v=T^{\perp}(x-y_1)/|T^{\perp}(x-y_1)|$)
\begin{equation}
\begin{split}
&|-\nabla^2(\eta_i\xi_i)\cdot S|\leq \xi_i|\nabla^2\eta_i\cdot S|
+2|(\nabla\eta_i\otimes\nabla\xi_i)\cdot S|+\eta_i|\nabla^2\xi_i\cdot S|\\
&\leq \xi_i k\|\nabla^2\eta_i\|+2|\xi_i'\nabla\eta_i\cdot S(v)|+\eta_i
\left|S\cdot\left\{\left(\xi_i''-\frac{\xi_i'}{|T^{\perp}(x-y_1)|}\right)v\otimes v+
\frac{\xi_i'}{|T^{\perp}(x-y_1)|}T^{\perp}\right\}\right|.
\end{split}
\label{mcom1}
\end{equation}
Note that (see Lemma \ref{re-ea}) we have
\begin{equation}
\begin{split}
&|S(T^{\perp}(v))|\leq \|S-T\||v|, \\
& S\cdot(T^{\perp}(v)\otimes
T^{\perp}(v))=S(T^{\perp}(v))\cdot S(T^{\perp}(v))\leq \|S-T\|^2 |v|^2,\\
&S\cdot T^{\perp}\leq k\|S-T\|^2.
\end{split}
\label{linear}
\end{equation}
Thus with \eqref{xicond}, \eqref{chicond}, \eqref{linear}, 
and the fact that $|T^{\perp}(x-y_1)|\geq \sigma/4$ on $\{\nabla\xi_i\neq 0\}$,
\eqref{mcom1} is estimated as 
\begin{equation}
|-\nabla^2(\eta_i\xi_i)\cdot S|\leq c(n,k)\sigma^{-2}\chi_{{\rm spt}\,|\nabla
\eta_i|\cap\,{\rm spt}\,\xi_i}+c(k)\sigma^{-2}\|S-T\|^2\chi_{{\rm spt}\,\eta_i
\cap\,{\rm spt}\,|\nabla\xi_i|}
\label{mcom2}
\end{equation}
with a suitable choice of constants $c(n,k)$ and $c(k)$. We next estimate
$\rho_i$ on ${\rm spt}|\nabla
\eta_i|\cap\,{\rm spt}\,\xi_i$ and ${\rm spt}\,\eta_i
\cap\,{\rm spt}|\nabla\xi_i|$ respectively. 

On ${\rm spt}|\nabla
\eta_i|\cap\,{\rm spt}\,\xi_i$, since ${\rm spt}|\nabla\eta_i|
\cap B_{(R^{3/2}-1)\sigma}(y_i)=\emptyset$ by \eqref{chicond}, 
\begin{equation}
\rho_i\leq \frac{\exp\left(-\frac{(R^{3/2}-1)^2\sigma^2}{4(s_i-t)}\right)}{
(4\pi(s_i-t))^{k/2}}
\label{esrho1}
\end{equation}
which takes the maximum value at $s_i-t=(R^{3/2}-1)^2\sigma^2/(2k)$ and 
monotone increasing until that point. Note that $s_i-t$ varies over the interval
$(0,s_i-{\bar s}+R^2\sigma^2)\subset(0,(R^2+\frac{l^2}{2})\sigma^2)$. 
Thus as long as \eqref{mass} is satisfied, the maximum value of \eqref{esrho1}
is estimated as 
\begin{equation}
\rho_i\leq \frac{\exp\left(-\frac{(R^{3/2}-1)^2}{4(R^2+l^2/2)}\right)}{
(4\pi)^{k/2}(R^2+l^2/2)^{k/2}\sigma^k}\leq \frac{\exp\left(-\frac{(R^{3/2}-1)^2}{6R^2}\right)}{
(4\pi)^{k/2}R^k\sigma^k}\leq \frac{c(k)\exp(-R^{3/4})}{\sigma^k}
\label{esrho2}
\end{equation}
on $\{{\rm spt}|\nabla
\eta_i|\cap{\rm spt}\,\xi_i\}$. The last two inequalities follow from 
$R>\max\{2,l\}$ and due to the exponential decay for large $R$. 

For $x\in \{{\rm spt}\,\eta_i
\cap{\rm spt}|\nabla\xi_i|\}$, we have $|x-y_i|\geq \sigma/4$, and 
one can check that the maximum value can be estimated as 
\begin{equation}
\rho_i\leq \frac{\exp\left(-\frac{\sigma^2}{64(s_i-t)}\right)}{
(4\pi(s_i-t))^{k/2}}\leq \frac{(32k)^{k/2}\exp(-k/2)}{(4\pi)^{k/2}
\sigma^k}.
\label{esrho3}
\end{equation}
Hence we may estimate $I_2$ by combining \eqref{mcom2}, \eqref{esrho2} and
\eqref{esrho3} as (with yet another suitable choice of constants)
\begin{equation}
\begin{split}
I_2&\leq \frac{c(n,k)\exp(-R^{3/4})}{\sigma^{k+2}}\int_{{\bar s}-R^2\sigma^2}^{s_i}
\|V_t\|(B_{R^{3/2}\sigma}(y_i))dt \\ &+\frac{c(k)}{\sigma^{k+2}}
\int_{{\bar s}-R^2\sigma^2}^{s_i}\int_{B_{R^{3/2}\sigma}(y_i)}\|S-T\|^2
\, dV(\cdot,S)dt.
\end{split}
\label{mcom3}
\end{equation}
By \eqref{dbd}, $|s_i-{\bar s}+R^2\sigma^2|\leq (R^2+l^2/2)\sigma^2$, 
$B_{R^{3/2}\sigma}(y_i)\times(s_i,{\bar s}-R^2\sigma^2)\subset P_{R^{3/2}\sigma}
(y_i,s_i)$ and \eqref{tiltde}, with a suitable choice of constants we have from
\eqref{mcom3}
\begin{equation}
I_2\leq c(n,k)E_1 \exp(-\sqrt{R})+c(k)\gamma R^{3(k+2)/2}.
\label{mcom4}
\end{equation}

{\bf Estimate for $I_3$}.
\newline
$I_3$ can be estimated as in \eqref{Bcom5} and \eqref{Bcom6}, so we have
\begin{equation}
I_3\leq c_2(k,p,q)\|u\|_{L^{p,q}(C(T,1)\times(0,1))}^2 E_1^{1-\frac{2}{p}}
(3R^2\sigma^2/2)^{\varsigma}.
\label{I3}
\end{equation}

{\bf Estimate for $I_4$}.
\newline
We have 
\begin{equation}
|\rho_i\nabla(\eta_i\xi_i)\cdot u^{\perp}|=|\rho_i S^{\perp}(\nabla(\eta_i\xi_i))\cdot u|
\leq \frac{\sigma^{-\varsigma}}{2} |u|^2\rho_i\chi_{{\rm spt}\,\eta_i}+\frac{\sigma^{\varsigma}}{2} |\nabla(\eta_i
\xi_i)|^2\rho_i.
\label{mcom5}
\end{equation}
The first term of the right-hand side of \eqref{mcom5} can be handled just like
$I_3$ (with multiplication of $\sigma^{-\varsigma}$). The second term is
\begin{equation}
\frac{\sigma^{\varsigma}}{2}|\nabla(\eta_i\xi_i)|^2\rho_i\leq 
\frac{\sigma^{\varsigma}}{2}(|\nabla\eta_i|+|\nabla\xi_i|)^2\rho_i\chi_{{\rm spt}\, \eta_i} \\
\leq \frac{\sigma^{\varsigma}}{2}\left(\frac{10}{\sigma}\right)^2\rho_i\chi_{{\rm spt}\, \eta_i}
\label{mcom55}
\end{equation}
by \eqref{chicond} and \eqref{xicond}.
Thus \eqref{mcom5} and \eqref{mcom55} combined with \eqref{dbd} and 
$|s_i-\bar{s}+R^2\sigma^2|\leq (R^2+l^2/2)\sigma^2$ show
\begin{equation}
I_4\leq c_2 \|u\|^2_{L^{p,q}(C(T,1)\times(0,1))}E_1^{1-\frac{2}{p}}(3R^2/2)^{\varsigma} \sigma^{\varsigma}
+75 \sigma^{\varsigma} R^2 E_1.
\label{mcom55.1}
\end{equation}

{\bf Estimate for $I_1$}.
\newline
It is important to note that 
\begin{equation}
\sum_{i=1}^2\left.\int\rho_i\eta_i\xi_i\, d\|V_t\|\right|_{t={\bar s}-R^2\sigma^2}
\leq \left.\int_{B_{R(l,\sigma)}({\bar y})}\xi_1\rho_1+\xi_2\rho_2\, d\|V_t\|
\right|_{t={\bar s}-R^2\sigma^2}
\label{mcom6}
\end{equation}
by \eqref{ball} and that the supports of $\xi_1$ and $\xi_2$ are disjoint. 
We next compute the ratio of $\rho_i$ and $\rho_{({\bar y},{\bar s})}$.
With $|s_i-{\bar s}|\leq l^2 \sigma^2/2$
\begin{equation}
\frac{\rho_{(y_i,s_i)}(x,{\bar s}-R^2\sigma^2)}{\rho_{({\bar y},{\bar s})}(x,
{\bar s}-R^2\sigma^2)}=\frac{R^k\sigma^k\exp\left(\frac{|{\bar y}-x|^2}{4R^2\sigma^2}-\frac{|y_i-x|^2}{4(s_i-{\bar s}
+R^2\sigma^2)}\right)}{(s_i-{\bar s}+R^2
\sigma^2)^{k/2}} 
\leq \frac{\exp\left(\frac{|{\bar y}-x|^2}{4R^2\sigma^2}-\frac{|y_i-x|^2}{
4(R^2+l^2/2)\sigma^2}\right)}{\left(1-\frac{l^2}{2R^2}\right)^{k/2}}.
\label{mcom7}
\end{equation}
One checks that when $|{\bar y}-x|\leq |{\bar y}-y_i|(\leq \sigma\sqrt{l^2+1}/2)$, the argument of
exponential function in \eqref{mcom7} is
\begin{equation}
\frac{\left(R^2+\frac{l^2}{2}\right)|{\bar y}-x|^2-R^2|y_i-x|^2}{
4R^2\left(R^2+\frac{l^2}{2}\right)\sigma^2}\leq \frac{|{\bar y}-x|^2}{
4R^2\sigma^2}\leq \frac{l^2+1}{16R^2}.
\label{mcom8}
\end{equation}
When $|{\bar y}-x|>|{\bar y}-y_i|$, we have
$|y_i-x|\geq |{\bar y}-x|-|{\bar y}-y_i|>0$ and $-|y_i-x|^2\leq -|{\bar y}
-x|^2+2|{\bar y}-x||{\bar y}-y_i|$. Thus
\begin{equation}
\begin{split}
&\frac{\left(R^2+\frac{l^2}{2}\right)|{\bar y}-x|^2-R^2|y_i-x|^2}{
4R^2\left(R^2+\frac{l^2}{2}\right)\sigma^2}\leq \frac{\frac{l^2}{2}|{\bar y}
-x|^2+2R^2|{\bar y}-x||{\bar y}-y_i|}{4R^4\sigma^2} \\
&\leq \frac{\frac{l^2}{2}(R^{3/2}+\sqrt{l^2+1}/2)^2+2R^2(R^{3/2}+\sqrt{l^2+1}/2)
(\sqrt{l^2+1}/2)}{4R^4}\leq c(l)/\sqrt{R}
\end{split}
\label{mcom9}
\end{equation}
since $|{\bar y}-x|\leq (R^{3/2}+\sqrt{l^2+1}/2)\sigma$ and by comparing the
exponents of $R$. Thus, combining \eqref{mcom7}, \eqref{mcom8} and \eqref{mcom9},
we obtain
\begin{equation}
\rho_{(y_i,s_i)}(x,{\bar s}-R^2\sigma^2)
\leq \frac{\exp(c(l)/\sqrt{R})}{\left(1-\frac{l^2}{2R^2}
\right)^{k/2}}\rho_{({\bar y},{\bar s})}(x,
{\bar s}-R^2\sigma^2).
\label{mcom10}
\end{equation}
Combining \eqref{mcom6} and \eqref{mcom10}, and the fact that $\xi_1$ and
$\xi_2$ have disjoint supports, we obtain
\begin{equation}
\sum_{i=1}^2\left.\int\rho_i\eta_i\xi_i\, d\|V_t\|\right|_{t={\bar s}-
R^2\sigma^2}\leq \frac{\exp(c(l)/\sqrt{R})}{\left(1-\frac{l^2}{2R^2}\right)^{k/2}}
\left.\int_{B_{R(l,\sigma)}({\bar y})}\rho_{({\bar y},{\bar s})}(\cdot,t)\, d\|V_t\|
\right|_{t={\bar s}-R^2\sigma^2}.
\label{mcom11}
\end{equation}

{\bf Final step}.
\newline
Finally summing \eqref{lip3} for $i=1,2$, and by the estimates \eqref{mcom4},
\eqref{I3}, \eqref{mcom55.1} and \eqref{mcom11}, we obtain \eqref{finest} with
a suitable choice of constants.
\hfill{$\Box$}

In the next Lemma \ref{liphelp} and Theorem \ref{lipapp} we use 
\begin{equation*}
\hat{\rho}_{(y,s)}(x,t)=\eta(7x/5)\rho_{(y,x)}(x,t),
\label{difrho}
\end{equation*}
where $\eta$ is defined in Section \ref{monosubsec}. In the statement of Proposition \ref{cam1},
this corresponds to choosing $R=5/7$. For $\phi_T$ as in \eqref{Tphi}, we have
\begin{equation}
{\rm spt}\, \hat{\rho}_{(y,s)}(\cdot,t)\subset B_{2/3}\subset\{\phi_T=1\},\,\,
\hat{\rho}_{(y,s)}(\cdot,t)=\rho_{(y,s)}(\cdot, t)\,\, \mbox{on $B_{13/21}$}.
\label{phisup}
\end{equation}
\begin{lemma}
Corresponding to $1\leq E_1<\infty$, $p$ and $q$ with \eqref{expcond}
there exist $0<r_1<1$ and $0<\varepsilon_2<1$ with
the following property. For $U=C(T,1)$, $\Lambda=1$ and $E_1$, 
suppose $\{V_t\}_{0\leq t\leq 1}$ and $\{u(\cdot,t)\}_{
0\leq t\leq 1}$ satisfy (A1)-(A4). Let $\phi_T$ and
${\bf c}$ be as in
\eqref{Tphi}.
\newline
Assume
\begin{equation}
\int_{C(T,1)\times(0,1)}|h(V_t,\cdot)|^2\phi_T^2\, d\|V_t\|dt\leq \varepsilon_2 r_1^2/4,
\label{hhyp}
\end{equation}
\begin{equation}
\left|\|V_t\|(\phi_T^2)-{\bf c}\right|\leq \varepsilon_2,\hspace{.5cm}0<\forall t<1,
\label{chika1}
\end{equation}
\begin{equation}
{\rm spt}\,\|V_t\|\cap C(T,1)\subset \{|T^{\perp}(x)|\leq  \varepsilon_2\},\hspace{.5cm}
0<\forall t<1,
\label{chika2}
\end{equation}
\begin{equation}
\|u\|_{L^{p,q}(C(T,1)\times(0,1))}\leq 1.
\label{chiu}
\end{equation}
Then for any $y\in C(T,1/2)\cap\{|T^{\perp}(x)|\leq \varepsilon_2\}$, $t\in (r_1^2, 1)$ and $0<r<r_1/\sqrt{2}$, 
we have
\begin{equation}
\int_{C(T,1)}{\hat{\rho}}_{(y,t+r^2)}(\cdot,t)\, d\|V_t\|\leq \frac{26}{25}.
\label{chika3}
\end{equation}
\label{liphelp}
\end{lemma}
\begin{rem}
In this Lemma we take advantage that quantities \eqref{hhyp}-\eqref{chiu} can
be made small when the $L^2$-height is small due to $L^2$-$L^{\infty}$
estimate and energy estimate. We then conclude that the weighted density
ratio can be estimated from above by a number close to 1 \eqref{chika3}
for all small radii.  
\end{rem}
{\it Proof}. We fix $0<r_1<1$ so that 
\begin{equation}
c_2 E_1^{1-\frac{2}{p}}(r_1^2)^{\varsigma}+c_1 E_1 \o_k r_1^2 (7/5)^2<\frac{1}{100},
\label{ap1b}
\end{equation}
\begin{equation}
\sup_{0<r\leq r_1}\,\sup_{y\in B_{1/2}\cap T}\left|\int_T
{\hat{\rho}}_{(y,r^2)}(x,0)\, d{\mathcal H}^k(x)-1\right|<\frac{1}{100}.
\label{ap1c}
\end{equation}
In regard to \eqref{ap1c}, we may choose such $r_1$ since $\hat{\rho}_{(y,r^2)}(\cdot,0)
=\rho_{(y,r^2)}(\cdot,0)$ 
on $B_{13/21}$ due to \eqref{phisup} and $\hat{\rho}_{(y,r^2)}(\cdot,0)$ approaches 
to the delta mass on $T$ at $y$ as $r\rightarrow 0$. The choice of $r_1$ depends
on $E_1$, $p$ and $q$. By \eqref{monoest}, for all $y\in B_{4/7}$,
$0<t_1<t<1$ and $t-t_1<r_1^2$ we have with \eqref{chiu} and \eqref{ap1b}
\begin{equation}
\int_{B_1}{\hat{\rho}}_{(y,t+r^2)}(\cdot,t)\, d\|V_t\|-\int_{B_1}{\hat{\rho}}_{(y,t+r^2)}(\cdot,
t_1)\, d\|V_{t_1}\|
<\frac{1}{100}.
\label{ap2}
\end{equation}
We next claim that there exists a positive constant $0<\varepsilon_2<1$
depending on $r_1$ and $E_1$ with the following properties. Assume that
$V$ is an integral varifold with
\begin{equation}
\|V\|(B_1)\leq E_1,
\label{ap3a}
\end{equation}
\begin{equation}
\int_{C(T,1)}|h(V,\cdot)|^2\phi_T^2 \, d\|V\|\leq \varepsilon_2,
\label{ap3b}
\end{equation}
\begin{equation}
{\rm spt}\,\|V_t\|\cap C(T,1)\subset \{|T^{\perp}(x)|\leq  \varepsilon_2\},
\label{ap3c}
\end{equation}
\begin{equation}
\left|\|V\|(\phi_T^2)-{\bf c}\right|
\leq \varepsilon_2.
\label{ap3d}
\end{equation}
Then for all $y\in C(T,1/2)\cap\{|T^{\perp}(x)|\leq  \varepsilon_2\}$
and $r_1/\sqrt{2}\leq r\leq r_1$, we have
\begin{equation}
\left|\int_{B_1}{\hat{\rho}}_{(y,r^2)}(\cdot,0)\, d\|V\|-1\right|
<\frac{1}{50}.
\label{ap3e}
\end{equation}
If false, then for each $m\in {\mathbb N}$
there should exist an integral varifold $V_m$, 
$y_m\in C(T,1/2)\cap\{|T^{\perp}(x)|\leq 1/m\}$ and $r_1 /\sqrt{2}\leq r_m\leq r_1$ satisfying
\eqref{ap3a}-\eqref{ap3d} for $\varepsilon_2=1/m$ but not satisfying
\eqref{ap3e}. Then by \eqref{ap3a} we would have a converging subsequence
with limit $V$. Due to \eqref{ap3b} and the compactness 
theorem of integral varifold, $V$ is a stationary integral varifold on 
$\{\phi_T>0\}$. Due to \eqref{ap3c}, ${\rm spt}\, \|V\|\subset T$. 
By the constancy theorem of integral varifold and due to \eqref{ap3d},
we may conclude that $\|V\|={\mathcal H}^k\lfloor_T$ on $\{\phi_T>0\}$. 
We may also assume that $y_m\rightarrow
y\in \overline{C(T,1/2)}\cap T$ and $r_m\rightarrow r\in [r_1/\sqrt{2},r_1]$
and we have 
\begin{equation*}
\left|\int_{B_1}{\hat{\rho}}_{(y,r^2)}(x,0)\, d\|V\|(x)-1\right|\geq \frac{1}{50}.
\end{equation*}
This contradicts with \eqref{ap1c} since ${\rm spt}\,\hat{\rho}_{(y,r^2)}(\cdot,0)
\subset\{\phi_T>0\}$  by \eqref{phisup} and $\|V\|={\mathcal H}^k\lfloor_T$. 

Now given $y\in C(T,1/2)\cap\{|T^{\perp}(x)|\leq \varepsilon_2\}$, 
$t\in (r_1^2,1)$ and $0<r<r_1/\sqrt{2}$, consider
the interval $[t+r^2-r_1^2, t+r^2-r_1^2/2]$ which has the length $r_1^2/2$. 
Then by \eqref{hhyp}, there exists some $t_1$ in this interval such that
\eqref{ap3b} holds for $V_{t_1}$. The conditions \eqref{ap3a}, \eqref{ap3c} 
and \eqref{ap3d} are satisfied for $V_{t_1}$ by the assumptions of the present
lemma. Since $r_1^2/2\leq t+r^2-t_1\leq r_1^2$, we have from \eqref{ap3e}
\begin{equation}
\left|\int_{B_1}{\hat{\rho}}_{(y,t+r^2)}(\cdot,t_1)\, d\|V_{t_1}\|-1\right|<\frac{
1}{50}.
\label{ap5}
\end{equation}
Since $0\leq t-t_1\leq r_1^2$, combining \eqref{ap5} and \eqref{ap2}, we obtain 
the desired result.
\hfill{$\Box$}

\begin{thm}
Corresponding to $1\leq E_1<\infty$, $p$ and $q$ with \eqref{expcond} there
exists $0<\varepsilon_3<1$ and $c_7$ with the following property. For $U=C(T,1)$,
$\Lambda=1$ and $1\leq E_1<\infty$, suppose $\{V_t\}_{0\leq t\leq 1}$ and $\{u(\cdot,t)\}_{0\leq t\leq 1}$
satisfy (A1)-(A4). Let $\phi_T$, $\phi_{T,1/2}$ and ${\bf c}$ be as in \eqref{Tphi} and let $r_1$ be as in Lemma \ref{liphelp}. Write $|M_t|:=V_t$ for a.e$.$ $t$
and identify $T$ with ${\mathbb R}^k\times\{0\}$. 
Suppose that we have
\begin{equation}
\int_{C(T,1)\times(0,1)}|h(V_t,\cdot)|^2\phi_T^2 \, d\|V_t\|dt\leq \varepsilon_3 r_1^2 /4,
\label{hhypsub}
\end{equation}
\begin{equation}
\left|\|V_t||(\phi_T^2)-{\bf c}\right|\leq \varepsilon_3,\,\,0<\forall t<1,
\label{chika1sub}
\end{equation}
\begin{equation}
{\rm spt}\, \|V_t\|\cap C(T,1)\subset\{|T^{\perp}(x)|\leq \varepsilon_3\},\,\, 0<\forall t<1, 
\label{chika2sub}
\end{equation}
\begin{equation}
\|u\|_{L^{p,q}(C(T,1)\times (0,1))}\leq 1.
\label{chika2sss}
\end{equation}
Set 
\begin{equation}
\beta^2:=\int_{G_k(C(T,1))\times(0,1)}\|S-T\|^2\phi_T^2 \, dV_t(\cdot,S)dt
\label{excess}
\end{equation}
and 
\begin{equation}
\kappa^2:=\left| \int_0^1 \|V_t\|(\phi_{T,1/2}^2)\, dt - \frac{\bf c}{2^k}\right|.
\label{excesslip}
\end{equation}
Then there exist maps $f:B^k_{1/3}
\times(1/4,3/4)\rightarrow {\mathbb R}^{n-k}$
and $F: B_{1/3}^k
\times(1/4, 3/4)\rightarrow {\mathbb R}^n\times (0,1)$
such that for all $(x,s),\, (y,t)\in B_{1/3}^k\times(1/4,
3/4)$, 
\begin{equation}
\begin{split}
& F(x,s)=(x,f(x,s),s),\\
&|f(x,s)-f(y,t)|\leq c(n,k)\max\{|x-y|,|s-t|^{1/2}\},\hspace{.5cm}
|f(x,s)|\leq \varepsilon_3.
\end{split}
\label{lipok}
\end{equation}
Define $X=\{\cup_{t\in (1/4,3/4)}(M_t\cap C(T,1/3))\times\{t\}
\}\cap {\rm image}\, F$ and $Y=(T\times {\rm Id}_{\mathbb R})(X)$. 
Then
\begin{equation}
\begin{split}
(\|V_t\|\times dt)&(\{C(T,1/3)\times(1/4,3/4)\}\setminus X)\\
&+{\mathcal L}^{k+1}(\{ B^k_{1/3}\times(1/4,3/4)\}
\setminus Y)\leq \kappa^2+c_7\beta^2.
\end{split}
\label{lipok2}
\end{equation}
\label{lipapp}
\end{thm}
{\it Proof}.
We use Lemma \ref{lip1} with $l=1$ and $R,\,\gamma,\,\sigma$ chosen as follows. 
We choose and fix $R$ sufficiently large depending only on $n,\, k,\, E_1$ so
that \eqref{mass} is satisfied and 
\begin{equation}
1<\frac{\exp(c_4/\sqrt{R})}{\left(1-\frac{1}{2R^2}\right)^{k/2}}<1+\frac{1}{100},
\hspace{1.cm}c_5 E_1 \exp(-\sqrt{R})<\frac{1}{100}.
\label{ap4a}
\end{equation}
After fixing such $R$, we choose and fix $\gamma>0$ sufficiently small 
depending only on $n,\, k,\, E_1$ so that 
\begin{equation}
c_5 \gamma R^{3(k+2)/2}<\frac{1}{100}.
\label{ap4b}
\end{equation}
We choose and fix $\sigma$ sufficiently small depending only on $n,\, k,\,
E_1,\, p,\, q$ so that
\begin{equation}
c_6  E_1^{1-\frac{2}{p}}R^{2\varsigma}( \sigma^{2\varsigma}+\sigma^{\varsigma})
 <\frac{1}{100},\hspace{1.cm} c_5 E_1 R^2 \sigma^{\varsigma}<\frac{1}{100},
\label{ap4c}
\end{equation}
\begin{equation}
(R^{3/2}+\sqrt{2}/2)\sigma\leq r_1,\hspace{1.cm}
R^2\sigma^2\leq r_1^2/2.
\label{ap4d}
\end{equation}
We set 
\begin{equation}
\varepsilon_3=\min\{\varepsilon_2,\sigma/2\}.
\label{ap5a}
\end{equation}
We then define
\begin{equation}
\begin{split}
A=\{(x,t)\, :\, x & \in C(T,1/2)\cap M_t,\, |T^{\perp}(x)|\leq \varepsilon_3,
\,\, t\in (1/4,3/4), \\
& \sup_{0<r<r_1}\frac{1}{r^{k+2}}\int_{P_r(x,t)}\|S-T\|^2\, dV_s(\cdot,S)ds\leq 
\gamma\}.
\end{split}
\label{apdef}
\end{equation}
We apply Lemma \ref{lip1} and show that $A$ can be contained in a 
Lipschitz graph. To do so, assume that we have distinct 
$(y_1,s_1)$ and $(y_2,s_2)$ in $A$
with $|T^{\perp}(y_1-y_2)|\geq \max\{|T(y_1-y_2)|,|s_1-s_2|^{1/2}\}$.
Set $\bar{\sigma}=|T^{\perp}(y_1-y_2)|$. By \eqref{ap5a} and \eqref{apdef},
$\bar{\sigma}\leq |T^{\perp}(y_1)|+|T^{\perp}(y_2)|\leq \sigma$. 
Since $R^{3/2}\bar{\sigma}\leq R^{3/2}\sigma<r_1$ by \eqref{ap4d}, we have 
\eqref{tiltde} satisfied by \eqref{apdef} for $i=1,\, 2$. 
We also have
\begin{equation}
P_{(R^{3/2}+\frac{\sqrt{2}}{2})\bar{\sigma}}({\bar y},{\bar s})\subset
C(T,13/21)\times(0,1)
\label{ap8}
\end{equation}
due to \eqref{ap4d} and restricting $r_1$ if necessary (depending on absolute constant). 
Thus we have \eqref{denhyp}-\eqref{mass}
satisfied for $(y_1,s_1)$ and $(y_2,s_2)$. Then we should have \eqref{finest}.
By \eqref{ap4a}-\eqref{ap4c} and \eqref{chika2sss} we have
\begin{equation}
2\leq (1+\frac{1}{100})\left.\int_{B_{(R^{3/2}+\sqrt{2}/2)\bar{\sigma}}({\bar y})}
\rho_{({\bar y},{\bar s})}(\cdot,t)\, d\|V_t\|\right|_{t={\bar s}-R^2{\bar{\sigma}}^2}
+\frac{4}{100}.
\label{ap10}
\end{equation}
The inclusion \eqref{ap8} implies that we may replace $\rho$ in \eqref{ap10}
by $\hat{\rho}$ since $\eta=1$ on $B_{13/21}$. 
Note that we have $R^2\bar{\sigma}^2\leq R^2\sigma^2\leq r_1^2/2$ by \eqref{ap4d}.
Due to \eqref{hhypsub}-\eqref{chika2sss} and \eqref{ap5a},
conditions \eqref{hhyp}-\eqref{chiu} are satisfied. Thus Lemma \ref{liphelp} gives
\begin{equation*}
2\leq (1+\frac{1}{100})\frac{26}{25}+\frac{4}{100}
\end{equation*}
which is a contradiction. Thus for any two distinct points $(y_1,s_1)$ and
$(y_2,s_2)$ in $A$, we have
\begin{equation}
|T^{\perp}(y_1-y_2)|\leq \max\{|T(y_1-y_2)|,|s_1-s_2|^{1/2}\}.
\label{finlip}
\end{equation}
We next consider the projection of $A$ on $T$. Define
\begin{equation}
A'=\{(x,s)\in B_{1/3}^k\times 
(1/4,3/4)\, :\,x=T(y)\mbox{ for some $(y,s)\in A$}\}.
\label{defAd}
\end{equation}
The inequality \eqref{finlip} shows that $(T^{-1}(x),s)\cap A$ 
consists of a single element $\{y\}$ for each $(x,s)\in A'$, thus we may
define the maps $\tilde{f}(x,s)=T^{\perp}(y)$ and $\tilde{F}(x,s)=(x,\tilde{f}(x,s),s)$.
We also have from \eqref{finlip}
and \eqref{apdef} that
\begin{equation*}
|{\tilde f}(x,s)-{\tilde f}(y,t)|\leq \max\{|x-y|,\, |s-t|^{1/2}\}
\hspace{.2cm}\mbox{ for $(x,s),\,(y,t)\in A'$},
\end{equation*}
\begin{equation*}
\sup_{(x,s)\in A'}|{\tilde f}(x,s)|\leq \varepsilon_3.
\end{equation*}
By the standard Lipschitz extension lemma applied with the natural
parabolic metric, we may extend ${\tilde f}$ and ${\tilde F}$ to
be defined on $T\times{\mathbb R}$. We denote them by $f$ and $F$,
respectively, and they satisfy
\begin{equation}
|f(x,s)-f(y,t)|\leq c(n,k) \max\{|x-y|,\, |s-t|^{1/2}\}
\hspace{.2cm}\mbox{ for $(x,s),\,(y,t)\in T\times{\mathbb R}$},
\label{tflip1}
\end{equation}
\begin{equation}
\sup_{(x,s)\in T\times{\mathbb R}}|f(x,s)|\leq \varepsilon_3.
\label{tflip2}
\end{equation}
This proves the claim of \eqref{lipok} by restricting to $C(T,1/3)$ from $C(T,1/2)$.
We next estimate the measures of 
\begin{equation}
\begin{split}
&B:=\cup_{t\in (1/4,
3/4)}\{(M_t\cap C(T,1/2))\times\{t\}\}\setminus A, \\
&B':=(T\times {\rm Id}_{\mathbb R})(B)\subset T\times{\mathbb R},\\
&B'':=B_{1/3}^k\times(1/4,3/4)\setminus A'.
\end{split}
\label{defb}
\end{equation}
$B$ is the subset of $\{M_t\}_{0\leq t\leq 1}$ which may not be covered by 
the image of $F$,
$B'$ is the projection to $T\times{\mathbb R}$, and $B''$ is the subset of
$T\times{\mathbb R}$ such that $T^{-1}(x)$ may not intersect with the
image of $F$. 

For each $(x,s)\in B$, there exists some $0<r(x,s)<r_1$ such that
\begin{equation}
\int_{\overline{P_{r(x,s)}(x,s)}}\|S-T\|^2\, dV_t(\cdot,S)dt\geq \gamma
(r(x,s))^{k+2}
\label{gam}
\end{equation}
by the definition of $A$, \eqref{apdef}. 
Thus $\{\overline{P_{r(x,s)}(x,s)}\}_{(x,s)\in B}$ is a covering of $B$.
By \cite[2.8.14]{Federer}, which is a generalized version of Besicovitch
covering theorem, there exists a finite number of subfamilies ${\mathcal B}_1,
\, \cdots, {\mathcal B}_{{\bf B}(n)}$ each of which consists of mutually disjoint
parabolic cylinders and $B\subset \cup_{i=1}^{{\bf B}(n)}\cup_{{\mathcal B}_i} \overline{P_{r(x,s)}(x,s)}$. Then using \eqref{dbd}, \eqref{gam} and the disjointness
of ${\mathcal B}_i$ we have
\begin{equation}
\begin{split}
&(\|V_t\|\times dt)(B)\leq \sum_{i=1}^{{\bf B}(n)}\sum_{{\mathcal B}_i}
(\|V_t\|\times dt)(\overline{P_{r(x,s)}(x,s)}) \\
& \leq \sum_{i=1}^{{\bf B}(n)}\sum_{{\mathcal B}_i}2E_1 \omega_k (r(x,s))^{k+2}
\leq
\sum_{i=1}^{{\bf B}(n)}\sum_{{\mathcal B}_i}2E_1\omega_k \gamma^{-1}\int_{
\overline{P_{r(x,s)}(x,s)}}\|S-T\|^2\, dV_t(\cdot,S)dt \\
&\leq 2{\bf B}(n)E_1\omega_k \gamma^{-1}\int_{C(T,13/24)\times(0,1)}\|S-T\|^2
\, dV_t(\cdot,S)dt.
\end{split}
\label{best}
\end{equation}
Note that $C(T,13/24)\subset\{\phi_T=1\}$ due to \eqref{pops1} and \eqref{Tphi} .

For $B'$, note that $T(P_r(x,s))\subset T\times{\mathbb R}$ is a 
$(k+1)$-dimensional parabolic cylinder with 
${\mathcal L}^{k+1}(T(P_r(x,s)))=
2\omega_k r^{k+2}$. Since $\cup_{i=1}^{{\bf B}(n)}\cup_{{\mathcal B}_i}T(P_{r(x,s)}
(x.s))$ covers $B'$, we have
\begin{equation}
{\mathcal L}^{k+1}(B')\leq \sum_{i=1}^{{\bf B}(n)}\sum_{{\mathcal B}_i}
{\mathcal L}^{k+1}(T(P_{r(x,s)}(x,s)))\leq \sum_{i=1}^{{\bf B}(n)}\sum_{{\mathcal B}_i}
2\omega_k (r(x,s))^{k+2}
\label{bbest}
\end{equation}
and the rest is estimated as in \eqref{best}, the result being without $E_1$ 
of \eqref{best}. 

We next estimate $B''$.
Since $f$ is Lipschitz in the space variables,
it is differentiable a.e$.$ and the area formula shows
(write $C_{*}:=C(T,1/2)\times(1/4,3/4)$ here for simplicity)
\begin{equation}
\begin{split}
\int_{C_*}\phi_{T,1/2}^2
\, d\|V_t\|dt& =\int_{A} \phi_{T,1/2}^2
\, d\|V_t\|dt + \int_{B}\phi_{T,1/2}^2
\, d\|V_t\|dt \\
&=\int_{A'}\phi_{T,1/2}^2(x)|\Lambda_k \nabla(x,f(x,t))|\, d{\mathcal H}^{k}(x)
dt+\int_{B}\phi_{T,1/2}^2
\, d\|V_t\|dt \\
&=\int_{A'}\phi_{T,1/2}^2(|\Lambda_k\nabla(x,f(x,t))|-1)\, d{\mathcal H}^k dt
+\int_{C_*\cap T}\phi_{T,1/2}^2\, d{\mathcal H}^k dt \\ &-\int_{(C_* \cap
T)\setminus
A'}\phi_{T,1/2}^2\, d{\mathcal H}^k dt+\int_B \phi_{T,1/2}^2
\, d\|V_t\|dt.
\end{split}
\label{cest}
\end{equation}
Since $\phi_{T,1/2}=1$ on $C(T,1/3)\times(1/4,3/4)$
and $B''\subset (C_*\cap T)\setminus A'$, 
\eqref{cest} shows
\begin{equation}
\begin{split}
&{\mathcal L}^{k+1}(B'')\leq \int_{(C_*\cap T)\setminus A'}\phi_{T,1/2}^2\, 
d{\mathcal H}^k dt 
 \leq \left|\int_{C_*}\phi_{T,1/2}^2\, d\|V_t\|dt-
\int_{C_*\cap T}\phi_{T,1/2}^2\, d{\mathcal H}^k dt\right| \\
&+\int_{A'}\phi_{T,1/2}^2(|\Lambda_k\nabla(x,f(x,s))|-1)\, d{\mathcal H}^k dt
+\int_B \phi_{T,1/2}^2\, d\|V_t\|dt.
\end{split}
\label{cest2}
\end{equation}
Note that for $f$ with its Lipschitz constant bounded by $c(n,k)$, we have
$|\Lambda_k \nabla (x,f(x,s))|-1\leq c(n,k)\|({\rm image}\, \nabla(x,f(x,s))
-T\|^2$, where ${\rm image}\, \nabla(x,f(x,s))$ coincides with the approximate
tangent space of $M_t$ a.e$.$ on $A'$. By \eqref{excesslip}, \eqref{best} and
\eqref{cest2}, we have
\begin{equation}
{\mathcal L}^{k+1}(B'')\leq 2{\bf B}(n)\omega_k E_1\gamma^{-1}\beta^2+
c(n,k)\beta^2+\kappa^2.
\label{cest3}
\end{equation}
Since $\cup_{t\in (1/4,3/4)}(M_t\times\{t\})\cap \{(C(T,1/3)\times(1/4,3/4)\}\setminus X\subset B$ 
and $\{B_{1/3}^k\times(1/4,3/4)\}\setminus Y
\subset (B'\cup B'')$, we proved \eqref{lipok2} by \eqref{best}, \eqref{bbest} and \eqref{cest3}.
\hfill{$\Box$}
\section{H\"{o}lder estimate for gradient}
In this section we prove that if the $L^2$-height of varifolds 
is small, then near the center of the domain,
the properly scaled $L^2$-height with respect to a slightly tilted plane shows at least a fixed
amount of decay. By iteration, this proves
that the spacial gradient is H\"{o}lder continuous. The method of proof, the blow-up technique, goes back
to Almgren \cite{Almgren} and the proof 
is analogous to that of Allard's regularity theory \cite{Allard}. 
On the other hand, there are a few subtle and interesting
differences from elliptic case. The first point is that $L^{\infty}$ estimate is rather essential to
show that the blow-up limit satisfies the heat equation. Note that the test function that we can
use in the weak inequality \eqref{maineq} has to be non-negative. Thus, we need to know 
the height of varifolds in $L^{\infty}$ norm, instead of, say, $L^2$ norm. The second point
is that we need to capitalize on the monotone decreasing property (with respect to time
variable) of $L^2$ norm of blow-up sequence to show the strong convergence of space-time $L^2$ norm. Recall that in elliptic case, Rellich's compactness theorem
shows the strong $L^2$ convergence immediately. 
Here, since we do not have any control of time-derivatives
of blow-up sequence, we need to use some special feature of parabolic problem.

\begin{prop}
Corresponding to $1\leq E_1<\infty$, $0<\nu<1$, $p$, $q$ with \eqref{expcond} there exist $0<\varepsilon_4<1$,
$2<\Lambda_*<\infty$, $0<\theta_*<1/4$, $1<c_{14}<\infty$ with the following property. For $W\in {\bf G}(n,k)$, $0<R<\infty$, $U=C(W,2R)$ and 
$(0,\Lambda)$ replaced by $(-\Lambda_* R^2,\Lambda_*R^2 )$, suppose $\{V_t\}_{-\Lambda_*R^2 \leq t\leq \Lambda_* R^2}$
and $\{u(\cdot,t)\}_{-\Lambda_*R^2 \leq t\leq \Lambda_*R^2}$ satisfy (A1)-(A4). 
For $W\in {\bf G}(n,k)$ let $\phi_{W,R}$ be
as in \eqref{Tphi}. 
Suppose 
\begin{equation}
T\in {\bf G}(n,k)\mbox{ satisfies }\|T-W\|<\varepsilon_4,
\label{blo0}
\end{equation}
\begin{equation}
\mbox{$A\in {\bf A}(n,k)$ is parallel to $T$,}
\label{blo0.5}
\end{equation}
\begin{equation}
\mu:=\left(R^{-(k+4)}\int_{-\Lambda_* R^2}^{\Lambda_* R^2}\int_{C(W,2R)}{\rm dist}\, (x,A)^2 \,
d\|V_t\|dt\right)^{\frac12}<\varepsilon_4,
\label{blo1}
\end{equation}
\begin{equation}
\|u\|:=\|u\|_{L^{p,q}(C(W,2R)\times(-\Lambda_*R^2,\Lambda_*R^2)) },
\label{blo1.5}
\end{equation}
\begin{equation}
(-\Lambda_*+1)R^2\leq \exists t_1\leq (-\Lambda_*+2)R^2\,\,:\,\, R^{-k}\|V_{t_1}\|(\phi_{W,R}^2)<
(2-\nu){\bf c},
\label{blo2}
\end{equation}
\begin{equation}
(\Lambda_*-2)R^2\leq \exists t_2\leq (\Lambda_*-1)R^2\,\,:\,\, R^{-k}\|V_{t_2}\|(\phi_{W,R}^2)>
\nu{\bf c}.
\label{blo3}
\end{equation}
Then there are $\tilde{T}\in {\bf G}(n,k)$ and $\tilde{A}\in {\bf A}(n,k)$ such that
\begin{equation}
\mbox{$\tilde{A}$ is parallel
to $\tilde{T}$, }
\label{imblo1}
\end{equation}
\begin{equation}
\|T-\tilde{T}\|\leq c_{14} \mu,
\label{imblo2}
\end{equation}
\begin{equation}
\left((\theta_*R)^{-(k+4)}\int_{-\theta_*^2 \Lambda_*R^2}^{\theta_*^2\Lambda_*R^2}
\int_{C(W,2\theta_*R)}{\rm dist}\, (x,\tilde{A})^2\, d\|V_t\|dt\right)^{\frac12} \leq
\theta_*^{\varsigma}\max\{\mu, R^{\varsigma}c_{14}\|u\|\}.
\label{imblo3}
\end{equation}
Moreover, if $R^{\varsigma}\|u\|<\varepsilon_4$, we have
\begin{equation}
(-\Lambda_*+1)(\theta_*R)^2\leq \exists \tilde{t}_1\leq (-\Lambda_*+2)(\theta_*R)^2
\,\,:\,\, (\theta_*R)^{-k}\|V_{\tilde{t}_1}\|(\phi_{W,\theta_*R}^2)<
(2-\nu){\bf c},
\label{blo2ag}
\end{equation}
\begin{equation}
(\Lambda_*-2)(\theta_*R)^2\leq \exists \tilde{t}_2\leq (\Lambda_*-1)(\theta_*R)^2
\,\,:\,\, (\theta_*R)^{-k}\|V_{\tilde{t}_2}\|(\phi_{W,\theta_*R}^2)>
\nu{\bf c}.
\label{blo3ag}
\end{equation}
\label{blowup1}
\end{prop}
{\it Proof}. We may assume that $R=1$ after a suitable change of variables.
We prove \eqref{imblo1}-\eqref{imblo3} by contradiction. 
We will show that 
\eqref{blo2ag} and \eqref{blo3ag}
follow from these at the end. 
By Theorem \ref{poptheorem}, corresponding to $E_1$ and $\nu$ there replaced by 
$\nu/2$, we fix
$\Lambda$, $K$, $\varepsilon_1$ first. We set $\Lambda_*:=\Lambda+5/2$.
We will fix $0<\theta_*<1/4$
later depending only on $E_1$, $\nu$, $p$ and $q$. 
If the claim were false, then for each $m\in {\mathbb N}$ there exist
$\{V^{(m)}_t\}_{-\Lambda_*\leq t\leq \Lambda_*}$, $\{u^{(m)}(\cdot,t)\}_{-\Lambda_*
\leq t\leq \Lambda_*}$ satisfying (A1)-(A4) on $C(W^{(m)},2)\times [-\Lambda_*,\Lambda_*]$ for $W^{(m)}\in {\bf G}(n,k)$ such that, after suitable translations and rotations,
we have $T={\mathbb R}^k\times\{0\}$ and 
\begin{equation}
\|T-W^{(m)}\|\leq \frac{1}{m},
\label{dame0}
\end{equation}
\begin{equation}
\mu^{(m)}:=\left(\int_{-\Lambda_*}^{\Lambda_*}\int_{C(W^{(m)},2)} |T^{\perp}(x)|^2
\, d\|V_t^{(m)}\|dt\right)^{\frac12}\leq \frac{1}{m},
\label{dame}
\end{equation}
and \eqref{blo2} and \eqref{blo3} satisfied for $W=W^{(m)}$, $t_1=t_1^{(m)}$ and $t_2=t_2^{(m)}$,
respectively, 
but for any $\tilde{T}\in {\bf G}(n,k)$ with $\|T-\tilde{T}\|\leq m\mu^{(m)}$ and $\tilde{A}\in
{\bf A}(n,k)$ which is parallel to $\tilde{T}$, 
\begin{equation}
\left(\theta_*^{-(k+4)}\int_{-\theta_*^2 \Lambda_*}^{\theta_*^2 \Lambda_*}\int_{C(W^{(m)},
2\theta_*)}{\rm dist}\, (x,\tilde{A})^2\, 
d\|V_t^{(m)}\|dt\right)^{\frac12}
> \theta_*^{\varsigma}\max\{\mu^{(m)}, m \|u^{(m)}\|\}.
\label{blo1no}
\end{equation}
By taking $\tilde{A}=\tilde{T}=T$ in \eqref{blo1no}, we obtain
\begin{equation}
\theta^{\varsigma}_* m\|u^{(m)}\|
<\theta_*^{-(k+4)/2}\mu^{(m)}.
\label{blo1no1}
\end{equation}
Thus \eqref{blo1no1} shows that 
\begin{equation}
\lim_{m\rightarrow\infty} (\mu^{(m)})^{-1}\|u^{(m)}\|=0.
\label{mudaex}
\end{equation} 
\begin{lemma}
For any given $0<\gamma<1$, there exists some $m_0$ such that 
\begin{equation}
{\rm spt}\, \|V^{(m)}_t\|\cap \{|T^{\perp}(x)|>\gamma\}\cap C(W^{(m)},3/2)=\emptyset,\hspace{.3cm}
-\Lambda_*+1\leq \forall t\leq \Lambda_*
\label{emp}
\end{equation}
for all $m\geq m_0$. 
\label{heiblo}
\end{lemma}
{\it Proof of Lemma \ref{heiblo}}.
Suppose that \eqref{emp} did not hold for some $\tilde{t}$ such that
$V^{(m)}_{\tilde{t}}$ is of unit density. Then there exists some $\tilde{x}\in \{|T^{\perp}(x)|>\gamma\}\cap
C(W^{(m)},3/2)$ with $\Theta^k(\|V^{(m)}_{\tilde{t}}\|,\tilde{x})=1$. Then by Proposition \ref{cam1} with
$U=B_{\gamma/2}(\tilde{x})$, $R=\gamma/2$, $t_2=\tilde{t}$, $t_1=\tilde{t}-t$ with
$0<t<1$, $y=\tilde{x}$ and
$s=\tilde{t}+\epsilon$, we have
\begin{equation}
\begin{split}
\int_{B_{\gamma/2}(\tilde{x})}\hat{\rho}_{(\tilde{x},\tilde{t}+\epsilon)}(\cdot,\tilde{t})\, d\|V^{(m)}_{\tilde{t}}\|
&\leq \int_{B_{\gamma/2}(\tilde{x})}\hat{\rho}_{(\tilde{x},\tilde{t}+\epsilon)}(\cdot, \tilde{t}-t)\,
d\|V^{(m)}_{\tilde{t}-t}\| \\
&+c_2 \|u^{(m)}\|^2E_1^{1-\frac{2}{p}}t^{\varsigma}+4\gamma^{-2}c_1\o_k
E_1 t.
\end{split}
\label{emp1}
\end{equation}
Set $\epsilon\rightarrow 0+$ and fix $0<t'<1$ such that the 
sum of the second and third term of the right-hand side of \eqref{emp1} is
less than $1/2$ when $t\leq t'$. We then obtain
\begin{equation}
\frac12\leq \int_{B_{\gamma/2}(\tilde{x})}\hat{\rho}_{(\tilde{x},\tilde{t})}(\cdot,
\tilde{t}-t)\, d\|V_{\tilde{t}-t}^{(m)}\|
\label{emp2}
\end{equation}
for all $t\in [t'/2,t']$. Since $\hat{\rho}_{(\tilde{x},\tilde{t})}(y,\tilde{t}-t)$ takes
maximum value $(4\pi t')^{-k/2}$ when $y=\tilde{x}$ and $t=t'/2$ in this interval, 
from \eqref{emp2} we
obtain
\begin{equation}
\frac12\leq (4\pi t')^{-k/2}\|V^{(m)}_{\tilde{t}-t}\|(B_{\gamma/2}(\tilde{x}))
\label{emp3}
\end{equation}
for all $t\in [t'/2,t']$. Since $B_{\gamma/2}(\tilde{x})\subset
\{|T^{\perp}(x)|>\gamma/2\}\cap C(W^{(m)},2)$, by \eqref{emp3}, we have
\begin{equation}
(\mu^{(m)})^2\geq \int_{\tilde{t}-t'}^{\tilde{t}-t'/2}\int_{B_{\gamma/2}(\tilde{x})}
|T^{\perp}(x)|^2\, d\|V^{(m)}_t\|dt\geq \frac{t'}{2}\left(\frac{\gamma}{2}\right)^2
\frac{(4\pi t')^{k/2}}{2}.
\label{emp4}
\end{equation}
Since $\mu^{(m)}\leq 1/m$, \eqref{emp4} cannot be true for sufficiently
large $m$. This shows \eqref{emp} holds true for a.e$.$ $t\in [-\Lambda_*
+1,\Lambda_*]$. By using \eqref{maineq}, we conclude that \eqref{emp}
is satisfied for all $t$. This ends the proof of Lemma \ref{heiblo}.
\hfill{$\Box$}

Next we use
Proposition \ref{el2elinf} for $U=B_2$ (or $R=2$ there) and $(0,\Lambda)$ there
replaced by $(-\Lambda_*, \Lambda_*)$, to obtain 
\begin{equation}
{\rm spt}\, \|V_t^{(m)}\|\cap B_{8/5}\subset\{|T^{\perp}(x)|\leq 2\tilde{\mu}^{(m)}\},\hspace{.5cm}
\forall t\in (-\Lambda_*+1, \Lambda_*)
\label{blo10}
\end{equation}
for 
\begin{equation}
\begin{split}
\tilde{\mu}^{(m)}:&=\left(\frac{c_3}{2^{k+4}}\int_{-\Lambda_*}^{\Lambda_*}\int_{B_2}
|T^{\perp}(x)|^2 \, d\|V_t^{(m)}\|dt+c_2\|u^{(m)}\|^2
E_1^{1-\frac{2}{p}}(2\Lambda_*)^{\varsigma}(2+\Lambda_*/2)\right)^{\frac12} \\
&\leq c_8(n,k)
\mu^{(m)}
\end{split}
\label{blo11}
\end{equation}
for all sufficiently large $m$ by \eqref{mudaex} for a suitable $c_8$. 
In the following, for $t\in (-\Lambda_*+1,\Lambda_*)$, we set $\|V^{(m)}_t\|=0$
on $C(T,1)\cap\{|T^{\perp}(x)|\geq 1\}$ for all sufficiently large $m$. Due to
Lemma \ref{heiblo} and \eqref{blo10}, this modified $V_t^{(m)}$ on $C(T,1)$
still satisfies (A1)-(A4). We do this modification for notational simplicity.
Otherwise, one can restrict the domain of integration from $C(T,1)$ to $C(T,1)
\cap \{|T^{\perp}(x)|<1\}$ in the following computations. With this convention, we have
\begin{equation}
\begin{split}
\mu_*^{(m)}:& =\sup_{-\Lambda_*+1\leq t\leq \Lambda_*-1} 
\left(\int_{C(T,1)}|T^{\perp}(x)|^2\, d\|V_t^{(m)}\|\right)^{\frac12}\leq (E_1\o_k (8/5)^k)^{1/2} 2\tilde{\mu}^{(m)} \\
& \leq c_9(E_1,n,k) \mu^{(m)}
\end{split}
\label{blo12}
\end{equation}
for a suitable $c_9$ by \eqref{blo10} and \eqref{blo11}. We next check that the 
assumptions for Theorem \ref{poptheorem} are all satisfied on $(-\Lambda_*+1,
\Lambda_*-1)$. By the definition of $\Lambda_*$, we have $(\Lambda_*-1)-
(-\Lambda_*+1)=2\Lambda_*-2=2\Lambda+3$. \eqref{pops31.1} and \eqref{pops31.2}
are satisfied for all sufficiently large $m$ 
due to \eqref{blo2}, \eqref{blo3} and \eqref{dame0}. \eqref{pops30.5} follows
from \eqref{blo10}. Thus for all sufficiently large $m$ depending on $\varepsilon_1$ so that \eqref{pops31}
and \eqref{pops31.3} are satisfied, we conclude from \eqref{pops31.4} and
\eqref{popsmean} as well as Corollary \ref{popscor} that 
\begin{equation}
\sup_{-1/2\leq t\leq 1/2}|\|V^{(m)}_t\|(\phi_T^2)-{\bf c}|\leq K((\mu_*^{(m)})^2+C(u^{(m)})),
\label{blo13}
\end{equation}
\begin{equation}
\int_{-1/2}^{1/2}\int_{C(T,1)}|h(V_t^{(m)},\cdot)|^2\phi_T^2\, d\|V^{(m)}_t\|dt\leq 12 K((\mu_*^{(m)})^2
+C(u^{(m)})),
\label{blo14}
\end{equation}
\begin{equation}
\int_{-1/2}^{1/2}|2^k \|V^{(m)}_t\|(\phi_{T,1/2}^2)-{\bf c}|\, dt\leq \tilde{K}((\mu_*^{(m)})^2 +C(u^{(m)})).
\label{blo15}
\end{equation}
The right-hand sides of \eqref{blo13}-\eqref{blo15} are all bounded by $c_{10}(E_1,\nu,n,k)(\mu^{(m)})^2$
due to \eqref{blo12}, the definition of $C(u^{(m)})$ and \eqref{mudaex} for all sufficiently larege $m$.  
By Proposition 
\ref{tiltexlem}, \eqref{blo12} and \eqref{blo14}, we also have
\begin{equation}
\int_{-1/2}^{1/2}\int_{G_k(C(T,1))}\|S-T\|^2 \phi_T^2\, dV_t^{(m)}(\cdot,S)dt\leq c_{10}(\mu^{(m)})^2.
\label{blo15.5}
\end{equation} 
We next use Theorem \ref{lipapp}. For all sufficiently large $m$ depending
on $r_1$ and $\varepsilon_3$, we may satisfy the assumptions 
\eqref{hhypsub}-\eqref{chika2sss}.
Since $\beta^2$ and $\kappa^2$ in Theorem \ref{lipapp} are bounded by $c_{10} 
(\mu^{(m)})^2$,
we have a corresponding Lipschitz functions $f^{(m)}$ and $F^{(m)}$
defined on $B_{1/3}^k\times(-1/4,1/4)$ with the corresponding estimates
\eqref{lipok} and \eqref{lipok2}. We also define $X^{(m)}$ and $Y^{(m)}$ to be 
the sets $X$ and $Y$ respectively corresponding to $V_t^{(m)}$. 
We now set
\begin{equation}
\tilde{f}^{(m)}:=f^{(m)}/\mu^{(m)}.
\label{tilf}
\end{equation}
In the following we write $\Omega:=C(T,1/3)\times (-1/4,1/4)$ and 
$\Omega':=B_{1/3}^k\times(-1/4,1/4)$. 
\begin{lemma}
There exists constant $c_{11}$ which depends only on $E_1$, $\nu$, $p$, $q$ such that 
for all sufficiently large $m$,
\label{el2}
\begin{equation}
\sup_{\Omega'}|\tilde{f}^{(m)}|+\int_{\Omega'}|\nabla\tilde{f}^{(m)}|^2\,
d{\mathcal H}^{k+1}\leq c_{11}.
\label{tilfest}
\end{equation}
\label{tilfestlem}
\end{lemma}
{\it Proof of Lemma \ref{tilfestlem}}. The supremum estimate follows immediately from \eqref{blo10}, 
\eqref{blo12}, \eqref{lipok} and
\eqref{tilf}. We may proceed just like \cite[p.483]{Allard} for the second
term. We split the domain of integration into $Y^{(m)}$ and the complement,
and using that the spacial Lipschitz constants of $F^{(m)}$ are $\leq c(n,k)$, 
\begin{equation}
\begin{split}
&\int_{Y^{(m)}}|\nabla f^{(m)}|^2\, d{\mathcal H}^{k+1}+
\int_{\Omega'\setminus Y^{(m)}}|\nabla f^{(m)}|^2\, d{\mathcal H}^{k+1} \\
&\leq c(n,k)\int_{Y^{(m)}}|\nabla F^{(m)}|^2 \|{\rm image}\, \nabla F^{(m)}-T\|^2
|\Lambda_k \nabla F^{(m)}|\, d{\mathcal H}^{k+1}+c(n,k){\mathcal L}^{k+1}(
\Omega'\setminus Y^{(m)})\\
&\leq c(n,k)\int_{G_k(C(T,1/3))\times(-1/4,1/4)}\|S-T\|^2\, dV^{(m)}_t(\cdot,S)dt
+c(n,k){\mathcal L}^{k+1}(\Omega'\setminus Y^{(m)}) \\
& \leq c(n,k)c_{10}(\mu^{(m)})^2
\end{split}
\label{tilfest1}
\end{equation}
where we used \eqref{blo15.5} and \eqref{lipok2} at the end of \eqref{tilfest1}.
This shows \eqref{tilfest}.
\hfill{$\Box$}
\begin{lemma}
There exist a subsequence $\{\tilde{f}^{(m_j)}\}_{j=1}^{\infty}$ and 
$\tilde{f}\in C^{\infty}(\Omega')$ such that 
\newline
$\lim_{j\rightarrow\infty}\|\tilde{f}^{(m_j)}-\tilde{f}\|_{L^2(\Omega')}=0$ and $\tilde{f}$ satisfies 
$\frac{\partial \tilde{f}}{\partial t}-\Delta \tilde{f}=0$ on $\Omega'$.
\label{heatsol}
\end{lemma}
{\it Proof of Lemma \ref{heatsol}}. Let $l\in
\{1,\cdots,n-k\}$ be fixed. We show that $l$-th components $\tilde{f}^{(m)}_l,\,
\tilde{f}_l$ of $\tilde{f}^{(m)},\, \tilde{f}\in {\mathbb R}^{n-k}$
respectively satisfy the desired properties. Due to \eqref{tilfest} 
we can extract a subsequence (denoted by the same index) 
so that 
\begin{equation}
\tilde{f}^{(m)}_l\rightharpoonup \tilde{f}_l\mbox{ weakly
in $L^2(\Omega')$},\hspace{.5cm}\nabla\tilde{f}^{(m)}_l\rightharpoonup 
\nabla\tilde{f}_l\mbox{ weakly in $L^2(\Omega')$}.
\label{weak1}
\end{equation}
Due to \eqref{blo15.5} and \eqref{lipok2}, we may further assume that for a.e$.$
$s\in (-1/4,1/4)$
\begin{equation}
\lim_{m\rightarrow \infty}\int_{C(T,1/3)}\|S-T\|^2\, dV^{(m)}_s(\cdot,S)=0,
\label{weak3}
\end{equation}
\begin{equation}
\lim_{m\rightarrow\infty}
\|V^{(m)}_s\|(C(T,1/3)\setminus X^{(m)}\lfloor_{t=s})=0,\hspace{.5cm}
\lim_{m\rightarrow\infty}{\mathcal L}^k(B_{1/3}^k\setminus Y^{(m)}\lfloor_{t=s})=0.
\label{weak4}
\end{equation}
Again by \eqref{tilfest}, the Rellich compactness theorem and diagonal argument, 
we may assume that 
\begin{equation}
\{\tilde{f}_l^{(m)}(\cdot,s_j)\}_{m\in {\mathbb N}}\mbox{ is a Cauchy sequence
in $L^2(B_{1/3}^k)$}
\label{weak2}
\end{equation}
for a countable dense set $\{s_i\}_{i=1}^{\infty}$ in $(-1/4,1/4)$.
We may choose these $s_i$ satisfying \eqref{weak3} and \eqref{weak4} as well.
We fix such a subsequence in the following.
For arbitrary non-negative $\tilde{\phi}\in C^{\infty}_c(\Omega')$ we set $\phi^{(m)}(x)
:=(x_l+c_{11}\mu^{(m)})\tilde{\phi}(T(x))$ for $x\in {\mathbb R}^n$. 
By definition $\phi^{(m)}$ is non-negative on 
$\{x\, :\, |x_l|\leq c_{11}\mu^{(m)}\}$. We may use $\phi^{(m)}$ 
as a test function in \eqref{maineq} by a slight modification away from $T$.
In the following similar computations using \eqref{maineq} we implicitly assume 
that we do this modification of test functions which does not affect the computations. Now
we have
\begin{equation}
0\leq \int_{\Omega}(-h(V_t^{(m)},\cdot)\phi^{(m)}+\nabla\phi^{(m)})
\cdot(h(V_t^{(m)},\cdot)+(u^{(m)})^{\perp})+\frac{\partial \phi^{(m)}}{
\partial t}\, d\|V_t^{(m)}\|dt.
\label{tai1}
\end{equation}
By the Cauchy-Schwarz inequality and by dropping a negative term, we obtain from \eqref{tai1}
\begin{equation}
\begin{split}
0\leq &\int_{\Omega}|u^{(m)}|^2\phi^{(m)}+|u^{(m)}||\nabla\phi^{(m)}|
+\frac{\partial \phi^{(m)}}{\partial t}
+(x_l+c_{11}\mu^{(m)})\nabla\tilde{\phi}\cdot h(V_t^{(m)},\cdot) \\ &+\tilde{\phi}
h(V_t^{(m)},\cdot)_l\, d\|V_t^{(m)}\|dt
=: I_1^{(m)}+I_2^{(m)}+I_3^{(m)}+I_4^{(m)}+I_5^{(m)}
\end{split}
\label{tai2}
\end{equation}
where $h(V_t^{(m)},x)_l$ is the $(k+l)$-th component of $h(V_t^{(m)},x)\in {\mathbb R}^n$.
In the following we estimate $\lim_{m\rightarrow \infty}(\mu^{(m)})^{-1}I_j^{(m)}$ for each $j=1,\cdots,5$.
\newline
{\bf Estimate of $I_1^{(m)}$}.
\newline
Since $|\phi^{(m)}|\leq 2c_{11}\mu^{(m)} \sup|\tilde{\phi}|$, by the H\"{o}lder
inequality and \eqref{mudaex},
\begin{equation}
\lim_{m\rightarrow\infty}|(\mu^{(m)})^{-1}I_1^{(m)}|\leq c_{11}c(\tilde{\phi},p,q)
 \lim_{m\rightarrow\infty}\|u^{(m)}\|^2=0.
\label{tai3}
\end{equation}
\newline
{\bf Estimate of $I_2^{(m)}$}.
\newline
By the H\"{o}lder inequality and \eqref{mudaex},
\begin{equation}
\begin{split}
\lim_{m\rightarrow \infty}|(\mu^{(m)})^{-1}I_2^{(m)}|
& \leq \lim_{m\rightarrow \infty} c(p,q)\|u^{(m)}\|
(\sup|\tilde{\phi}|(\mu^{(m)})^{-1} 
+c_{11}\sup |\nabla\tilde{\phi}|) \\
&\leq \lim_{m\rightarrow \infty}o(1)c(\tilde{\phi},p,q)(1+c_{11}\mu^{(m)})=0.
\end{split}
\label{taitai}
\end{equation}
\newline
{\bf Estimate of $I_3^{(m)}$}.
\newline
We estimate as
\begin{equation}
\begin{split}
\int_{\Omega} & \frac{\partial\phi^{(m)}}{\partial t}\, d\|V_t^{(m)}\|dt \\
&=\int_{\Omega\setminus X^{(m)}}\frac{\partial\phi^{(m)}}{\partial t}
\, d\|V_t^{(m)}\|dt+\int_{Y^{(m)}}(f^{(m)}_l+c_{11}\mu^{(m)})\frac{\partial
\tilde{\phi}}{\partial t}(|\Lambda_k \nabla F^{(m)}|-1)\, d{\mathcal H}^{k+1} \\
& -\int_{\Omega'\setminus Y^{(m)}}(f^{(m)}_l+c_{11}\mu^{(m)})\frac{\partial\tilde{\phi}}
{\partial t}\, d{\mathcal H}^{k+1}
+\int_{\Omega'}(f^{(m)}_l+c_{11}\mu^{(m)})\frac{\partial\tilde{\phi}}{\partial t}
\, d{\mathcal H}^{k+1}.
\end{split}
\label{tai4}
\end{equation}
Using \eqref{tilfest}, \eqref{lipok2} and \eqref{blo15.5}, we can show that the 
first three terms of \eqref{tai4} are of order $(\mu^{(m)})^3$. 
Thus, using the weak convergence \eqref{weak1}, we obtain
\begin{equation}
\lim_{m\rightarrow\infty}(\mu^{(m)})^{-1}I_3^{(m)}=\int_{\Omega'}(\tilde{f}_l+c_{11})
\frac{\partial\tilde{\phi}}{\partial t}\, d{\mathcal H}^{k+1}.
\label{tai5}
\end{equation}
\newline
{\bf Estimate of $I_4^{(m)}$}.
\newline
By \eqref{tilfest} and \eqref{blo14}, we have
\begin{equation*}
\begin{split}
|I_4^{(m)}|&\leq 2c_{11}\mu^{(m)}\left(\int_{\Omega}|\nabla\tilde{\phi}|^2\, d\|V_t^{(m)}\|
dt\right)^{1/2}\left(\int_{\Omega}|h(V_t^{(m)},\cdot)|^2\, d\|V_t^{(m)}\|dt
\right)^{1/2} \\ &\leq 2c_{11}c(E_1,\tilde{\phi}) c_{10}^{1/2}(\mu^{(m)})^2
\end{split}
\end{equation*}
so that 
\begin{equation}
\lim_{m\rightarrow\infty}|(\mu^{(m)})^{-1}I_4^{(m)}|=0.
\label{tai6}
\end{equation}
\newline
{\bf Estimate of $I_5^{(m)}$}.
\newline
Let ${\bf e}_l$ be the unit vector with 1 in the $l$-th component. We  may write
\begin{equation}
\begin{split}
&\int_{\Omega}\tilde{\phi}h(V_t^{(m)},\cdot)_l\, d\|V_t^{(m)}\|dt =
-\int_{\Omega}(\nabla\tilde{\phi}\otimes {\bf e}_l)\cdot S\, dV_t^{(m)}dt \\
&=-\int_{\Omega\setminus X^{(m)}}(\nabla\tilde{\phi}\otimes {\bf e}_l)\cdot S\, dV_t^{(m)}dt\\
&-\int_{Y^{(m)}}{\rm image}\, \nabla F^{(m)}\cdot (\nabla\tilde{\phi}\otimes{\bf e}_l)|\Lambda_k
\nabla F^{(m)}|-\nabla f_l^{(m)}\cdot \nabla\tilde{\phi}\, d{\mathcal L}^{k+1} \\
&+\int_{\Omega'\setminus Y^{(m)}}\nabla f_l^{(m)}\cdot \nabla\tilde{\phi}\, d{\mathcal L}^{k+1}
-\int_{\Omega'}\nabla f_l^{(m)}\cdot \nabla\tilde{\phi}\, d{\mathcal L}^{k+1}:=J_1^{(m)}+
J_2^{(m)}+J_3^{(m)}+J_4^{(m)}.
\end{split}
\label{fvto}
\end{equation}
By \eqref{lipok2} $|J_1^{(m)}|$ is $O((\mu^{(m)})^2)$ and by \eqref{tilfest} and 
\eqref{lipok2}, the same for $|J_3^{(m)}|$.  
By \cite[8.14]{Allard}, \eqref{lipok} and \eqref{tilfest},
\begin{equation}
|J_2^{(m)}|\leq c(n,k,|\nabla\tilde{\phi}|)\int_{Y^{(m)}}|\nabla f_l^{(m)}|^2\, d{\mathcal L}^{k+1}
\leq O((\mu^{(m)})^2).
\label{fvto2}
\end{equation}
By \eqref{weak1}, it follows from above estimates that we have
\begin{equation}
\lim_{m\rightarrow\infty}(\mu^{(m)})^{-1}I_5^{(m)}=\lim_{m\rightarrow\infty}
(\mu^{(m)})^{-1}J_4^{(m)}=-\int_{\Omega'}
\nabla\tilde{f}_l\cdot\nabla\tilde{\phi}\, d{\mathcal H}^{k+1}.
\label{tai7}
\end{equation}
\newline
{\bf Summary}.
\newline
Due to \eqref{tai2}-\eqref{tai7}, we obtain
\begin{equation}
0\leq \int_{\Omega'}(\tilde{f}_l+c_{11})\frac{\partial \tilde{\phi}}{\partial t}
-\nabla\tilde{f}_l\cdot\nabla\tilde{\phi}\, d{\mathcal H}^{k+1}.
\label{tai8}
\end{equation}
We may carry out the similar argument with $\phi=(c_{11}\mu^{(m)}-x_l)\tilde{\phi}$,
and obtain
\begin{equation}
0\leq \int_{\Omega'}(c_{11}-\tilde{f}_l)\frac{\partial\tilde{\phi}}{\partial t}
+\nabla\tilde{f}_l\cdot\nabla\tilde{\phi}\, d{\mathcal H}^{k+1}.
\label{tai9}
\end{equation}
Combining \eqref{tai8} and \eqref{tai9} and since $\int_{\Omega'}\frac{ \partial\tilde{\phi}}
{\partial t}\, d{\mathcal H}^{k+1}=0$
, we conclude that 
\begin{equation}
0=\int_{\Omega'}\tilde{f}_l\frac{\partial\tilde{\phi}}{\partial t}
-\nabla\tilde{f}_l\cdot\nabla\tilde{\phi}\, d{\mathcal H}^{k+1}
\label{tai10}
\end{equation}
for all non-negative $\tilde{\phi}\in C^{\infty}_c(\Omega')$. 
This also shows that \eqref{tai10} is satisfied for all (not necessarily
non-negative) $\tilde{\phi}\in C^{\infty}_c(\Omega')$, and  $\tilde{f}_l$ 
satisfies the heat equation in a weak sense. The standard regularity
theory shows that $\tilde{f}_l$ is $C^{\infty}(\Omega')$ satisfying
the heat equation in the classical sense. 
\newline
Next we show that for a.e$.$ $s\in (-1/4,1/4)$, the weak convergent limit
is unique and
\begin{equation}
\tilde{f}^{(m)}_l(\cdot,s)\rightharpoonup 
\tilde{f}_l(\cdot,s)\mbox{ 
weakly in $L^2(B_{1/3}^k)$.}
\label{yweak}
\end{equation}
Take any $s$ with the properties
\eqref{weak3} and \eqref{weak4}. Since $|\tilde{f}^{(m)}_l(\cdot,s)|\leq c_{11}$, it
is bounded in particular in $L^2(B_{1/3}^k)$. Let $g\in L^2(B_{1/3}^k)$ be
any weak limit of a subsequence $\{\tilde{f}^{(m_j)}_l(\cdot,s)\}_{j=1}^{\infty}$.
The similar estimate as in \eqref{tai4} using \eqref{weak3} and \eqref{weak4} 
shows that for any $\tilde{\phi}\in C^{\infty}(B_{1/3}^k)$
\begin{equation}
\lim_{j\rightarrow\infty}\int_{C(T,1/3)}(\mu^{(m_j)})^{-1}
(x_l+c_{11}\mu^{(m_j)})\tilde{\phi}\, d\|V_{s}^{m_j}\|
=\int_{B_{1/3}^k}(g+c_{11})\tilde{\phi}\, d{\mathcal H}^k.
\label{fat1}
\end{equation}
We now proceed just as the first part of the proof. The differences 
this time are that
the domain of integration is changed to $C(T,1/3)\times(-1/4,s)$ and
that we take the subsequence $\{\tilde{f}^{(m_j)}_l\}_{j=1}^{\infty}$.
Then using \eqref{fat1}, we obtain for any non-negative $\tilde{\phi}
\in C^{\infty}_c(\Omega')$ 
\begin{equation}
\int_{B_{1/3}^k}(g+c_{11})\tilde{\phi}(\cdot,s)\, d{\mathcal H}^k
\leq \int_{\Omega'\cap \{t\leq s\}}(\tilde{f}_l+c_{11})\frac{\partial\tilde{\phi}}{
\partial t}-\nabla\tilde{f}_l\cdot\nabla\tilde{\phi}\, d{\mathcal H}^k
=\int_{\Omega'\cap\{t=s\}}(\tilde{f}_l+c_{11})\tilde{\phi}\, d{\mathcal H}^k,
\label{gdt1}
\end{equation}
the last equality follows from $\tilde{f}_l$ being the classical solution for
the heat equation. Similarly, we have
\begin{equation}
\int_{B_{1/3}^k}(c_{11}-g)\tilde{\phi}(\cdot,s)\, d{\mathcal H}^k
\leq \int_{\Omega'\cap\{t\leq s\}}(c_{11}-\tilde{f}_l)\frac{\partial\tilde{\phi}}{
\partial t}+\nabla\tilde{f}_l\cdot\nabla\tilde{\phi}\, d{\mathcal H}^k
=\int_{\Omega'\cap\{t=s\}}(c_{11}-\tilde{f}_l)\tilde{\phi}\, d{\mathcal H}^k.
\label{gdt2}
\end{equation}
Thus \eqref{gdt1} and \eqref{gdt2} show $\int_{B_{1/3}^k}(g-\tilde{f}_l(\cdot,s))\tilde{\phi}(\cdot,s)\, d{\mathcal H}^k=0$
for all non-negative function $\tilde{\phi}\in C^{\infty}_c(\Omega')$, which shows that
$g(\cdot)=\tilde{f}_l(\cdot,s)$ a.e$.$ on $B_{1/3}^k$. Since the limit is determined 
independent of the choice of subsequence, the whole sequence 
$\{\tilde{f}^{(m)}_l(\cdot,s)\}_{m=1}^{\infty}$ must converge to $\tilde{f}_l
(\cdot,s)$ weakly in $L^2(B_{1/3}^k)$, proving \eqref{yweak}. The same argument
shows that the $L^2$ Cauchy sequence in \eqref{weak2} also
converges to $\tilde{f}_l(\cdot,s_j)$. 
The lower semicontinuity under weak convergence shows that
\begin{equation}
\|\tilde{f}_l(\cdot,s)\|_{L^2(B_{1/3}^k)}\leq \liminf_{m\rightarrow\infty}
\|\tilde{f}_l^{(m)}(\cdot,s)\|_{L^2(B_{1/3}^k)}
\label{lsc1}
\end{equation}
for a.e$.$ $s\in (-1/4,1/4)$. We next show that for any $-1/4<s_j<s<1/4$
with $s_j$ satisfying \eqref{weak2}, $s$ satisfying \eqref{weak3} and
\eqref{weak4}, and for any non-negative $\tilde{\phi}\in C^{\infty}_c(B_{1/3}^k)$, 
\begin{equation}
\limsup_{m\rightarrow\infty}\|\tilde{\phi}\tilde{f}_l^{(m)}(\cdot,s)\|^2_{L^2(B_{1/3}^k)}
\leq \|\tilde{\phi}\tilde{f}_l(\cdot,s_j)\|^2_{L^2(B_{1/3}^k)}+c(\tilde{\phi})(s-s_j)^{1/2}.
\label{lsc2}
\end{equation}
We use $(x_l)^2\tilde{\phi}$ as a test function in \eqref{maineq} with time
interval $[s_j,s]$ and divide both sides by $(\mu^{(m)})^2$. 
By the Cauchy-Schwarz inequality and abbreviating the notations,
substitution of $(x_l)^2\tilde{\phi}$ into \eqref{maineq} gives
\begin{equation}
\begin{split}
&\int_{C(T,1/3)}(x_l)^2\tilde{\phi}\, d\|V_s^{(m)}\|-
\int_{C(T,1/3)}(x_l)^2\tilde{\phi}\, d\|V_{s_j}^{(m)}\| \\
&\leq \int_{C(T,1/3)\times(s_j,s)}(-h(V_t^{(m)},\cdot)(x_l)^2\tilde{\phi}
+\nabla((x_l)^2\tilde{\phi}))\cdot(h(V_t^{(m)},\cdot)+(u^{(m)})^{\perp})\,
d\|V_t^{(m)}\|dt \\
&\leq \int_{C(T,1/3)\times(s_j,s)}\nabla((x_l)^2\tilde{\phi})\cdot(h(V_t^{(m)},\cdot)
+(u^{(m)})^{\perp})+|u^{(m)}|^2(x_l)^2\tilde{\phi}\, d\|V_t^{(m)}\|dt\\
&\leq c(\tilde{\phi})\left\{\left(\int(x_l)^2\right)^{1/2}+\left(\int (x_l)^4\right)^{1/2}\right\}\left\{\left(\int |h|^2\right)^{1/2}+\left(\int|u^{(m)}|^2\right)^{1/2}\right\}+\int(x_l)^2|u^{(m)}|^2.
\end{split}
\label{hisa}
\end{equation}
Due to \eqref{blo10} and \eqref{blo14}, we can estimate the right-hand side of 
\eqref{hisa} by $c(\mu^{(m)})^2 ((s_j-s)^{1/2}+\mu^{(m)})$. On the other hand one can
estimate the left-hand side of \eqref{hisa} (devided by $(\mu^{(m)})^2$) just like \eqref{tai4} using 
\eqref{weak1}-\eqref{weak2}. It is important to note that we have the strong
convergence in $L^2(B_{1/3}^k)$ at $t=s_j$. This shows \eqref{lsc2}. 
Since $\tilde{f}_l$ is smooth and $\{s_j\}_{j=1}^{\infty}$ is dense, \eqref{lsc1}
and \eqref{lsc2} show that for a.e$.$ $s\in (-1/4,1/4)$, 
$\tilde{f}_l^{(m)}(\cdot,s)$ converges strongly in $L^2(B_{1/3}^k)$
to $\tilde{f}_l(\cdot,s)$ as $m\rightarrow\infty$. Since $\|\tilde{f}_l^{(m)}(\cdot,s)
-\tilde{f}_l(\cdot,s)\|_{L^2(B_{1/3}^k)}$ is bounded uniformly in $s\in
(-1/4,1/4)$, the dominated convergence theorem shows that $\|\tilde{f}_l^{(m)}
-\tilde{f}_l\|_{L^2(\Omega')}$ converges to 0. This ends the proof of Lemma \ref{heatsol}. 
\hfill{$\Box$}

Next we define $T^{(m)}\in {\bf G}(n,k)$ as the image of the map $x\longmapsto 
\mu^{(m)}x\cdot\nabla\tilde{f}(0,0)$, which is the tangent space of the
graph $(x,\mu^{(m)}\tilde{f}(x,0))$ at $x=0$. Also define $A^{(m)}\in {\bf A}(n,k)$ by
$A^{(m)}
=T^{(m)}+(0,\mu^{(m)}\tilde{f}(0,0))$.
\begin{lemma}
There exists a constant $c_{12}$ which depends only on $c_{11}$ and $\Lambda_*$
with the following property. 
For $0<\theta\leq (16\Lambda_*)^{-1/2}$ we have $\|T-T^{(m)}\|\leq c_{12}\mu^{(m)}$ and 
\begin{equation}
\limsup_{m\rightarrow\infty}(\mu^{(m)})^{-2}\int_{C(T,2\theta)\times
(-\theta^2 \Lambda_* ,\theta^2 \Lambda_*)} {\rm dist}\, (A^{(m)},x)^2\, d\|V_t^{(m)}\|dt
\leq c_{12}\theta^{k+6}.
\label{dec1}
\end{equation}
\label{decestfin}
\end{lemma}
{\it Proof of Lemma \ref{decestfin}}. We first note that any interior partial derivatives of $\tilde{f}$ 
may be estimated depending only on the order of differentiations and $c_{11}$
due to the fact that $\tilde{f}$ is a solution of the heat equation and by the
standard linear regularity theory. In particular, by the second order
Taylor expansion at the origin, we have (with $|f(x,t)-f(x,0)|\leq c|t|\leq c\theta^2$)
\begin{equation}
\int_{B_{2\theta}^k\times(-\theta^2\Lambda_* ,\theta^2\Lambda_*)}|\tilde{f}(x,t)
-\tilde{f}(0,0)-x\cdot\nabla\tilde{f}(0,0)|^2\, d{\mathcal H}^{k+1}\leq
c_{13}\theta^{k+6}
\label{heatest}
\end{equation}
for any $0<\theta\leq (16\Lambda_*)^{-1/2}$ and for some $c_{13}$ depending only on $c_{11}$
and $\Lambda_*$. 
Next, by the Lipschitz approximation of $V_t^{(m)}$ and ${\rm dist}\, (A^{(m)},x)\leq c(c_{11})\mu^{(m)}$ on the support of $\|V_t^{(m)}\|$, we may prove
(just like \eqref{tai4}) that
\begin{equation}
\begin{split}
& (\mu^{(m)})^{-2}\int_{C(T,2\theta)\times(-\theta^2\Lambda_* ,\theta^2\Lambda_* )}
{\rm dist}\,(A^{(m)},x)^2\, d\|V_t^{(m)}\| dt \\
&=o(1)+(\mu^{(m)})^{-2}\int_{B^k_{2\theta}
\times(-\theta^2\Lambda_* ,\theta^2\Lambda_* )}|(T^{(m)})^{\perp}((x,f^{(m)}(x,t))
-(0,\tilde{f}(0,0)\mu^{(m)}))|^2\, d{\mathcal H}^{k+1}.
\end{split}
\label{heatest2}
\end{equation}
Since $(x,x\cdot\nabla\tilde{f}(0,0)\mu^{(m)})\in T^{(m)}$,
\begin{equation}
\begin{split}
&|(T^{(m)})^{\perp}((x,f^{(m)}(x,t))-(0,\tilde{f}(0,0)\mu^{(m)}))|\\
&=|(T^{(m)})^{\perp}((x,f^{(m)}(x,t))-(0,\tilde{f}(0,0)\mu^{(m)})
-(x,x\cdot\nabla\tilde{f}(0,0)\mu^{(m)}))| \\
&\leq |f^{(m)}(x,t)-\tilde{f}(x,t)\mu^{(m)}|+|\tilde{f}(x,t)-\tilde{f}(0,0)
-x\cdot\nabla\tilde{f}(0,0)|\mu^{(m)}
\end{split}
\label{heatest3}
\end{equation}
where we added and subtracted $(x,\tilde{f}(x,t)\mu^{(m)})$ and used the triangle 
inequality in the second line. Now we substitute \eqref{heatest3} into \eqref{heatest2}
and use \eqref{heatest}. Since we have the strong $L^2(\Omega')$ convergence
(Proposition \ref{heatsol}) we obtain \eqref{dec1} for a suitable $c_{12}$.
\hfill{$\Box$}

Fix $\theta_*$ so that 
\begin{equation}
0<\theta_*\leq (16\Lambda_*)^{-1/2},\hspace{.5cm}3c_{12}\theta_*^2\leq \theta_*^{2\varsigma}.
\label{theta}
\end{equation}
For $1<\gamma<2$ we have \eqref{dec1}
for $\theta=\gamma \theta_*$ and hence for all large $m$, 
\begin{equation}
\begin{split}
\int_{C(T,2\gamma\theta_*)\times (-\theta_*^2\Lambda_*,\theta_*^2 \Lambda_*)}&
{\rm dist}\, (A^{(m)},x)^2 \, d\|V_t^{(m)}\|dt\leq (c_{12}(\gamma\theta_*)^{k+6}+o(1))
(\mu^{(m)})^2 \\ &<\frac12 \gamma^{k+6}\theta_*^{2\varsigma +k+4}(\mu^{(m)})^2
\end{split}
\label{dame2}
\end{equation}
by \eqref{theta}, with $\|T^{(m)}-T\|\leq c_{12}\mu^{(m)}$ and $A^{(m)}$ is 
parallel to $T^{(m)}$. 
For all sufficiently large $m$, due to Lemma \ref{heiblo}, $C(W^{(m)},2\theta_*)\cap{\rm spt}\|V_t^{(m)}\|\subset
C(T,2\gamma\theta_*)$. Thus we have from \eqref{dame2}
\begin{equation}
\begin{split}
\limsup_{m\rightarrow\infty}\theta_*^{-(k+4)}(\mu^{(m)})^{-2} & \int_{C(W^{(m)},2\theta_*)\times(-\theta_*^2\Lambda_*,
\theta_*^2\Lambda_*)} {\rm dist}\, (A^{(m)},x)^2\, d\|V_t^{(m)}\|dt \\
& \leq \frac12 \gamma^{k+6}\theta_*^{2\varsigma}.
\end{split}
\label{dame4}
\end{equation}
By letting $\gamma\rightarrow 1$, \eqref{dame4} contradicts with \eqref{blo1no}, and the first part of the
Proposition \ref{blowup1} is proved.  

For \eqref{blo2ag} and \eqref{blo3ag}, just as we obtained \eqref{blo15} by Corollary \ref{popscor}, we
may obtain 
\begin{equation}
\int_{-1/2}^{1/2}|\theta_*^{-k}\|V_t\|(\phi_{T,\theta_*}^2)-{\bf c}|\, dt\leq \tilde{K}
(\mu^2+C(u)).
\label{mous}
\end{equation}
Since $\|W-T\|\leq \varepsilon_4$, $\|u\|<\varepsilon_4$ and $\mu<\varepsilon_4$, by further restricting $\varepsilon_4$ 
so that the right-hand side of \eqref{mous} is sufficiently small, and choosing
generic $\tilde{t}_1$ and $\tilde{t}_2$ satisfying the same inequality as \eqref{mous} (with $W$ and a different constant), we may conclude that such $t$ can be chosen. 
This concludes the proof of Proposition \ref{blowup1}
\hfill{$\Box$}
\begin{cor}(cf. \cite[8.17]{Allard})
Corresponding to $E_1$, $\nu$, $p$, $q$ there exist $0<\varepsilon_5<1$ and $1<c_{15}<\infty$ with the following 
property. Under the assumptions of Proposition \ref{blowup1} where $\varepsilon_4$ is
replaced by $\varepsilon_5$, and with $R^{\varsigma}\|u\|<\varepsilon_5$, 
\begin{itemize}
\item[(1)] there exists a unique element $a\in {\rm spt}\, \|V_0\|\cap\{x\,:\, W(x)=0\}$,
\item[(2)] there exists $T_{\infty}\in {\bf G}(n,k)$ such that 
\begin{equation}{\rm Tan}\, ({\rm spt}\, \|V_0\|,a)
\subset T_{\infty}\label{taninc1}
\end{equation}
and
\begin{equation}
\|T_{\infty}-T\|\leq c_{15}\max\{\mu,c_{14} R^{\varsigma}\|u\|\},
\label{taninc2}
\end{equation}
\item[(3)] whenever $0<s<R$, there are $T_s\in {\bf G}(n,k)$ and $A_s\in {\bf A}(n,k)$ such
that $A_s$ is parallel to $T_s$ and
\begin{equation}
\begin{split}
\|T_s-T_{\infty}\|+&\left(s^{-k-4}\int_{-s^2\Lambda_*}^{s^2\Lambda_*}\int_{C(W,2s)}
{\rm dist}\, (x,A_s)^2\, d\|V_t\|dt\right)^{\frac12} \\
&\leq c_{15}(s/R)^{\varsigma}\max\{\mu,c_{14}R^{\varsigma}\|u\|\}
.
\end{split}
\label{taninc3}
\end{equation}
\end{itemize}
\label{decor}
\end{cor}
{\it Proof}. We may assume $R=1$ after a change of variables. 
We choose $0<\varepsilon_5<1$ and $1<c_{15}<\infty$ so that
\begin{equation}
c_{14}\varepsilon_5<\varepsilon_4,
\label{ve1}
\end{equation}
\begin{equation}
\varepsilon_5+\sum_{j=1}^{\infty}(c_{14})^2 \theta_*^{(j-1)\varsigma}\varepsilon_5<\varepsilon_4,
\label{ve2}
\end{equation}
\begin{equation}
2c_{14}\sum_{j=0}^{\infty}\theta_*^{(j-1)\varsigma}\leq c_{15},
\label{ve3}
\end{equation}
\begin{equation}
2\theta_*^{-(k+4)/2-\varsigma}\leq c_{15}.
\label{ve4}
\end{equation}
We inductively prove the following. Set $T_0=T$ and $A_0=A$. Suppose for $j=1,\cdots,l$ that
there are $T_j\in {\bf G}(n,k)$ and $A_j\in {\bf A}(n,k)$ such that $A_j$ is parallel to $T_j$ and
\begin{equation}
\|T_j-T_{j-1}\|\leq c_{14}\theta_*^{(j-1)\varsigma}\max\{\mu,c_{14}\|u\|\},
\label{ve5}
\end{equation}
\begin{equation}
\mu_j:=\left(\theta_*^{-j(k+4)}\int_{-\theta_*^{2j}\Lambda_*}^{\theta_*^{2j}\Lambda_*}
\int_{C(W,2\theta_*^j)}{\rm dist}\, (x,A_j)^2\, d\|V_t\|dt\right)^{\frac12}\leq \theta_*^{j\varsigma}\max\{\mu,
c_{14}\|u\|\},
\label{ve6}
\end{equation}
\begin{equation}
(-\Lambda_*+1)\theta_*^{2j}\leq \exists t_{1,j}\leq (-\Lambda_*+2)\theta_*^{2j}\,\,:\,\,
\theta_*^{-kj}\|V_{t_{1,j}}\|(\phi_{W,\theta_*^j}^2)<(2-\nu){\bf c},
\label{ve7}
\end{equation}
\begin{equation}
(\Lambda_*-2)\theta_*^{2j}\leq \exists t_{2,j}\leq (\Lambda_*-1)\theta_*^{2j}\,\,:\,\,
\theta_*^{-kj}\|V_{t_{2,j}}\|(\phi_{W,\theta_*^j}^2)>\nu{\bf c}.
\label{ve8}
\end{equation}
Since $\varepsilon_5<\varepsilon_4$, Proposition \ref{blowup1} provides the proof for
the validity of $j=1$ case. Under the inductive assumption up to $j=l$, we have
\begin{equation}
\mu_l\leq \theta_*^{l\varsigma}\max\{\mu,c_{14}\|u\|\}<\varepsilon_4
\label{ve9}
\end{equation}
by \eqref{blo1}, $\|u\|<\varepsilon_5$ and \eqref{ve1},
\begin{equation}
\begin{split}
\|T_l-W\| & \leq \|T_0-W\|+\sum_{j=1}^l \|T_j-T_{j-1}\| \leq \|T-W\| \\ & + \sum_{j=1}^{l}c_{14}
\theta_*^{(j-1)\varsigma} \max\{\mu,c_{14}\|u\|\} 
\leq \varepsilon_5+\sum_{j=1}^l (c_{14})^2\varepsilon_5 \theta_*^{(j-1)\varsigma}<\varepsilon_4
\end{split}
\label{ve10}
\end{equation}
by \eqref{blo0}, \eqref{blo1}, $\|u\|<\varepsilon_5$, \eqref{ve5} and \eqref{ve2}. 
For $R=\theta_*^l$, the assumptions \eqref{blo0}-\eqref{blo3} are now satisfied due to
\eqref{ve9}, \eqref{ve10}, \eqref{ve7}, \eqref{ve8} as well as $\theta_*^{l\varsigma}\|u\|<\varepsilon_4$.
Thus we have $T_{l+1}\in{\bf G}(n,k)$ and $A_{l+1}\in {\bf A}(n,k)$ such that 
$A_{l+1}$ is parallel to $T_{l+1}$,
\begin{equation}
\|T_{l+1}-T_l\|\leq c_{14}\mu_l\leq c_{14}\theta_*^{l\varsigma}\max\{\mu,c_{14}\|u\|\}
\label{ve11}
\end{equation}
by \eqref{imblo2} and \eqref{ve6},
\begin{equation}
\begin{split}
\mu_{l+1}:&=\left(\theta_*^{-(l+1)(k+4)}\int_{-\theta_*^{2(l+1)}\Lambda_*}^{\theta_*^{2(l+1)}\Lambda_*}
\int_{C(W,2\theta_*^{l+1})}{\rm dist}\, (x,A_{l+1})^2\, d\|V_t\|dt\right)^{\frac12}\\
& \leq \theta_*^{\varsigma}\max\{\mu_l,c_{14}\theta_*^{l\varsigma}\|u\|\}
\leq \theta_*^{(l+1)\varsigma}\max\{\mu,c_{14}\|u\|\}
\end{split}
\end{equation}
by \eqref{imblo3} and \eqref{ve6}. \eqref{ve7} and \eqref{ve8} for $j=l+1$ are also satisfied due to 
\eqref{blo2ag} and \eqref{blo3ag}. Hence the next inductive assumptions \eqref{ve5}-\eqref{ve8}
are satisfied. It is now clear from \eqref{ve5} that there exists $T_{\infty}=\lim_{j\rightarrow\infty}
T_j\in {\bf G}(n,k)$ satisfying \eqref{taninc2} due to \eqref{ve5} and \eqref{ve3}. 
To prove (1), first assume that ${\rm spt}\, \|V_0\|\cap \{x\, :\, W(x)=0\}=\emptyset$. 
By definition, there exists some $\gamma>0$ such that $\|V_0\|C(W,\gamma)=0$. 
By the unit density assumption (A1) and Proposition \ref{cam1}, we may then prove
that there exists $0<\gamma'<\gamma$ such that $\|V_t\|(C(W,\gamma'))=0$ for $0\leq 
t\leq \gamma'$. But this contradicts with \eqref{ve8} for all sufficiently large $j$. 
Thus we prove ${\rm spt}\, \|V_0\|\cap \{x\, :\, W(x)=0\}\neq \emptyset$.
To prove the uniqueness, we observe that
\begin{equation}
{\rm spt}\, \|V_t\|\cap \{x\, :\, {\rm dist}\, (x,A_j)>o(1)\theta_*^j\}\cap C(W,3\theta_*^j/2)=\emptyset,
\,\,\, (-\Lambda_*+1)\theta_*^{2j}\leq \forall t\leq \Lambda_*\theta_*^{2j}
\label{ve12}
\end{equation}
where $o(1)$ here means $0<o(1)\rightarrow 0$ as $j\rightarrow\infty$. 
The proof is almost identical to that of Lemma \ref{heiblo} after
a change of variables and with some notational modifications. Let $a_j:=
A_j\cap \{x\,:\, W(x)=0\}$. One can prove $\lim_{j\rightarrow\infty} a_j=a:=
{\rm spt}\, \|V_0\|\cap \{x\, :\, W(x)=0\}$, by means of \eqref{ve12}, which proves (1). Also 
\eqref{ve12} shows that ${\rm Tan}\, ({\rm spt}\,\|V_0\|, a)\subset T_{\infty}$
(recall \eqref{deftan}), 
proving (2). Finally for $0<s<1$, choose $j$ such that $\theta_*^{j+1}\leq s<\theta_*^j$,
and let $T_s=T_j$, $A_s=A_j$. Then we have
\begin{equation}
\begin{split}
\|T_s-T_{\infty}\|&\leq c_{14}\max\{\mu,c_{14}\|u\|\}\sum_{l=0}^{\infty} \theta_*^{(j+l)\varsigma} \\
&=c_{14}\sum_{l=0}^{\infty}\theta_*^{(l-1)\varsigma} \max\{\mu,c_{14}\|u\|\}\theta_*^{(j+1)\varsigma}
\leq c_{15}\max\{\mu,c_{14}\|u\|\}s^{\varsigma}/2 
\end{split}
\label{ve13}
\end{equation}
by \eqref{ve3} and 
\begin{equation}
\begin{split}
&\left(s^{-(k+4)}\int_{-s^2\Lambda_*}^{s^2\Lambda_*}
\int_{C(W,2s)}{\rm dist}\, (x,A_s)^2\, d\|V_t\|dt\right)^{\frac12} \leq
\theta_*^{-(k+4)/2}\mu_j \\
&\leq \theta_*^{-(k+4)/2}\max\{\mu,c_{14}\|u\|\} \theta_*^{j\varsigma} \leq \theta_*^{-(k+4)/2-\varsigma}\max\{\mu,c_{14}\|u\|\}s^{\varsigma} \\
&\leq c_{15}\max\{\mu,c_{14}\|u\|\} s^{\varsigma}/2
\end{split}
\label{ve14}
\end{equation}
by \eqref{ve4}. Combining \eqref{ve13} and \eqref{ve14}, we obtain \eqref{taninc3}.
\hfill{$\Box$}
\begin{thm}
Corresponding to $1\leq E_1<\infty$, $0<\nu<1$, $p$ and $q$ with \eqref{expcond}, there
exist $0<\varepsilon_6<1$, $0<\sigma\leq 1/2$, $2<\Lambda_3<\infty$ 
and $1<c_{16}<\infty$ with the following
property. For $T\in {\bf G}(n,k)$, $0<R<\infty$, $U=C(T,3R)$ and $(0,\Lambda)$
replaced by $(-\Lambda_3 R^2,\Lambda_3R^2)$, suppose $\{V_t\}_{-\Lambda_3 R^2\leq t\leq
\Lambda_3 R^2}$ and $\{u(\cdot, t)\}_{-\Lambda_3R^2\leq t\leq \Lambda_3R^2}$ satisfy (A1)-(A4).
Suppose
\begin{equation}
\mu:=\left(R^{-(k+4)}\int_{-\Lambda_3R^2}^{\Lambda_3R^2}\int_{C(T,3R)}
|T^{\perp}(x)|^2\, d\|V_t\|dt\right)^{\frac12}<\varepsilon_6,
\label{regth1}
\end{equation}
\begin{equation}
R^{\varsigma}\|u\|:=R^{\varsigma}\|u\|_{L^{p,q}(C(T,3R)\times(-\Lambda_3R^2,\Lambda_3R^2))}
<\varepsilon_6,
\label{regth2}
\end{equation}
\begin{equation}
(-\Lambda_3+3/2)R^2\leq \exists t_1\leq (-\Lambda_3+2)R^2\,\, \, :\,\, \, R^{-k}\|V_{t_1}\|
(\phi_{T,R}^2)<(2-\nu){\bf c},
\label{regth3}
\end{equation}
\begin{equation}
(\Lambda_3-2)R^2\leq \exists t_2\leq (\Lambda_3-3/2)R^2\,\,\,:\,\,\, R^{-k}\|V_{t_2}\|
(\phi_{T,R}^2)>\nu{\bf c}.
\label{regth4}
\end{equation}
Denote $\tilde{D}:=(T\cap B_{\sigma R})\times (-R^2/4,R^2/4)$. Then there are $f\,:\,\tilde{D}\rightarrow
T^{\perp}$ and $F\,:\, \tilde{D}\rightarrow {\mathbb R}^n$ such that
$T(F(y,t))=y$ and $T^{\perp}(F(y,t))=f(y,t)$ for all $(y,t)\in \tilde{D}$,
\begin{equation}
{\rm spt}\, \|V_t\|\cap C(T,\sigma R)={\rm image}\, F(\cdot, t)\hspace{.3cm}\forall t\in (-R^2/4,R^2/4),
\label{regth5}
\end{equation}
\begin{equation}
\mbox{$f(x,t)$ is differentiable with respect to $x$ at every point of $\tilde{D}$,}
\label{regth6}
\end{equation}
\begin{equation}
R^{-1}|f(y,s)|+\|\nabla f(y,s)\|\leq c_{16}\max\{\mu,R^{\varsigma}\|u\|\},\hspace{.3cm}\forall (y,s)\in \tilde{D},
\label{regth6.5}
\end{equation}
\begin{equation}
\begin{split}
\|\nabla f(y_1,s_1)-\nabla f(y_2,s_2)\| \leq c_{16}\max\{\mu,R^{\varsigma}\|u\|\}&\left(R^{-1}\max \{|y_1-y_2|,
|s_1-s_2|^{1/2}\}\right)^{\varsigma},\\
\forall (y_1,s_1),\,\forall(y_2,s_2)\in \tilde{D},
\end{split}
\label{regth7}
\end{equation}
\begin{equation}
\begin{split}
|f(y,s_1)-f(y,s_2)|\leq c_{16} \max\{\mu,R^{\varsigma}\|u\|\}&(R^{-1}|s_1-s_2|^{1/2})^{1+\varsigma}, \\
\forall (y,s_1),\,\forall (y,s_2)\in \tilde{D}.
\end{split}
\label{regth7.5}
\end{equation}
\label{regth}
\end{thm}
{\it Proof}. Let $\Lambda_*$, $\theta_*$, $c_{15}$, $\varepsilon_5$ be fixed by Proposition \ref{blowup1} and
Corollary \ref{decor} corresponding to $E_1$ and $\nu$ there replaced by $\nu/2$. We set $\Lambda_3=\Lambda_*+1/4$. 
Without loss of generality we may 
assume $R=1$. Set $\sigma:= \min\{1/2,{\bf c}\nu(2^{k+1}  E_1\sup |\nabla(\phi_{T}^2)|)^{-1}\}$.
Then $\sigma\leq 1/2$ and we have for any $y\in B^k_{\sigma}$
\begin{equation}
\begin{split}
\int_{C(T,2)}\phi_T^2(\cdot-y)\,d\|V_{t_1}\|& \leq \int_{C(T,2)}\phi_T^2\, d\|V_{t_1}\|
+|y| \|V_{t_1}\|(C(T,3/2))\sup\, |\nabla(\phi_T^2)|\\
&\leq (2-\nu){\bf c}+|y|\|V_{t_1}\|(B_2)\sup\, |\nabla(\phi_T^2)|\leq (2-\nu/2){\bf c}
\end{split}
\label{regth8}
\end{equation}
where we used \eqref{regth3}, $C(T,3/2)\cap {\rm spt}\, \|V_{t_1}\|\subset B_2$ for sufficiently small 
$\varepsilon_6$ which follows from Proposition \ref{el2elinf}, (A2) and the definition of $\sigma$. 
The similar bound for $t_2$ may be obtained. Note that for any $s\in (-1/4,1/4)$, $t_1\in (-\Lambda_*+s+1,
-\Lambda_*+s+2)$ and $t_2\in (\Lambda_*+s-2,\Lambda_*+s-1)$ due to the intervals in \eqref{regth3}
and \eqref{regth4}. Thus all the assumptions for Corollary \ref{decor} are satisfied for domain 
centered at $(y,s)$, i.e., for $C(T,y,2)\times (-\Lambda_*+s,\Lambda_*+s)$. We may conclude that
there exists a unique element which we define to be $F(y,s)\in {\rm spt}\, \|V_s\|\cap \{x\,:\,
T(x)=y\}$ for each $(y,s)\in \tilde{D}$. Define $f(y,s):=T^{\perp}(F(y,s))$. It is also clear that
$F(y,s)$ is differentiable with respect to space variables and that ${\rm image}\, \nabla F(y,s)
=T_{\infty}(y,s)$, the latter being $T_{\infty}$ for $(y,s)$ in Corollary \ref{decor}. Then \eqref{regth6.5}
follows from \eqref{taninc2} and Proposition \ref{el2elinf} for a suitable $c_{16}$. 
For $(y_1,s_1),\,(y_2,s_2)\in \tilde{D}$, set $d:=\max\{|y_1-y_2|,|s_1-s_2|^{1/2}\}$ and $w:=2d$. 
We have $T_w(y_1,s_1)\in {\bf G}(n,k)$
and $A_w(y_1,s_1)\in {\bf A}(n,k)$ which is parallel to $T_w(y_1,s_1)$ satisfying 
\begin{equation}
\begin{split}
\|T_w(y_1,s_1)-T_{\infty}(y_1,s_1)\|&+\left(w^{-(k+4)}\int_{s_1-w^2\Lambda_*}^{s_1+w^2\Lambda_*}
\int_{C(T,y_1,2w)}{\rm dist}\, (x,A_w(y_1,s_1))^2\, d\|V_t\|dt\right)^{\frac12}\\
&\leq c_{15}w^{\varsigma}
\max\{\mu, c_{14}\|u\|\}.
\end{split}
\label{regth9}
\end{equation}
By the definition of $w$ and $d$, one can check that
\begin{equation}
C(T,y_2,2d)\times(s_2-d^2\Lambda_*, s_2+d^2\Lambda_*)\subset C(T, y_1, 2w)
\times(s_1-w^2\Lambda_*,s_1+w^2\Lambda_*).
\label{regth10}
\end{equation}
By \eqref{regth10},
\begin{equation}
\begin{split}
&\left(d^{-(k+4)}\int_{s_2-d^2\Lambda_*}^{s_2+d^2\Lambda_*}\int_{C(T,y_2,2d)}
{\rm dist}\, (x,A_w(y_1,s_1))^2\, d\|V_t\|dt\right)^{\frac12}\\
&\leq 2^{(k+4)/2} \left(w^{-(k+4)}\int_{s_1-w^2\Lambda_*}^{s_1+w^2\Lambda_*}
\int_{C(T,y_1,2w)}{\rm dist}\, (x,A_w(y_1,s_1))^2\, d\|V_t\|dt\right)^{\frac12}.
\end{split}
\label{regth11}
\end{equation}
Since we may choose $\varepsilon_6$ small so that $\|T_{w}(y_1,s_1)-T\|\leq  \varepsilon_5$, and
conditions  \eqref{blo2} and \eqref{blo3} corresponding to the domain $C(T,y_2,2d)\times
(s_2-d^2\Lambda_*, s_2+d^2\Lambda_*)$ are satisfied, we may apply Corollary \ref{decor}
to this domain by restricting $\varepsilon_6$ small. 
Thus we obtain with \eqref{regth9} and \eqref{regth10} that 
\begin{equation}
\|T_{w}(y_1,s_1)-T_{\infty}(y_2,s_2)\|\leq c_{15}\max\{2^{(k+4)/2}c_{15}w^{\varsigma}
\max\{\mu,c_{14}\|u\|\},c_{14}d^{\varsigma}\|u\|\}.
\label{regth12}
\end{equation}
\eqref{regth9} and \eqref{regth12} shows with an appropriate $c_{16}$ that 
\begin{equation}
\|T_{\infty}(y_1,s_1)-T_{\infty}(y_2,s_2)\|\leq c_{16} d^{\varsigma}\max\{\mu,\|u\|\}.
\label{regth13}
\end{equation}
Since ${\rm image}\, \nabla F(y,s)=T_{\infty}(y,s)$, \eqref{regth13} shows \eqref{regth7}
via \cite[8.9(5)]{Allard} with a new suitable $c_{16}$. 
For \eqref{regth7.5}, we observe that \eqref{regth9} combined with Proposition \ref{el2elinf}
gives the estimate with a suitable $c_{16}$. 
\hfill{$\Box$}
\begin{rem}
Note that Theorem \ref{regth} shows that ${\rm spt}\, \|V_t\|\cap C(T,\sigma R)$ for 
$t\in (-R^2/4,R^2/4)$ is a set of 
$C^{1,\varsigma}$ regular points. 
\end{rem}
\label{holder}
\section{Partial regularity}
In this section, we prove the main partial regularity result. The main idea of proof
comes from \cite[6.12]{Brakke}, though we greatly simplify the overall computations.
First we give a more convenient form of regularity criterion which is suitable for 
establishing partial regularity. 
\begin{prop}
Corresponding to $1\leq E_1<\infty$, $0<\nu<1$, $p,\, q$ with \eqref{expcond} there exists $2<L<\infty$, $3<\Lambda_4<
\infty$, $0<\varepsilon_7<1$ with the following property. 
Assume $\{V_t\}_{0<t<\Lambda}$
and $\{u(\cdot,t)\}_{0<t<\Lambda}$ satisfy (A1)-(A4). For $(a,s)\in U\times(0,\Lambda)$, assume that 
for some $R>0$, $B_{RL}(a)\times(s-R^2\Lambda_4,
s+R^2\Lambda_4)\subset U\times (0,\Lambda)$ and for some $T\in {\bf G}(n,k)$,
\begin{equation}
\left(R^{-(k+2)}\int_{B_{RL}(a)}|T^{\perp}(x-a)|^2\, d\|V_{s-R^2\Lambda_4}\|\right)^{\frac12}<\varepsilon_7,
\label{pa1}
\end{equation}
\begin{equation}
R^{\varsigma}\|u\|_{L^{p,q}(B_{RL}(a)\times(s-R^2\Lambda_4,s+R^2\Lambda_4))}<\varepsilon_7,
\label{pa2}
\end{equation}
\begin{equation}
s-R^2(\Lambda_4-5/2)\leq \exists t_1\leq s-R^2(\Lambda_4-3)\,\,:\,\, R^{-k}\|V_{t_1}\|(\phi_{R}^2(\cdot-a))
<(2-\nu){\bf c},
\label{pa3}
\end{equation}
\begin{equation}
s+R^2(\Lambda_4-3)\leq \exists t_2\leq s+R^2 (\Lambda_4-5/2)\,\, :\,\, R^{-k}\|V_{t_2}\|(\phi_{R}^2(\cdot-a))
>\nu{\bf c}.
\label{pa4}
\end{equation}
Here, $\phi_{R}(x):=\phi(|x|/R)$ with $\phi$ defined in \eqref{pops1} and ${\bf c}$ defined in \eqref{pops1.5}.
Then $(a,s)$ is a $C^{1,\varsigma}$ regular point.
\label{prprop}
\end{prop}
\label{partialreg}
{\it Proof}. We check that Theorem \ref{regth} is applicable for the conclusion. We may set $R=1$, 
$a=0$ and $s=0$ after a change of variables. Let $\Lambda_3$ and 
$\varepsilon_6$ be constants corresponding to $E_1$ and $\nu$ there replaced by $\nu/2$. 
Set $\Lambda_4=\Lambda_3+1$. We use 
Proposition \ref{uniel2} with $\Lambda=2\Lambda_4$ and obtain $c_{17}$. 
Fix a large $L>2$ so that we have
\begin{equation}
c_{17}L^{k+2}\exp(-(L/4 -1)^2/(8\Lambda_4))<\varepsilon_6^2/(8\Lambda_4).
\label{pa4.5}
\end{equation}
Then restrict $\varepsilon_7\leq \varepsilon_6$ so that 
\begin{equation}
4^{k+2}c_{17}(4^{2\varsigma}\varepsilon_7^2+4^{\varsigma}\varepsilon_7)(L/4)^2<\varepsilon_6^2/(8\Lambda_4).
\label{pa4.7}
\end{equation}
By \eqref{uniel2a} with $R=4$ and $L$ there replaced by $L/4$ with $L$ given in the assumption, and 
\eqref{pa4.5}, \eqref{pa4.7}, \eqref{pa2}, we obtain
\begin{equation}
\int_{B_4}|T^{\perp}(x)|^2\, d\|V_t\|\leq \exp(1/(8\Lambda_4))\int_{B_{L}}|T^{\perp}(x)|^2\, d\|V_{-
\Lambda_4}\|+\varepsilon_6^2/(4\Lambda_4)
\label{pa5}
\end{equation}
for all $t\in [-\Lambda_4,\Lambda_4]$. By integraring \eqref{pa5} and
by \eqref{pa1},
we have
\begin{equation}
\int_{-\Lambda_4}^{\Lambda_4}\int_{B_4}|T^{\perp}(x)|^2\, d\|V_t\|dt\leq 2\Lambda_4
\exp(1/(8\Lambda_4))\varepsilon_7^2+\varepsilon_6^2/2.
\label{pa6}
\end{equation}
By \eqref{pa6} and Proposition 
\ref{el2elinf}, by further restricting $\varepsilon_7$ if necessary, we have
\begin{equation}
{\rm spt}\, \|V_t\|\cap B_{16/5}\subset \{|T^{\perp}(x)|\leq 1/5\}
\label{pa7}
\end{equation}
for $t\in [-\Lambda_4+1,\Lambda_4-1]=[-\Lambda_3,\Lambda_3]$.
By re-defining $V_t=0$ on $C(T,3)\setminus B_{16/5}$ for $t\in [-\Lambda_3,\Lambda_3]$, 
one can check that the newly defined $\{V_t\}$ satisfies \eqref{maineq} on $C(T,3)\times
[-\Lambda_3,\Lambda_3]$. By restricting $\varepsilon_7$ further, we may guarantee that 
\eqref{regth1} is satisfied by \eqref{pa6}. For conditions \eqref{regth3} and \eqref{regth4},
note that $\int |\phi^2_T-\phi^2_1|\, d\|V_t\|$ can be made arbitrarily small by Proposition 
\ref{el2elinf} and choosing small $\varepsilon_7$, i.e., by making ${\rm dist}\, (T,{\rm spt}\,
\|V_t\|)$ sufficiently small. Thus, having assumed \eqref{pa3} and \eqref{pa4}, we can 
guarantee that \eqref{regth3} and \eqref{regth4} hold for $\nu/2$ instead of $\nu$. 
We now have all the assumptions for Theorem
\ref{regth}, and hence $(a,s)$ is a $C^{1,\varsigma}$ regular point. 
\hfill{$\Box$}
\begin{lemma}
Under the assmuptions (A1)-(A4), there exists a co-countable set $G\subset(0,\Lambda)$
such that $\|V_t\|(\phi)$ is continuous at $t\in G$ for any $\phi\in C^2_c(U;{\mathbb R}^+)$.
\label{contime}
\end{lemma}
{\it Proof}. Fix $\phi\in C^2_c(U;{\mathbb R}^+)$. Computing as in \eqref{hhhp}, we have
for any $0<t_2<t_2<\Lambda$
\begin{equation}
\begin{split}
\|V_{t_2}\|(\phi)-\|V_{t_1}\|(\phi)&\leq \int_{t_1}^{t_2}\left(\int_U\frac{|\nabla\phi|^2}{\phi}
+|u|^2\phi+|\nabla\phi||u|\, d\|V_t\|\right)dt\\
&=:\int_{t_1}^{t_2}\Phi(t)dt,
\end{split}
\label{pa8}
\end{equation}
where $\Phi$ is integrable due to (A2) and (A4). \eqref{pa8} shows $\|V_t\|(\phi)-\int_0^t\Phi(s)\, ds$ is monotone
decreasing, and there are at most countably many discontinuities for such function. 
Since $\int_0^t\Phi(s)\, ds$ is continuous, such discontinuities are caused only by $\|V_t\|(\phi)$.
By choosing a dense and countable set of functions in $C^2_c(U;{\mathbb R}^+)$, we may conclude the proof.
\hfill{$\Box$}

We focus our attention at time $t\in (0,\Lambda)$ when $V_t$ is unit density and $t\in G$. By (A1) and Lemma \ref{contime},
such time is a.e$.$ on $(0,\Lambda)$. Without loss of generality we assume $t=0\in (-t, \Lambda-t)$ in the following, and 
assume $V_0=|M|$. Without loss of generality, we may assume $M\subset {\rm spt}\, \|V_0\|$.
By the standard measure-theoretic argument (\cite[2.10.19]{Federer}),
for ${\mathcal H}^k$ a.e$.$ $x\in {\rm spt}\, \|V_0\|\setminus M$, we have $\Theta^k(\|V_0\|,x)=0$.
For ${\mathcal H}^k$ a.e$.$ $x\in M$, $\Theta^k(\|V_0\|,x)=1$ and a unique approximate tangent
space ${\rm Tan}_x M$ exists. 
\begin{define}
Set 
\begin{equation}
A_1:={\rm spt}\, \|V_0\|\setminus M,\,\,\, A_2:=M\setminus \{\mbox{$C^{1,\varsigma}$ regular points}\}
\label{cot1}
\end{equation}
\end{define}
We aim to prove that ${\mathcal H}^k(A_1)={\mathcal H}^k( A_2)=0$. By definition, $A_1\cup A_2$ is closed.
\begin{lemma} We have
${\mathcal H}^k(A_1)=0$.
\label{a1is0}
\end{lemma}
{\it Proof}. It suffices to prove ${\mathcal H}^k(A_1\cap B_r(a))=0$ for arbitrary $B_{3r}(a)\subset\subset U$.
After a change of variables, let $a=0$.
Assume for a contradiction that ${\mathcal H}^k(A_1\cap
B_r)>0$. 
Let $c_{18}$ and $c_{19}$ be constants given in Corollary \ref{clearing}.
Fix a sufficiently small $0<R_0<r/2$ so that 
\begin{equation}
R_0^{\varsigma}\|u\|_{L^{p,q}(U\times(-t,\Lambda-t))}\leq 1\,\mbox{ and }\, B_{R_0}(x)\times(-c_{18}R_0^2,0)\subset
\subset U\times(-t,\Lambda-t)
\label{a1is1}
\end{equation}
for all $x\in B_r$. Then by Corollary \ref{clearing} with \eqref{a1is1}, for any $x\in A_1\cap B_r$ and $0<R\leq R_0$, we have
\begin{equation}
\|V_{-c_{18}R^2}\|(B_{14R/15}(x))>c_{19}R^k,
\label{a1is2}
\end{equation}
since we would have $\|V_0\|(B_{4R/5}(x))=0$ otherwise, contradicting $x\in {\rm spt}\, \|V_0\|$. 
Next let
\begin{equation}
A_{1,m}:=\{x\in A_1\cap B_r\,\,:\,\, \|V_0\|(B_R (x))\leq c_{19}R^k/2,\,\, 0<\forall R\leq r/m\}.
\label{a1is3}
\end{equation}
Since $\Theta^k(\|V_0\|,x)=0$ for ${\mathcal H}^k$ a.e$.$ $x\in A_1\cap B_r$, we have
${\mathcal H}^k(A_1\cap B_r)={\mathcal H}^k (\cup_{m=1}^{\infty}A_{1,m})$. Since $A_{1,m}$
is an increasing sequence and since ${\mathcal H}^k(A_1\cap B_r)>0$, there exists some
$m_0$ such that ${\mathcal H}^k(A_{1,m_0})>0$. By the definition of the Hausdorff measure,
there exists some $\delta_0>0$ such that
\begin{equation}
{\bf b}:={\mathcal H}^k_{\delta_0}(A_{1,m_0}):=\inf\left\{\sum_{j=1}^{\infty}\omega_k\left(\frac{{\rm diam}\, G_j}{2}\right)^k
\,:\, A_{1,m_0}\subset\cup_{j=1}^{\infty}G_j,\, {\rm diam}\, G_j<2\delta_0\right\}>0,
\label{a1is4}
\end{equation}
where $\inf$ is taken over all such covering of $A_{1,m_0}$. Note that ${\bf b}<\infty$ even if ${\mathcal H}^k
(A_{1,m_0})=\infty$. In the following we fix such $R_0$, $m_0$, $\delta_0$,
$A_{1,m_0}$ and ${\bf b}$.

For any $R$ with $0<R<\min\{R_0,r/m_0,\delta_0\}$, consider the covering $\{\overline{B_R(b)}\}_{b\in A_{1,m_0}}$ of
$A_{1,m_0}$. 
By the Besicovitch covering theorem, there exists a family of collections ${\mathcal B}_1,\cdots,{\mathcal B}_{{\bf B}(n)}$
each of which consists of mutually disjoint balls and which satisfies $A_{1,m_0}\subset \cup_{j=1}^{{\bf B}(n)}
\cup_{\overline{B_R(b)}\in {\mathcal B}_j}\overline{B_R(b)}$. Due to \eqref{a1is4}
and $R<\delta_0$, we have for this covering
\begin{equation}
{\bf b}\leq \omega_k\sum_{j=1}^{{\bf B}(n)}(\mbox{$\#$ of elements of }{\mathcal B}_j )R^k
\label{a1is5}
\end{equation}
By \eqref{a1is5}, there exists at least one of the collections, say, ${\mathcal B}_{j_0}$, consisting of
closures of mutually disjoint $B_R(b_1),\,\cdots, B_R(b_N)$, where $N$ is the number of elements of ${\mathcal B}_{j_0}$, and such that
\begin{equation}
{\bf b}\leq \omega_k {\bf B}(n)N R^k
\label{a1is6}
\end{equation}
holds. Since $R<r/m_0$ and $b_j\in A_{1,m_0}$, \eqref{a1is3} shows
\begin{equation}
\sum_{j=1}^N\|V_0\|(B_R(b_j))=\|V_0\|(\cup_{j=1}^N B_R(b_j))\leq Nc_{19}R^k/2.
\label{a1is7}
\end{equation}
On the other hand, by \eqref{a1is2} and $R<R_0$, we have
\begin{equation}
\sum_{j=1}^N\|V_{-c_{18}R^2}\|(B_{14R/15}(b_j))=\|V_{-c_{18}R^2}\| (\cup_{j=1}^N
B_{14R/15}(b_j))>Nc_{19}R^k.
\label{a1is8}
\end{equation} 
Hence \eqref{a1is6}-\eqref{a1is8} show
\begin{equation}
\|V_0\|(\cup_{j=1}^N B_R(b_j))-\|V_{-c_{18}R^2}\|(\cup_{j=1}^N B_{14R/15}(b_j))<
-\frac{c_{19}{\bf b}}{2\omega_k {\bf B}(n)},
\label{a1is9}
\end{equation}
where we note that the right-hand side is a negative constant independent of $R$. 
Let $\phi\in C_c^{\infty}(B_1)$ be a radially symmetric function such that
$0\leq \phi\leq 1$, $\phi=1$ on $B_{14/15}$ and $|\nabla\phi |\leq 30$. Then define
$\phi_0(x):=\phi(x/(2r))$ and for $j=1,\cdots, N$, $\phi_j(x):=\phi((x-b_j)/R)$. Since
$\cup_{j=1}^N B_{R}(b_j)\subset B_{3r/2}$ and $\{B_R(b_j)\}_{j=1}^N$ are mutually disjoint, we have
\begin{equation}
1\geq \tilde{\phi}:=\phi_0-\sum_{j=1}^N\phi_j\geq 0.
\label{a1is10}
\end{equation}
Due to the definition of $\phi_j$ and by \eqref{a1is9}, we also have
\begin{equation}
\|V_0\|\big(\sum_{j=1}^N\phi_j\big)-\|V_{-c_{18}R^2}\|\big(\sum_{j=1}^N\phi_j\big)<-\frac{c_{19}{\bf b}}{2\omega_k
{\bf B}(n)}.
\label{a1is10.5}
\end{equation}
We have by \eqref{a1is10} and
\eqref{a1is10.5}
\begin{equation}
\|V_0\|(\phi_0)-\|V_{-c_{18}R^2}\|(\phi_0)<\|V_0\|(\tilde{\phi})-\|V_{-c_{18}R^2}\|(\tilde{\phi})
-\frac{c_{19}{\bf b}}{2\omega_k
{\bf B}(n)}.
\label{a1is11}
\end{equation}
Since $\tilde{\phi}\in C^{\infty}_c(U;{\mathbb R}^+)$, by \eqref{maineq} and writing $-c_{18}R^2=: s$ (and omitting
$d\|V_t\|dt$),
\begin{equation}
\begin{split}
\|V_0\|(\tilde{\phi})-\|V_{s}\|(\tilde{\phi})& \leq \int_{s}^0\int (-\tilde{\phi}|h|^2+|\nabla\tilde{\phi}||h|
+\tilde{\phi}|h||u|+|\nabla\tilde{\phi}||u|)\\
& \leq \left\{\left(\int_s^0\int |\nabla\tilde{\phi}|^2\right)^{1/2}+\left(\int_s^0\int_{B_{2r}} |u|^2\right)^{1/2}\right\}\left(\int_s^0\int_{B_{2r}}|h|^2\right)^{1/2} \\
& +\left(\int_s^0\int |\nabla\tilde{\phi}|^2\right)^{1/2}\left(\int_s^0\int_{B_{2r}}|u|^2\right)^{1/2}.
\end{split}
\label{a1is12}
\end{equation}
Since $|\nabla\tilde{\phi}|\leq 30/R$ and by (A3), we have
\begin{equation}
\int_{-c_{18}R^2}^0\int|\nabla\tilde{\phi}|^2\leq 900\omega_k c_{18}(2r)^k E_1,\,\,\,\, \lim_{R\rightarrow 0}
\int_{-c_{18}R^2}^0\int_{B_{2r}}|u|^2=0.
\label{a1is13}
\end{equation}
Combining \eqref{a1is11}-\eqref{a1is13}, for all small $R>0$, we obtain
\begin{equation}
\|V_0\|(\phi_0)-\|V_{-c_{18}R^2}\|(\phi_0)\leq (900\omega_k c_{18}(2r)^k E_1+1)^{1/2}\left(\int_{-c_{18}R^2}^0\int_{B_{2r}}
|h|^2\right)^{1/2}-\frac{c_{19}{\bf b}}{4\omega_k
{\bf B}(n)}.
\label{a1is14}
\end{equation}
Since $\|V_t\|(\phi_0)$ is continuous at $t=0$, for all small $R>0$, we finally obtain
\begin{equation}
c_{20}:=\left(\frac{c_{19}{\bf b}}{8\omega_k {\bf B}(n)}\right)^2(900\omega_k c_{18}(2r)^k E_1+1)^{-1}\leq  \int_{-c_{18}R^2}^0
\int_{B_{2r}}|h|^2.
\label{a1is15}
\end{equation}
Note that $c_{20}$ is independent of $R$. Now consider $\hat{\phi}(x)=\phi(x/3r)$ so that $\hat{\phi}\in C^{\infty}_c(
B_{3r})$ and $\hat{\phi}=1$ on $B_{2r}$. We have by \eqref{maineq}, \eqref{hhhp} and \eqref{a1is15}
\begin{equation}
\|V_0\|(\hat{\phi})-\|V_{-R^2}\|(\hat{\phi})\leq -\int_{-R^2}^0\int \frac{\hat{\phi}|h|^2}{2}+o(1)\leq -\frac{c_{20}}{2}
+o(1)
\label{a1is16}
\end{equation}
for all sufficiently small $R>0$. 
But this contradicts the continuity of $\|V_t\|(\hat{\phi})$ as $t\nearrow 0$. Hence we complete the proof.
\hfill{$\Box$}
\begin{lemma}
We have ${\mathcal H}^k(A_2)=0$.
\label{a2is0}
\end{lemma}
{\it Proof}. Just as in the beginning of proof for Lemma \ref{a1is0}, for a contradiction, assume $B_{3r}\subset
\subset U$ and ${\mathcal H}^k(A_2\cap B_r)>0$. 
Since $A_2\subset M$, it is a finite value. 
Corresponding to $\nu=1/2$, we fix $L$, $\Lambda_4$ and $\varepsilon_7$ by Proposition \ref{prprop}.
For each $m\in{\mathbb N}$ with $(L+1)\leq m$, define
\begin{equation}
\begin{split}
A_{2,m}:= \Big\{& x\in A_2\cap B_r\,\,: \,\, |R^{-k}\|V_0\|(\phi_R^2(\cdot-x))-{\bf c}|\leq {\bf c}/4,\,\, 0<\forall R\leq r/m \,\,\mbox{ and}\\
&  R^{-k-2}\int_{B_{(L+1)R}(x)} |({\rm Tan}_x M)^{\perp}(y-x)|^2\, d\|V_0\|(y)\leq
\varepsilon_7^2 /2,\,\, 0<\forall R\leq r/m \Big\},
\end{split}
\label{a2is1}
\end{equation}
where $\phi_R(x)$ and ${\bf c}$ are as in Proposition \ref{prprop}. Since there exists an approximate tangent space ${\rm Tan}_x M$
for ${\mathcal H}^k$ a.e$.$ $x\in A_2\cap B_r$, we have ${\mathcal H}^k(A_2\cap B_r)={\mathcal H}^k(\cup_{(L+1)\leq m\in
{\mathbb N}}
A_{2,m})$. Since ${\mathcal H}^k (A_2\cap B_r)>0$, for some $(L+1)\leq m_0\in {\mathbb N}$, we have
${\mathcal H}^k(A_{2,m_0})>0$. Similar to \eqref{a1is4}, we may choose $\delta_0>0$ such that
${\mathcal H}^k_{\delta_0}(A_{2,m_0})>0$. 
Choose small $R_0>0$ so that 
\begin{equation}
R_0^{\varsigma}\|u\|_{L^{p,q}(U\times(-t,\Lambda-t))}\leq \varepsilon_7,\,\,\,\mbox{and}\,\,\,
B_{(L+1)R_0}(x)\times(-\Lambda_4 R_0^2,\Lambda_4 R_0^2)\subset\subset U\times(-t,\Lambda-t)
\label{a2is1.5}
\end{equation}
for all $x\in B_r$. For each $0<R\leq R_0$, we define
\begin{equation}
\begin{split}
&A_{2,m_0}^1(R):=\{x\in A_{2,m_0}\,:\,\eqref{pa1}\mbox{ fails for any $T\in {\bf G}(n,k)$ }\},\\
&A_{2,m_0}^2(R):=\{x\in A_{2,m_0}\,:\, \eqref{pa3}\mbox{ fails}\}, \,\,\, A_{2,m_0}^3(R):=\{x\in A_{2,m_0}\,:\, \eqref{pa4}\mbox{ fails}\}
\end{split}
\label{a2is3}
\end{equation}
with $\nu=1/2$ in \eqref{pa3} and \eqref{pa4}. Since any point of $A_{2,m_0}$ is not $C^{1,\varsigma}$ regular, and by \eqref{a2is1.5}, Proposition \ref{prprop} shows $A_{2,m_0}\subset \cup_{l=1}^3 A_{2,m_0}^l(R)$ for all $0<R\leq R_0$ (if not,
then $x\in A_{2,m_0}\setminus \cup_{l=1}^3 A_{2,m_0}^l(R)$ satisfies \eqref{pa1}-\eqref{pa4} and Proposition \ref{prprop} applies to $(x,0)$). 
Since ${\mathcal H}^k_{\delta_0}$ is sub-additive, we have ${\mathcal H}^k_{\delta_0}(A_{2,m_0})\leq\sum_{l=1}^3
{\mathcal H}^k_{\delta_0}(A_{2,m_0}^l(R))$ and we have for $0<R\leq R_0$
\begin{equation}
0<{\bf b}:={\mathcal H}^k_{\delta_0}(A_{2,m_0})/3\leq \max_{l=1,2,3} \{{\mathcal H}^k_{\delta_0} (A_{2,m_0}^l(R))\}.
\label{a2is4}
\end{equation}
We prove that we will have a contradiction whichever of three quantities takes the maximum value.
We proceed with the assumption that $0<R\leq \min\{R_0,r/(L+1)m_0,\delta_0/(L+1)\}$. 

{\bf Case 1}. When ${\bf b}\leq {\mathcal H}^k_{\delta_0}(A_{2,m_0}^1(R))$.
\newline
Consider the covering $\{\overline{B_{(L+1)R}(x)}\}_{x\in A_{2,m_0}^1(R)}$ of $A_{2,m_0}^1(R)$. 
Just as in Lemma \ref{a1is0}, there exist mutually disjoint $B_{(L+1)R}(b_1),\cdots,B_{(L+1)R}(b_N)$ with
$b_1,\cdots,b_N\in A_{2,m_0}^1(R)$ and with 
\begin{equation}
{\bf b}\leq \omega_k {\bf B}(n)(L+1)^k N R^k.
\label{a2is5}
\end{equation}
Since $b_j\in A_{2,m_0}$,
\eqref{a2is1} shows
\begin{equation}
\sum_{j=1}^N \int_{B_{(L+1)R}(b_j)}|T_j^{\perp}(x-b_j)|^2\, d\|V_0\| \leq \varepsilon_7 ^2 N R^{k+2}/2,
\label{a2is6}
\end{equation}
where we denote $T_j:={\rm Tan}_{b_j} M$. On the other hand, since $b_j\in A_{2,m_0}^1$, by
\eqref{a2is3}, we have
\begin{equation}
\sum_{j=1}^N \int_{B_{LR} (b_j)}|T_j^{\perp}(x-b_j)|^2\, d\|V_{-\Lambda_4 R^2}\|\geq \varepsilon_7^2 NR^{k+2}.
\label{a2is7}
\end{equation}
Let $\phi\in C_c^{\infty}(B_{L+1})$ be a radially symmetric function such that $0\leq \phi\leq 1$,
$\phi=1$ on $B_{L}$ and $|\nabla\phi|\leq 2$. For $j=1,\cdots, N$, define 
\begin{equation}
\phi_j(x)=
\phi((x-b_j)/R)\,\,\,\mbox{ and }\,\,\,\xi_j(x)=\{(L+1)R\}^{-2}\phi_j(x)|T_j^{\perp}(x-b_j)|^2.
\label{a2is8}
\end{equation}
Note that $\{\xi_j\}_{j=1}^N$ have mutually disjoint supports and that $0\leq \xi_j\leq 1$.  
Now combining \eqref{a2is5}-\eqref{a2is8}, we then obtain
\begin{equation}
 \|V_0\|\big(\sum_{j=1}^N\xi_j\big)-\|V_{-\Lambda_4 R^2}\|\big(\sum_{j=1}^N \xi_j\big)\leq -\frac{\varepsilon_7^2 {\bf b}}{2\omega_k
 {\bf B}(n)(L+1)^{k+2}}.
 \label{a2is9}
 \end{equation}
 Let $\phi_0\in C_c^{\infty}(B_{2r})$ be a radially symmetric function with $0\leq \phi_0\leq 1$, 
 $\phi_0=1$ on $B_{3r/2}$ and $|\nabla\phi_0|\leq 3/r$. Since ${\rm spt}\, \xi_j\subset
 B_{3r/2}$, we have
 \begin{equation}
 1\geq \tilde{\phi}:=\phi_0-\sum_{j=1}^N\xi_j\geq 0, 
 \label{a2is10}
 \end{equation}
 and by \eqref{a2is9},
 \begin{equation}
 \|V_0\|(\phi_0)-\|V_{-\Lambda_4 R^2}\|(\phi_0)\leq \|V_0\|(\tilde{\phi})
 -\|V_{-\Lambda_4 R^2}\|(\tilde{\phi})-\frac{\varepsilon_7^2 {\bf b}}{2\omega_k
 {\bf B}(n)(L+1)^{k+2}}.
\label{a2is11}
\end{equation}
Now the rest of the argument proceeds just like \eqref{a1is12}-\eqref{a1is16}
since newly defined $\tilde{\phi}$ also satisfies $|\nabla\tilde{\phi}|\leq 3/R$, and
we similarly obtain a contradiction to the continuity at $t=0$. Thus Case 1 cannot 
occur for sufficiently small $R>0$. 

{\bf Case 2}. When ${\bf b}\leq {\mathcal H}^k_{\delta_0}(A_{2,m_0}^2(R))$.
\newline
Consider the covering $\{\overline{B_R(x)}\}_{x\in A_{2,m_0}^2(R)}$ of $A_{2,m_0}^2(R)$.
Just as in Lemma \ref{a1is0}, there exist mutually disjoint $B_R(b_1),\cdots,B_R(b_N)$
with $b_1,\cdots,b_N\in A_{2,m_0}^2(R)$ with
\begin{equation}
{\bf b}\leq \omega_k {\bf B}(n)NR^k.
\label{a3is0}
\end{equation}
Since $b_j\in A_{2,m_0}$, \eqref{a2is1} shows
\begin{equation}
R^{-k}\|V_0\|(\phi_R^2(\cdot-b_j))\leq 5{\bf c}/4,
\label{a3is1}
\end{equation}
while $b_j\in A_{2,m_0}^2(R)$ implies
\begin{equation}
R^{-k}\|V_{-(\Lambda_4-3)R^2}\|(\phi_R^2(\cdot-b_j))\geq 3{\bf c}/2.
\label{a3is2}
\end{equation}
Hence \eqref{a3is0}-\eqref{a3is2} show 
\begin{equation}
\|V_0\|\big(\sum_{j=1}^N \phi_R^2(\cdot-b_j)\big)-\|V_{-(\Lambda_4-3)R^2}\|\big(\sum_{j=1}^N 
\phi_R^2(\cdot-b_j)\big)\leq -\frac{{\bf c}{\bf b}}{4\omega_k {\bf B}(n)}.
\label{a3is3}
\end{equation}
Define $\phi_0\in C^{\infty}_c(B_{2r})$ as in Case 1 and
\begin{equation}
1\geq \tilde{\phi}:=\phi_0-\sum_{j=1}^N\phi_R^2(\cdot-b_j)\geq 0.
\label{a3is4}
\end{equation}
Then by \eqref{a3is3} and \eqref{a3is4}, we have
\begin{equation}
\|V_0\|(\phi_0)-\|V_{-(\Lambda_4-3)R^2}\|(\phi_0)\leq \|V_0\|(\tilde{\phi})-\|V_{-(\Lambda_4
-3)R^2}\|(\tilde{\phi})-\frac{{\bf c}{\bf b}}{4\omega_k {\bf B}(n)}.
\label{a3is5}
\end{equation}
Note that $|\nabla \tilde{\phi}|\leq c/R$, and the rest of the argument proceeds just like 
\eqref{a1is12}-\eqref{a1is16}, leading to a contradiction for all sufficiently small $R>0$. 

{\bf Case 3}. When ${\bf b}\leq {\mathcal H}_{\delta_0}^k (A^3_{2,m_0}(R))$.
\newline
The argument is similar to Case 2. With the same kind of covering balls, we have 
\begin{equation}
R^{-k}\|V_0\|(\phi_R^2(\cdot-b_j))\geq 3{\bf c}/4
\label{a4is0}
\end{equation}
in place of \eqref{a3is1} and
\begin{equation}
R^{-k}\|V_{(\Lambda_4-3)R^2}\|(\phi_R^2(\cdot-b_j)\leq {\bf c}/2
\label{a4is1}
\end{equation}
in place of \eqref{a3is2}, leading to 
\begin{equation}
\|V_{(\Lambda_4-3)R^2}\|\big(\sum_{j=1}^N\phi_R^2(\cdot-b_j)\big)
-\|V_0\|\big(\sum_{j=1}^N\phi_R^2(\cdot-b_j)\big)\leq -\frac{{\bf c}{\bf b}}{4\omega_k {\bf B}(n)}.
\label{a4is2}
\end{equation}
Similar argument leads to a contradiction to the continuity of $\|V_t\|(\tilde{\phi})$ 
as $t\searrow 0$. 

Thus all three alternatives lead to a contradiction and this shows the claim of Lemma 
\ref{a2is0}.
\hfill{$\Box$}

By Lemma \ref{a1is0} and Lemma \ref{a2is0}, now we conclude the proof of Theorem \ref{gp}.

\section{Concluding remarks}
\subsection{$C^{\infty}$ regularity}
In case $u=0$, Brakke claimed that the support of moving varifold is a.e$.$ time and a.e$.$ everywhere $C^{\infty}$. 
Unfortunately, in examining his proof, there is an essential gap in the argument for going from $C^{1,\varsigma}$ to $C^{2}$ estimates. After showing that the spacial gradient is H\"{o}lder continuous in the middle
of page 202 of \cite{Brakke}, Brakke proceeds to prove that the graph can be approximated by some quadratic function, claiming
that the graph is $C^{2,1/8}$ in the space variables. For the proof, from (21) of page 203 to (22) of the next page, 
he substitutes $p+q=R^{\varepsilon}$ with $\varepsilon=1/100$, which one can easily see from (22) that it should be $p+q=R^{3/2}$ instead. 
In the 4th line of (25), Brakke substitutes $p+q=R^{\varepsilon}$, which instead can only be $R^{3/2}$ from the 
previous computations. This is a crucial step to
show the decay rate is good enough. If $p+q=R^{3/2}$, no improvement is achieved, and his argument
has an essential gap in this regard. To remedy the situation, the forthcoming paper \cite{regTonegawa}  gives a proof that, assuming 
that $u$ is $\alpha$-H\"{o}lder continuous, the support of moving varifold is a.e$.$ $C^{2,\alpha}$ (in parabolic sense)
and it satisfies the motion law classically. 
Thus by the standard parabolic regularity theory, if $u$ is $C^{\infty}$, and in particular $u=0$, the support is a.e$.$ $C^{\infty}$, which
is Brakke's original claim. The idea of proof is to carry out another blow up 
argument with second order approximation and is similar to the proof in Section 8 in spirit. 
\label{conrem1}
\subsection{Time independent case}
Suppose that $V_t$ and $u$ do not depend on $t$ while (A1)-(A4) are satisfied. Then \eqref{maineq}
reduces simply to 
\begin{equation}
0\leq {\mathcal B}(V,u,\phi),\,\,\, \forall \phi\in C^1_c(U;{\mathbb R}^+).
\label{steq}
\end{equation}
We have $V=|M|$ for some countably $k$-rectifiable set $M\subset U$ by (A1), $h(V)\in L^2(\|V\|)$ 
by \eqref{hhh} and $u\in L^p(\|V\|)$ with $p>k$ by \eqref{expcond}. We claim that
\begin{lemma}
$h(V)+u^{\perp}=0$, $\|V\|$ a.e$.$ on $U$.
\end{lemma}
{\it Proof}.
Let $\hat{x}\in M$ be any Lebesgue point of 
$h(V)$ and $u^{\perp}$ with respect to ${\mathcal H}^k$ measure. We may also assume
that the approximate tangent space for $M$ exists at $\hat{x}$ and $h(V)=h(V)^{\perp}$. Set of such
points is a full measure set on $M$. For any $\phi\in C^1_c({\mathbb R}^n;{\mathbb R}^+)$ and
$0<r<1$, define $\phi_r(y):=\phi((y-\hat{x})/r)$. For all sufficiently small $r>0$, $\phi_r\in C_c^1(U;
{\mathbb R}^+)$. By \eqref{steq} with this $\phi_r$, we have
\begin{equation}
0\leq \int \{-h\phi_r+r^{-1}\nabla\phi((\cdot-\hat{x})/r)\}\cdot(h+u^{\perp})\, d\|V\|.
\label{steq2}
\end{equation}
By change of variables $\tilde{y}=(y-\hat{x})/r$, dividing \eqref{steq2} by $r^{k-1}$, and
letting $r\rightarrow 0$, we obtain
\begin{equation}
0\leq \int_{{\rm Tan}_{\hat{x}} M}\nabla\phi(\tilde{y})\, d{\mathcal H}^k(\tilde{y})\cdot (h(V,\hat{x})+u^{\perp}(\hat{x})).
\label{steq3}
\end{equation}
Note that the vector $\int_{{\rm Tan}_{\hat{x}}M}\nabla\phi(\tilde{y})\, d{\mathcal H}^k(\tilde{y})$ is
perpendicular to ${\rm Tan}_{\hat{x}}M$, but otherwise it may be arbitrarily chosen depending on $\phi$.
Thus for \eqref{steq3} to be true, $h(V,\hat{x})+u^{\perp}(\hat{x})=0$. Thus this
holds for $\|V\|$ a.e$.$ $\hat{x}$ on U. 
\hfill{$\Box$}

Thus, in this case, we have $h(V,\cdot)\in L^p(\|V\|)$ with $p>k$. Since $V$ has 
density equal to 1 a.e$.$, Allard's regularity theorem applies. In this sense, 
the partial regularity theorem of the present paper is a natural generalization
of Allard's regularity theorem. 

\section{Appendix}
\subsection{Inequalities for ${\bf G}(n,k)$}
In this subsection we collect some well-known facts which are used in this paper.
For completeness we include their proofs.
\begin{lemma}
For $S,\, T\in {\bf G}(n,k)$ and $v\in {\mathbb R}^n$, we have the following.
\begin{equation}
I\cdot T=k,\hspace{.5cm}T^t=T,\hspace{.5cm} T\circ T=T, \hspace{.5cm} T\circ T^{\perp}=0.
\label{simin0}
\end{equation}
\begin{equation}
0\leq k-S\cdot T=S^{\perp}\cdot T\leq k\|S-T\|^2.
\label{simin1}
\end{equation}
\begin{equation}
0\leq \|S-T\|^2 \leq (S-T)\cdot (S-T)=2T^{\perp}\cdot S.
\label{simin2}
\end{equation}
\begin{equation}
|T(S^{\perp}(v))|\leq \|T-S\||v|.
\label{simin3}
\end{equation}
\begin{equation}
|T(S^{\perp}(T(v)))|\leq \|T-S\|^2|v|.
\label{simin4}
\end{equation}
\label{re-ea}
\end{lemma}
{\it Proof}. By taking a set of orthonormal vectors $\{v_1,\cdots,v_n\}$
such that $v_1,\cdots,v_k\in T$, $I\cdot T={\rm tr}\,(T)=\sum_{i=1}^k v_i\cdot T(v_i)=k$.
It is clear that $T^t=T$ if represented for such basis, and thus it is symmetric matrix.
$T\circ T=T$ and $T\circ T^{\perp}=0$ are clear and this ends the proof of \eqref{simin0} \newline
By \eqref{simin0}, and for the orthonormal vectors $\{v_1,\cdots,v_n\}$ as above,
\begin{equation*}
\begin{split}
&k-S\cdot T=I\cdot T-S\cdot T=S^{\perp}\cdot T
={\rm tr}\, (S^{\perp}\circ S^{\perp}\circ T\circ T) ={\rm tr}\, (S^{\perp}\circ T\circ T\circ S^{\perp}) \\ &= {\rm tr}\, ((T-S)\circ T\circ
T\circ (T-S)) 
=\sum_{i=1}^n v_i^t\cdot (T\circ(T-S)\circ(T-S)\circ T)(v_i) \\
&=\sum_{i=1}^k v_i^t\cdot ((T-S)\circ(T-S))(v_i)\leq k\|S-T\|^2.
\end{split}
\end{equation*}
It is also clear from the intermediate expression that the above is nonnegative. Thus
this proves \eqref{simin1}.
\newline
By choosing a unit vector $v$ such that $\|S-T\|=|(S-T)(v)|$, and noticing 
\begin{equation*}
|(S-T)(v)|^2=v^t\cdot((S-T)\circ(S-T))(v)\leq (S-T)\cdot(S-T),
\end{equation*}
we have
$\|S-T\|^2\leq (S-T)\cdot(S-T)$. Also
\begin{equation*}
(S-T)\cdot(S-T)=2k-2T\cdot S=2(I\cdot S-T\cdot S)=2T^{\perp}\cdot S
\end{equation*}
proves \eqref{simin2}.
\newline
For \eqref{simin3} and \eqref{simin4}, we have
\begin{equation*}
\begin{split}
& |T(S^{\perp}(v))|=|T((I-S)(v))|=|T((T-S)(v))|\leq \|T-S\||v|,\\
& |T(S^{\perp}(T(v)))|=|T((I-S)(T-S)(v))|=|T((T-S)^2(v))|\leq \|T-S\|^2 |v|.
\end{split}
\end{equation*}
\hfill{$\Box$}
\subsection{Additional results}
In this subsection we include
some results from \cite{Brakke} along with their proofs
for the convenience of readers. 
\begin{lemma}(Tilt of tangent planes \cite[5.5]{Brakke},\cite[8.13]{Allard}) If $U\subset {\mathbb R}^n$ is open,  
$V\in {\bf IV}_k(U)$, $T\in {\bf G}(n,k)$, $\phi\in C^1_c(U;{\mathbb R}^+)$,
\begin{equation*}
\begin{split}
&\alpha=\left(\int_U |h(V,x)|^2 \phi^2(x)\, d\|V\|(x)\right)^{\frac12},\hspace{.5cm}
\mu=\left(\int_U |T^{\perp}(x)|^2\phi^2(x)\, d\|V\|(x)\right)^{\frac12},\\
&\tilde{\mu}=\left(\int_U |T^{\perp}(x)|^2|\nabla\phi(x)|^2\, d\|V\|(x)\right)^{\frac12},\hspace{.5cm}
\beta=\left(\int_{G_k(U)}\|S-T\|^2 \phi^2(x)\, dV(x,S)\right)^{\frac12},
\end{split}
\end{equation*}
then we have
\begin{equation}
\beta^2\leq 4\alpha\mu+16\tilde{\mu}^2.
\label{append1}
\end{equation}
\label{tiltexlem}
\end{lemma}
{\it Proof}.
Let $g(x)=\phi^2(x)T^{\perp}(x)$ for $x\in U$. Then for $S\in {\bf G}(n,k)$ we have
\begin{equation}
\nabla g(x)\cdot S=2\phi(x)S(T^{\perp}(x))\cdot \nabla\phi(x)+\phi^2(x)T^{\perp}\cdot S.
\label{append2}
\end{equation}
By  \eqref{simin2} and \eqref{simin0}, \eqref{append2} implies
\begin{equation}
\begin{split}
&\frac12 \phi^2(x)\|S-T\|^2\leq \nabla g(x)\cdot S+2\phi(x)|S(T^{\perp}(x))\cdot \nabla\phi(x)| \\
&\leq \nabla g(x)\cdot S+2\phi(x)\|S-T\||T^{\perp}(x)||\nabla\phi(x)| \\
&\leq \nabla g(x)\cdot S+\frac14 \phi^2(x)\|S-T\|^2+4|T^{\perp}(x)|^2|\nabla \phi(x)|^2.
\end{split}
\label{append3}
\end{equation}
Since 
\begin{equation}
\int_{G_k(U)}\nabla g(x)\cdot S\, dV(x,S)=-\int_U\phi^2(x) h(V,x)\cdot T^{\perp}(x)\, d\|V\|(x)
\leq \alpha\mu
\label{append4}
\end{equation}
by Schwarz' inequality, \eqref{append3} and \eqref{append4} give
\eqref{append1}.
\hfill{$\Box$}
\begin{thm}(Cylindrical growth rates \cite[6.4]{Brakke})
Suppose
$T\in {\bf G}(n,k)$, $0<R_1<R_2<\infty$, $0\leq \alpha<\infty$, $0\leq \beta<\infty$,
$V\in {\bf IV}_k(C(T,R_2))$ is finite and ${\rm spt}\, \|V\|$ is bounded. Suppose $\phi\in
C^3(C(T,1);{\mathbb R}^+)$ depends only on $|T(x)|$ and ${\rm spt}\, \phi
\subset C(T,1)$. Moreover suppose
\begin{equation}
\int_{C(T,R_2)}|h(V,x)|^2 \phi(x/r)\, d\|V\|(x)\leq \alpha^2 r^k,\hspace{.5cm} R_1\leq
\forall r\leq R_2,
\label{append5}
\end{equation}
\begin{equation}
\int_{G_k(C(T,R_2))} \|S-T\|^2 \phi(x/r)\, dV(x,S)\leq \beta^2 r^k,\hspace{.5cm}
R_1\leq \forall r\leq R_2.
\label{append6}
\end{equation}
Then we have
\begin{equation}
\left|R_2^{-k}\|V\|(\phi(x/R_2))-R_1^{-k}\|V\|(\phi(x/R_1))\right|\leq 
k\beta^2 \log(R_2/R_1)+\alpha\beta(R_2-R_1)+\beta^2.
\label{append7}
\end{equation}
\label{cylgrolem}
\end{thm}
{\it Proof}. For each fixed $r\in [R_1,R_2]$ we use the vector field $g(x)=r^{-1}\phi(x/r)T(x)$
in the first variation. Since ${\rm spt}\, \|V\|$ is bounded in $C(T,R_2)$, $g$ can be modified
suitably so that it has a compact support in $C(T,R_2)$ and so that it does not affect the computations.
We have
\begin{equation}
\delta V(g)=r^{-1}\int_{G_k(C(T,R_2))}\phi(x/r)T\cdot S+T(x)\otimes \nabla\phi(x/r)\cdot S\,
dV(x,S).
\label{append8}
\end{equation}
Since $\phi$ depends only on $|T(x)|$, we can derive that
\begin{equation}
\nabla \phi(x/r)=-r\frac{\partial \phi(x/r)}{\partial r}\frac{T(x)}{|T(x)|^2}.
\label{append9}
\end{equation}
Using the perpendicularity of mean curvature \eqref{fvf2}, Schwarz' inequality, \eqref{simin3}, \eqref{append5} and 
\eqref{append6}, we have
\begin{equation}
|\delta V(g)|=\left|\int_{G_k(C(T,R_2))}h(V,x)\cdot S^{\perp}(T(x)/r)\phi(x/r)\, dV(x.S)\right|
\leq \alpha\beta r^k.
\label{append10}
\end{equation}
We use \eqref{append8}-\eqref{append10} and  \eqref{apsup1} of Lemma \ref{linalg} to derive
\begin{equation}
\begin{split}
&\Big| \frac{d}{dr}\int_{G_k(C(T,R_2))}|S(T(x))|^2|T(x)|^{-2}\phi(x/r)\, dV(x,S) \\
&-\frac{k}{r}\int_{G_k(C(T,R_2))}|S(T(x))|^2|T(x)|^{-2}\phi(x/r)\, dV(x,S)\Big|
\leq r^{-1}k\beta^2 r^k+\alpha\beta r^k.
\end{split}
\label{append11}
\end{equation}
Dividing both sides of \eqref{append11} by $r^k$ and integrating this from
$R_1$ to $R_2$, we obtain
\begin{equation}
\begin{split}
\Big| r^{-k}\int_{G_k(C(T,R_2))}|S(T(x))|^2 & |T(x)|^{-2}\phi(x/r)\, dV(x,S)
\big|_{r=R_1}^{R_2}\Big| \\ 
& \leq k\beta^2\log(R_2/R_1)+\alpha\beta(R_2-R_1).
\end{split}
\label{append12}
\end{equation}
By \eqref{apsup1.1} of Lemma \ref{linalg},  \eqref{append6} and \eqref{append12},
we obtain \eqref{append7}.
\hfill{$\Box$}
\begin{lemma}
For $x\in {\mathbb R}^n$, $S,\, T\in {\bf G}(n,k)$, we have
\begin{equation}
|T\cdot S-k|S(T(x))|^2 |T(x)|^{-2}|\leq k\|S-T\|^2,
\label{apsup1}
\end{equation}
\begin{equation}
0\leq 1-|S(T(x))|^2 |T(x)|^{-2}\leq \|S-T\|^2.
\label{apsup1.1}
\end{equation}
\label{linalg}
\end{lemma}
{\it Proof}. We verify \eqref{apsup1} first. We have
\begin{equation}
\begin{split}
|S(T(x))|^2& =T(x)\cdot T(S(T(x)))=T(x)\cdot T( (S-I)(T(x)))+|T(x)|^2\\ 
& =-T(x)\cdot
T(S^{\perp}(T(x)))+|T(x)|^2.
\end{split}
\label{apsup2}
\end{equation}
Thus we obtain from \eqref{apsup2}
\begin{equation}
k|S(T(x))|^2|T(x)|^{-2}-T\cdot S=k-T\cdot S-k \frac{T(x)\cdot T(S^{\perp}(T(x)))}{|T(x)|^2}.
\label{apsup3}
\end{equation}
Writing $v=T(x)$, we have
\begin{equation}
T(x)\cdot T(S^{\perp}(T(x)))=v\cdot S^{\perp}(v)=|S^{\perp}(v)|^2\geq 0.
\label{apsup4}
\end{equation}
Combining \eqref{simin1}, \eqref{simin4}, \eqref{apsup3}, and \eqref{apsup4}, we
prove \eqref{apsup1}. Observing \eqref{apsup2} and \eqref{simin4}, we see that
\eqref{apsup1.1} holds. 
\hfill{$\Box$}

\end{document}